\documentclass[final]{article}

\usepackage[legalpaper, margin=1in]{geometry}
\usepackage[latin1]{inputenc}
\usepackage{float}
\usepackage{fancyvrb}
\usepackage{upgreek}
\usepackage{mathdots}
\usepackage{BOONDOX-calo}
\usepackage[mathscr]{eucal}
\usepackage{color}
\usepackage[usenames,dvipsnames,svgnames,table]{xcolor}
\usepackage{amsmath,amssymb,amsopn,mathtools}
\usepackage{amsthm}
\usepackage{graphicx}
\usepackage{hyperref}       
\usepackage{relsize}
\usepackage{subcaption}
\usepackage{cite}

\usepackage{enumitem}

\newlist{Assumption}{enumerate}{1}
\setlist[Assumption]{label=A\arabic*}

\newcommand{\bmat}[1]{\begin{bmatrix}#1\end{bmatrix}} 
\definecolor{Blue}{rgb}{0,0,1}
\definecolor{Red}{rgb}{1,0,0}
\definecolor{Green}{rgb}{0,1,0}

\newcommand{\paramSymbol}{\mu}
\newcommand{\param}{\boldsymbol{\paramSymbol}}
\newcommand{\paramk}[1]{\param_{#1}}
\newcommand{\nparam}{n_{\paramSymbol}}
\newcommand{\statemat}{\boldsymbol{A}}

\newcommand{\ststatemat}{\statemat^{\text{st}}}

\newcommand{\basis}{\boldsymbol{\Phi}}

\newcommand{\spatialbasis}{\basis_{\text{s}}}

\newcommand{\spacetimebasis}{\basis_{\text{st}}}
\newcommand{\basisvec}{\boldsymbol{\phi}}
\newcommand{\spacetimebasisvec}{\basisvec^{\text{st}}}
\newcommand{\spacetimebasisveck}[1]{\spacetimebasisvec_{#1}}
\newcommand{\spatialbasisvec}{\basisvec^{\text{s}}}
\newcommand{\spatialbasisveck}[1]{\spatialbasisvec_{#1}}
\newcommand{\temporalbasis}{\basis_{\text{t}}}
\newcommand{\temporalbasisi}[1]{\temporalbasis^{#1}}
\newcommand{\temporalbasisvec}{\basisvec^{\text{t}}}
\newcommand{\temporalbasisveck}[2]{\temporalbasisvec_{#1#2}}

\newcommand{\identity}{\boldsymbol{I}}

\newcommand{\zero}{\boldsymbol{0}}
\newcommand{\inputmat}{\boldsymbol{B}}

\newcommand{\inputvec}{\boldsymbol{f}}
\newcommand{\stinputvec}{\inputvec^{\text{st}}}

\newcommand{\resSymbol}{r}
\newcommand{\res}{\boldsymbol{\resSymbol}}
\newcommand{\resk}[1]{\res^{(#1)}}
\newcommand{\stres}{\res^{\text{st}}}
\newcommand{\approxres}{\tilde{\res}}
\newcommand{\approxresk}[1]{\approxres^{(#1)}}
\newcommand{\stateSymbol}{u}
\newcommand{\stateSnapshot}{\boldsymbol{U}}
\newcommand{\temporalSnapshot}{\boldsymbol{\Upsilon}}
\newcommand{\temporalSnapshotk}[1]{\temporalSnapshot_{#1}}
\newcommand{\stateSnapshotk}[1]{\stateSnapshot_{#1}}
\newcommand{\state}{\boldsymbol{\stateSymbol}}
\newcommand{\ststate}{\state^{\text{st}}}
\newcommand{\statek}[1]{\state^{(#1)}}

\newcommand{\initialstate}{\state^0}

\newcommand{\stinitial}{\state_0^{\text{st}}}

\newcommand{\approxstate}{\tilde{\state}}
\newcommand{\approxstatek}[1]{\approxstate^{(#1)}}
\newcommand{\approxststate}{\approxstate^{\text{st}}}

\newcommand{\reducedstate}{\hat{\state}}
\newcommand{\reducedststate}{\reducedstate^{\text{st}}}

\newcommand{\timevar}{t}
\newcommand{\timevark}[1]{\timevar_{(#1)}}
\newcommand{\initialtime}{0}
\newcommand{\finaltime}{T}
\newcommand{\timestep}{\Delta \timevar}
\newcommand{\timestepk}[1]{\timestep^{(#1)}}
\newcommand{\inputveck}[1]{\inputvec^{(#1)}}

\newcommand{\leftsingularmat}{\boldsymbol{W}}
\newcommand{\rightsingularmat}{\boldsymbol{V}}

\newcommand{\leftsingularmattemporal}{\boldsymbol{\Lambda}}
\newcommand{\leftsingularmattemporalk}[1]{\leftsingularmattemporal_{#1}}
\newcommand{\rightsingularmattemporal}{\boldsymbol{\Psi}}
\newcommand{\rightsingularmattemporalk}[1]{\rightsingularmattemporal_{#1}}

\newcommand{\leftsingularvec}{\boldsymbol{w}}
\newcommand{\rightsingularvec}{\boldsymbol{v}}
\newcommand{\leftsingularveck}[1]{\leftsingularvec_{#1}}
\newcommand{\rightsingularveck}[1]{\rightsingularvec_{#1}}
\newcommand{\singularvalmat}{\boldsymbol{\Sigma}}
\newcommand{\singularvalmatk}[1]{\singularvalmat_{#1}}

\newcommand{\singularvaluesSymbol}{\sigma}

\newcommand{\singularvaluesArg}[1]{\singularvaluesSymbol_{#1}}

\newcommand{\NN}{\mathbb{N}}
\newcommand{\RR}[1]{\ensuremath{\mathbb{R}^{ #1 }}}

\newcommand{\paramspace}{\Omega_{\paramSymbol}}

\newcommand{\nspace}{N_s}
\newcommand{\nreducedspace}{n_s}
\newcommand{\nreducedtime}{n_t}
\newcommand{\ntime}{N_t}
\newcommand{\ninput}{N_i}

\newcommand{\natNo}{\NN} 
\newcommand{\nat}[1]{\natNo(#1)}

\newcommand{\reducedtimeindex}{j}
\newcommand{\reducedspaceindex}{i}
\providecommand{\keywords}[1]
{
  \small
  \textbf{\textit{Keywords---}} #1
}

\newcommand{\bds}[1]{{\boldsymbol{#1}}}

\newcommand{\constraintMatROM}{\bds{A}}
\def\mubold{\boldsymbol{\mu}}
\newcommand{\basismatspaceSymb}{\Phi}
\newcommand{\basismatspace}{\boldsymbol{\basismatspaceSymb}}

\newcommand{\jprimej}{(j^\prime,j)}

\newcommand{\Amatrix}{\constraintMatROM(\mubold) }
\newcommand{\Bhatmatrix}{\bds{B}(\mubold) }
\newcommand{\uinit}{\bds{u}_{0}(\mubold) }

\newcommand{\AhatGalerkin}{\hat{\constraintMatROM}^{st,g}(\mubold) }
\newcommand{\AhatPetrovGalerkin}{\hat{\constraintMatROM}^{st, pg}(\mubold) }
\newcommand{\AhatGalerkinBlock}{\hat{\constraintMatROM}^{st, g}_{\jprimej}(\mubold) }
\newcommand{\AhatPetrovGalerkinBlock}{\hat{\constraintMatROM}^{st, pg}_{\jprimej}(\mubold) }
\newcommand{\temporalBasisMatrixk}{\bds{D}_k}
\newcommand{\temporalBasisMatrixkpone}{\bds{D}_{k+1}}
\newcommand{\fhatPetrovGalerkin}{\hat{\bds{f}}^{st, pg}(\mubold)}
\newcommand{\fhatGalerkin}{\hat{\bds{f}}^{st,g}(\mubold)}
\newcommand{\uhatPetrovGalerkin}{\hat{\bds{u}_0}^{st, pg}(\mubold)}
\newcommand{\uhatGalerkin}{\hat{\bds{u}_0}^{st,g}(\mubold)}

\newcommand{\IdentityNs}{\bds{I}_{N_s}}
\newcommand{\Astparameter}{\constraintMatROM^{st} (\mubold) }

\newcommand{\ststatedummy}{\boldsymbol{v}}
\newcommand{\reducedststatedummy}{\hat{\ststatedummy}}
\newcommand{\approxstres}{\approxres^{\text{st}}}
\newcommand{\argmin}[1]{\underset{#1}{\text{argmin}}}

\newcommand{\Aparameter}{\constraintMatROM(\mubold)}
\newcommand{\Bparameter}{\bds{B}(\mubold)}

\newcommand{\errorSymbol}{e}
\newcommand{\error}{\boldsymbol{\errorSymbol}}
\newcommand{\errork}[1]{\error^{(#1)}}
\newcommand{\sterror}{\error^{st}}

\newcommand{\relError}{\text{relative error}}

\numberwithin{equation}{section}

\newtheorem{theorem}{Theorem}[section]

\hypersetup{
	colorlinks=true,   
	citecolor=blue,     
	filecolor=blue,     
	linkcolor=red,    
	urlcolor=blue}      

\title{Efficient space--time reduced order model for linear dynamical systems
       in Python using less than 120 lines of code}

\author{
Youngkyu Kim\thanks{Mechanical Engineering, University of California, Berkeley,
CA 94720 (youngkyu$\_$kim@berkeley.edu)} \footnotemark[3] \and
Karen May Wang\thanks{Design physicist, Lawrence Livermore National Laboratory,
Livermore, CA, USA (wang79@llnl.gov)} \thanks{These authors contributed equally 
to this work.} \and
Youngsoo Choi\thanks{Center for Applied Scientific Computing, Lawrence Livermore
National Laboratory, Livermore, CA 94550 (choi15@llnl.gov)}
}

\date{}

\begin{document}

\maketitle

\begin{abstract}
A classical reduced order model (ROM) for dynamical problems typically involves
  only the spatial reduction of a given problem. Recently, a novel space--time
  ROM for linear dynamical problems has been developed \cite{choi2020space},
  which further reduces the problem size by introducing a temporal reduction in
  addition to a spatial reduction without much loss in accuracy. The authors
  show an order of a thousand speed-up with a relative error of less than
  $10^{-5}$ for a large-scale Boltzmann transport problem. In this work, we
  present for the first time the derivation of the space--time Petrov--Galerkin
  projection for linear dynamical systems and its corresponding block
  structures. Utilizing these block structures, we demonstrate the ease of
  construction of the space--time ROM method with two model problems: 2D
  diffusion and 2D convection diffusion, with and without a linear source term.
  For each problem, we demonstrate the entire process of generating the full
  order model (FOM) data, constructing the space--time ROM, and predicting the
  reduced-order solutions, all in less than 120 lines of Python code. We compare
  our Petrov--Galerkin method with the traditional Galerkin method and show that
  the space--time ROMs can achieve $\mathcal{O}(10^2)$ speed-ups with
  $\mathcal{O}(10^{-3})$ to $\mathcal{O}(10^{-4})$ relative errors for these
  problems. Finally, we present an error analysis for the space--time
  Petrov--Galerkin projection and derive an error bound, which shows an
  improvement compared to traditional spatial Galerkin ROM methods.
\end{abstract}


\pagestyle{myheadings}
\thispagestyle{plain}
\markboth{Y.Kim \& K.Wang \& Y.Choi}{The space--time ROM}

\section{Introduction}\label{sec:intro}
Many computational models for physical simulations are formulated as linear
dynamical systems. Examples of linear dynamical systems include, but are not
limited to, the Schr\"odinger equation that arises in quantum mechanics, the
computational model for the signal propagation and interference in electric
circuits, storm surge prediction models before an advancing hurricane, vibration
analysis in large structures, thermal analysis in various media,
neuro-transmission models in the nervous system, various computational models
for micro-electro-mechanical systems, and various particle transport
simulations. These linear dynamical systems can quickly become large scale and
computationally expensive, which prevents fast generation of solutions. Thus,
areas in design optimization, uncertainty quantification, and controls where
large parameter sweeps need to be done can become intractable, and this
motivates the need for developing a Reduced Order Model (ROM) that can
accelerate the solution process without loss in accuracy.

Many ROM approaches for linear dynamical systems have been developed, and they
can be broadly categorized as data-driven or non data-driven approaches. We give
a brief background of some of the methods here. For the non data-driven
approaches, there are several methods, including: balanced truncation methods
\cite{mullis1976synthesis, moore1981principal, willcox2002balanced,
willcox2005fourier, heinkenschloss2008balanced, sandberg2004balanced,
hartmann2010balanced, petreczky2013balanced, ma2010snapshot}, moment-matching
methods \cite{bai2002krylov, gugercin2008h_2, astolfi2010model,
chiprout1992generalized, pratesi2006generalized}, and Proper Generalized
Decomposition (PGD) \cite{ammar2006new} and its extensions \cite{ammar2007new,
chinesta2010proper, pruliere2010deterministic, chinesta2011overview,
giner2013proper, barbarulo2014proper, amsallem2012stabilization,
amsallem2008interpolation, thomas2003three, hall2000proper}. The balanced
truncation method is by far the most popular method, but it requires the
solution of two Lyapunov equations to construct bases, which is a formidable
task in large-scale problems. Moment matching methods were originally developed
as non data-driven, although later papers extended the method to include it. They provide a
computationally efficient framework using Krylov subspace techniques in an
iterative fashion where only matrix-vector multiplications are required.  The
optimal $H_2$ tangential interpolation for nonparametric systems
\cite{gugercin2008h_2} is also available. Proper Generalized Decomposition was
first developed as a numerical method for solving boundary value problems. It
utilizes techniques to separate space and time for an efficient solution
procedure and is considered a model reduction technique. For the detailed
description of PGD, we refer to a short review paper \cite{chinesta2011short}.
Many data driven ROM approaches have been developed as well. When datasets are
available either from experiments or high-fidelity simulations, these datasets
can contain rich information about the system of interest and utilizing this in
the construction of a ROM can produce an optimal basis.  Although there are some
data-driven moment matching works available \cite{mayo2007framework,
scarciotti2017data}, two popular methods are Dynamic Mode Decomposition (DMD)
and Proper Orthogonal decomposition (POD). DMD generates reduced modes that
embed an intrinsic temporal behavior and was first developed by Peter Schmid
\cite{schmid2010dynamic}. The method has been actively developed and extended to
many applications \cite{chen2012variants, williams2015data,
takeishi2017learning, askham2018variable, schmid2011applications,
kutz2016multiresolution, li2017extended, proctor2016dynamic}.  For a more
detailed description about DMD, we refer to this preprint \cite{tu2013dynamic}
and book \cite{kutz2016dynamic}. POD utilizes the method of snapshots to obtain
an optimal basis of a system and typically applies only to spatial projections,
although temporal projection techniques have been developed as well
\cite{berkooz1993proper, gubisch2017proper, kunisch2001galerkin, hinze2008error,
kerschen2005method, bamer2012application, atwell2001proper, rathinam2003new,
kahlbacher2007galerkin, bonnet1994stochastic, placzek2008hybrid,
legresley2000airfoil, efe2003proper}. 

In our paper, we focus on building a space--time ROM where both spatial and
temporal projections are applied to achieve an optimal reduction. This method
has been developed by previous authors \cite{choi2019space, urban2014improved,
yano2014space1, yano2014space2}, and a space--time ROM for large-scale linear
dynamical systems has been recently introduced \cite{choi2020space}. The authors
show a speed-up of $ >8,000$ with good accuracy for a large-scale transport
problem. In our work, we present several new contributions on the space--time
ROM development:

\begin{itemize}
  \item We derive the block structures of least--squares Petrov--Galerkin
    space--time ROM operators and compare them with the Galerkin space--time ROM
    operators and show that the computational cost saving due to the block
    structure is a factor of the FOM spatial degrees of freedom. 
  \item We present an error analysis of both Galerkin and least--squares
    Petrov--Galerkin space--time ROMs and demonstrate the growth rate of
    the stability constant with the actual space--time operators used in our
    numerical results.
  \item  Utilizing the block structures derived, we demonstrate the ease of
    implementing both Galerkin and least--squares Petrov--Galerkin space--time
    ROM implementations and provide the source code for three canonical
    problems. For each problem, we cover the entire space--time ROM process in
    less than 120 lines of Python code, which includes sweeping a wide parameter
    space and generating data from the full order model, constructing the
    space--time ROM, and generating the ROM prediction in the online phase.
  \item Finally, we present our results for the two model problems and compare
    the speed up and relative error between the Galerkin and Petrov--Galerkin
    methods and show that they give similar results.
\end{itemize}
   
We hope that by providing full access to the Python source codes, researchers can
easily apply space--time ROMs to their linear dynamical problem of interest.
Furthermore, we have curated the source codes to be simple and short so that it
may be easily extended in various multi-query problem settings, such as design
optimization \cite{yoon2010structural, amir2010efficient, amsallem2015design,
gogu2015improving, choi2020gradient, choi2019accelerating, white2020dual},
uncertainty quantification \cite{najm2009uncertainty, walters2002uncertainty,
zang2002needs}, and optimal control problems \cite{petersson2020discrete,
choi2015practical, choi2012simultaneous}.

\subsection{Organization of the paper}
The paper is organized in the following way:
Section~\ref{sec:LinearDynamicalSystems} describes a parametric linear dynamical
systems and space--time formulation. Section~\ref{sec:SpaceTimeROM} introduces
linear subspace solution representation in Section~\ref{sec:SolRep} and
space--time ROM formulation using Galerkin projection in
Section~\ref{sec:GalerkinProjection} and least--squares Petrov--Galerkin
projection in Section~\ref{sec:LSPGProjection}. Then, both space--time ROMs are
compared in Section~\ref{sec:GvsPG}. Section~\ref{sec:BasisGen} describes how to
generate space--time basis. We investigate block structures of space--time ROM
basis in Section~\ref{sec:BlockSpaceTimeBasis}. We introduce the block
structures of Galerkin space--time ROM operators derived in \cite{choi2020space}
in Section~\ref{sec:BlockG}. In Section~\ref{sec:BlockPG}, we derive
least-squares Petrov--Galerkin space--time ROM operators in terms of the blocks.
Then, we compared Galerkin and least--squares Petrov--Galerkin block structures
in Section~\ref{sec:BlockGvsBlcokPG}. We compute computational complexity of
forming the space--time ROM operators in Section~\ref{sec:ComputationalCost}.
The error analysis is presented in Section~\ref{sec:ErrorAnalysis}. We
demonstrate the performance of both Galerkin and least-squares Petrov--Galerkin
space--time ROMs in two numerical experiments in
Section~\ref{sec:NumericalResults}. Finally, the paper is concluded with summary
and future works in Section~\ref{sec:Conclusion}. Note that we use
``least--squares Petrov--Galerkin" and ``Petrov--Galerkin" interchangeably
throughout the paper. Appendix~\ref{sec:codes} presents six Python codes with
less than 120 lines that are used to generate our numerical results.

\section{Linear dynamical systems}\label{sec:LinearDynamicalSystems}
We consider the parameterized linear dynamical system shown in Equation
\eqref{eq:Dynamicalstate}.
\begin{align}
\frac{\partial {\state}(\timevar;\param)}{\partial t} &= \statemat(\param)
  \state(\timevar;\param) + \inputmat(\param)\inputvec(\timevar;\param), \quad 
\state(0;\param) = \initialstate(\param),\label{eq:Dynamicalstate} 
\end{align}
where $\param\in\paramspace\subset\RR{\nparam}$ denotes a parameter vector,
$\state: [0,\finaltime]\times\RR{\nparam} \rightarrow \RR{\nspace}$ denotes a
time dependent state variable function,
$\initialstate:\RR{\nparam}\rightarrow\RR{\nspace}$ denotes an initial state,
and $\inputvec: [\initialtime,\finaltime]\times\RR{\nparam} \rightarrow
\RR{\ninput}$ denotes a time dependent input variable function. The operators
$\statemat:\RR{\nparam}\rightarrow\RR{\nspace\times\nspace}$,
$\inputmat:\RR{\nparam}\rightarrow\RR{\nspace\times\ninput}$, and are real
valued matrices that are independent of state variables.  
  
Although any time integrator can be used, for the demonstration purpose, we
choose to apply a backward Euler time integration scheme shown in Equation
\eqref{eq:backwardEuler}:
\begin{equation}\label{eq:backwardEuler}
\left (\IdentityNs- \timestepk{k} \statemat(\param) \right )\statek{k} =
  \statek{k-1} + \timestepk{k} \inputmat(\param) \inputveck{k}(\param),
\end{equation}
where $\IdentityNs\in\RR{\nspace\times\nspace}$ is the identity matrix,
$\timestepk{k}$ is the $k$th time step size with
$\finaltime=\sum_{k=1}^{\ntime}\timestepk{k}$ and
$\timevark{k}=\sum_{j=1}^{k}\timestepk{j}$, and
$\statek{k}(\param):=\state(\timevark{k};\param)$ and
$\inputveck{k}(\param):=\inputvec(\timevark{k};\param)$ are the state and input
vectors at $k$th time step where $k\in\nat{\ntime}$. The Full Order Model (FOM)
solves Equation \eqref{eq:backwardEuler} for every time step, where its spatial
dimension is $\nspace$ and the temporal dimension is $\ntime$. Each time step of
the FOM can be written out and put in another matrix system shown in Equation
\eqref{eq:st-fom}.  This is known as the space--time formulation.   
\begin{equation}\label{eq:st-fom}
\ststatemat(\param)\ststate(\param) = \stinputvec(\param) +
\stinitial(\param),
\end{equation}
where 
\begin{gather}
\ststatemat(\param) = \bmat{
  \IdentityNs - \timestepk{1} \statemat(\param) & & & \\
 -\IdentityNs & \IdentityNs - \timestepk{2} \statemat(\param) & & \\
   & \ddots & \ddots & \\
   &   & -\IdentityNs & \IdentityNs -\timestepk{\ntime} \statemat(\param)
   } \label{eq:stsystemmat}, \\
\ststate(\param)  = \bmat{ \statek{1}(\param) \\ \statek{2}(\param) \\
\vdots \\ \statek{\ntime}(\param) 
}\label{eq:ststate} , \\
\stinputvec(\param) = \bmat{
  \timestepk{1} \inputmat(\param)\inputveck{1}(\param)  \\ \timestepk{2}
  \inputmat(\param)\inputveck{2}(\param) \\ \vdots \\ \timestepk{\ntime}
  \inputmat(\param)\inputveck{\ntime}(\param) 
} \label{eq:stinputvec}, \\
\stinitial(\param) = \bmat{
  \initialstate(\param) \\ \zero \\ \vdots \\ \zero
}\label{eq:stinitialstate}.
\end{gather}

The space--time system matrix $\ststatemat$ has dimensions $\RR{\nparam}
\rightarrow\RR{\nspace\ntime \times \nspace\ntime}$, the space--time state
vector $\ststate$ has dimensions  $\RR{\nparam} \rightarrow\RR{\nspace\ntime}$,
the space--time input vector $\stinputvec$ has dimensions  $\RR{\nparam}
\rightarrow\RR{\nspace\ntime}$, and the space--time initial vector $\stinitial$
has dimensions  $\RR{\nparam} \rightarrow\RR{\nspace\ntime}$. Although it seems
that the solution can be found in a single solve, in practice there is no
computational saving gained from doing so since the block structure of the
space--time system will solve the system in a time-marching fashion anyways.
However, we formulate the problem in this way since our reduced order model
(ROM) formulation can reduce and solve the space--time system efficiently. In
the following sections, we describe the parametric Galerkin and least--squares
Petrov--Galerkin ROM formulations.

\section{Space--time reduced order models}\label{sec:SpaceTimeROM}
We investigate two projection-based space--time ROM formulations:
the Galerkin and least--squares Petrov--Galerkin formulations. Here, we use
``least-squares Petrov--Galerkin" and ``Petrov--Galerkin" interchangeably
throughout the paper.

\subsection{Linear subspace solution representation}\label{sec:SolRep}
Both the Galerkin and Petrov--Galerkin methods reduce the number of space--time
degrees of freedom by approximating the space--time state variables as a smaller
linear combination of space--time basis vectors:
\begin{equation}\label{eq:ststateapproximation}
\ststate(\param) \approx \approxststate(\param) \equiv 
\spacetimebasis\reducedststate(\param),
\end{equation}
where $\reducedststate(\param): \RR{\nparam} \rightarrow
\RR{\nreducedspace\nreducedtime}$ with $\nreducedspace \ll \nspace$ and
$\nreducedtime \ll \ntime$. The space--time basis,
$\spacetimebasis\in\RR{\nspace\ntime \times \nreducedspace\nreducedtime}$ is
defined as 
\begin{equation}\label{eq:spacetimebasis}
\spacetimebasis \equiv \bmat{\spacetimebasisveck{1} & \cdots &
  \spacetimebasisveck{\reducedspaceindex+\nreducedspace(\reducedtimeindex-1)}  &
  \cdots \spacetimebasisveck{\nreducedspace\nreducedtime}},
\end{equation}
where $\reducedspaceindex\in\nat{\nreducedspace}$,
$\reducedtimeindex\in\nat{\nreducedtime}$.  Substituting Equation
\eqref{eq:ststateapproximation} into the space--time formulation in Equation
\eqref{eq:st-fom} gives an over-determined system of equations:
\begin{equation}\label{eq:overdetermined-spacetime}
\ststatemat(\param)\spacetimebasis\reducedststate(\param) = \stinputvec(\param)
  + \stinitial(\param) 
\end{equation}
This over-determined system of equations can be closed by either the Galerkin or
Petrov--Galerkin projections. 

\subsection{Galerkin projection}\label{sec:GalerkinProjection}
In the Galerkin formulation, Equation \ref{eq:overdetermined-spacetime} is
closed by the Galerkin projection, where both sides of the equation is
multiplied by $\spacetimebasis^T$. Thus, we solve following reduced system for
the unknown generalized coordinates, $\ststate$:
\begin{equation}
    \spacetimebasis^{T} \Astparameter \spacetimebasis
    \reducedststate(\param)=\spacetimebasis^{T}
    \stinputvec(\param)+\spacetimebasis^{T}  \stinitial(\param).
\end{equation}

For notational simplicity, let us define the reduced space--time system matrix
as $\AhatGalerkin := \spacetimebasis^{T} \Astparameter \spacetimebasis $,
reduced space--time input vector as $\fhatGalerkin := \spacetimebasis^{T}
\stinputvec(\param)$, and reduced space--time initial state vector as
$\uhatGalerkin := \spacetimebasis^{T}  \stinitial(\param)$.

\subsection{Least-squares Petrov--Galerkin projection}\label{sec:LSPGProjection}
In the least-squares Petrov--Galerkin formulation, we first define the
space--time residual as
\begin{equation}\label{eq:LSPGRes}
    \approxstres(\reducedststate;\param):=\stinputvec(\param)
    +\stinitial(\param) - \ststatemat(\param)\spacetimebasis\reducedststate
\end{equation}
where $\approxstres: \RR{\nreducedspace\nreducedtime}\times\RR{\nparam}
\rightarrow \RR{\nspace\ntime}$. Note Equation~\eqref{eq:LSPGRes} is an
over-determined system. To close the system and solve for the unknown
generalized coordinates, $\reducedststate$, the least-squares Petrov--Galerkin
method takes the squared norm of the residual vector function and minimize it:
\begin{equation}\label{eq:LSPGResSol}
    \reducedststate =
    \argmin{\reducedststatedummy\in\RR{\nreducedspace\nreducedtime}} \quad
    \frac{1}{2} \left \|\approxstres(\reducedststatedummy;\param) \right \|_2^2.
\end{equation}
The solution to Equation~\eqref{eq:LSPGResSol} satisfies
\begin{equation}
    (\Astparameter\spacetimebasis)^T\approxstres(\reducedststate;\param)=\zero
\end{equation}
leading to
\begin{equation}
    \spacetimebasis^{T}  \Astparameter^{T} \Astparameter \spacetimebasis
    \reducedststate(\param)=\spacetimebasis^{T}  \Astparameter^{T}
    \stinputvec(\param)+ \spacetimebasis^{T} \Astparameter^T \stinitial(\param).
\end{equation}

For notational simplicity, let us define the reduced space--time system matrix
as 
\[ \AhatPetrovGalerkin := \spacetimebasis^{T}  \Astparameter^{T} \Astparameter
\spacetimebasis, \] 
reduced space--time input vector as $\fhatPetrovGalerkin :=
\spacetimebasis^{T}  \Astparameter^{T} \stinputvec(\param) $, and reduced
space--time initial state vector as $\uhatPetrovGalerkin := \spacetimebasis^{T}
\Astparameter^T \stinitial(\param)$.

\subsection{Comparison of Galerkin and Petrov--Galerkin
projections}\label{sec:GvsPG}
The reduced space--time system matrices, reduced space--time input vectors, and
reduced space--time initial state vectors for Galerkin and Petrov--Galerkin
projections are presented in Table \ref{tab:table_formulations}.
\begin{table} [H]
\caption{Comparison of Galerkin and Petrov--Galerkin projections}
\label{tab:table_formulations}
\centering
\begin{tabular}{ |c|c|c| }
\hline
\bf{Galerkin} & \bf{Petrov--Galerkin}  \\ 
\hline
$\AhatGalerkin = \spacetimebasis^{T} \Astparameter \spacetimebasis $ &
  $\AhatPetrovGalerkin= \spacetimebasis^{T}  \Astparameter^{T} \Astparameter
  \spacetimebasis $  \\ 
\hline
$\fhatGalerkin = \spacetimebasis^{T} \stinputvec (\param)$ &
  $\fhatPetrovGalerkin= \spacetimebasis^{T}  \Astparameter^{T}
  \stinputvec(\param) $  \\ 
\hline
$\uhatGalerkin = \spacetimebasis^{T} \stinitial(\param)$ &  $\uhatPetrovGalerkin
  = \spacetimebasis^{T} \Astparameter^T \stinitial (\param)$ \\
\hline
\end{tabular}
\end{table}

\section{Space-time Basis Generation}\label{sec:BasisGen}
In this section, we repeat Section 4.1 in \cite{choi2020space} to be
self-contained.

\begin{figure}[H]
\centering
\includegraphics[width=0.8\textwidth]{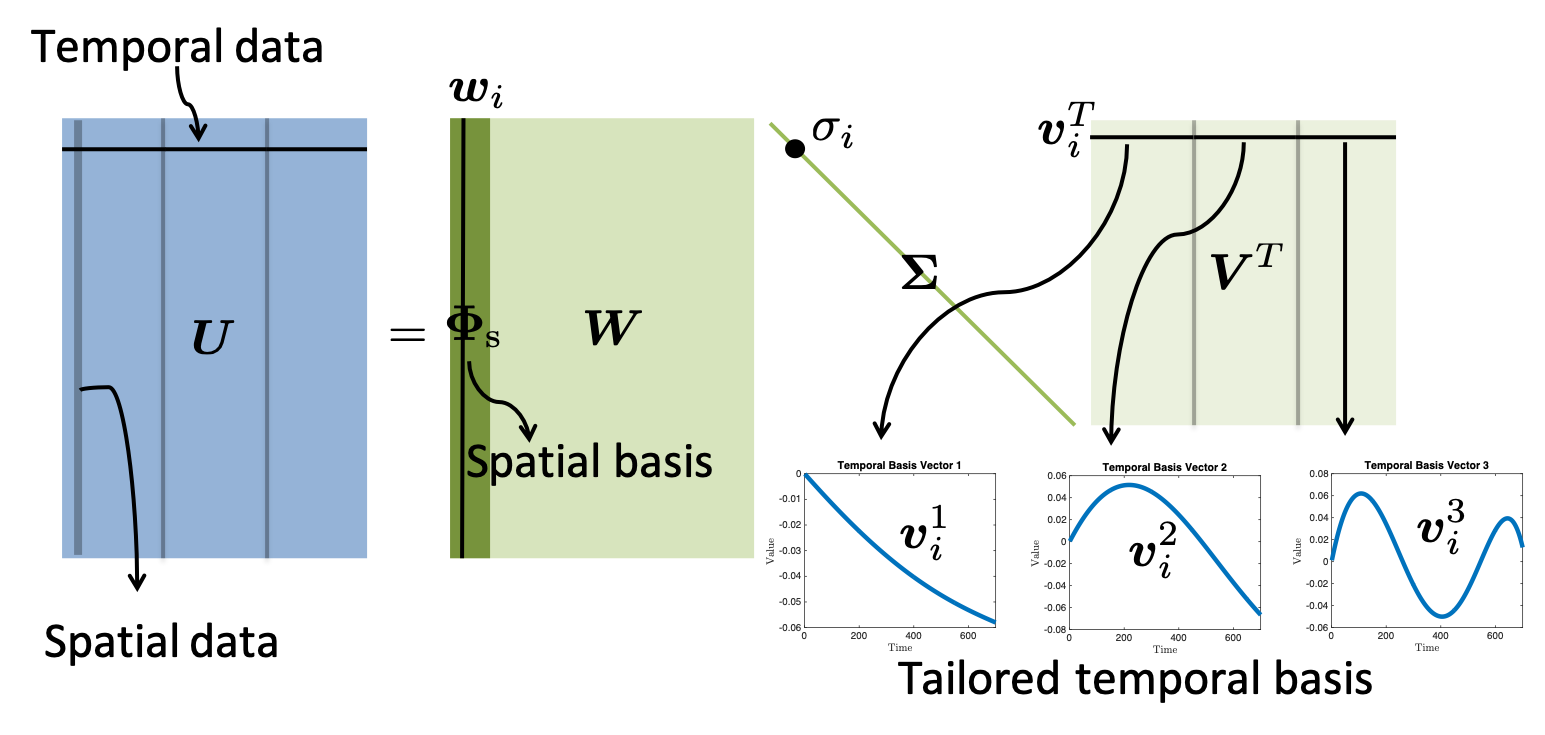}
\caption{Illustration of spatial and temporal bases construction, using SVD with
  $\nparam=3$. The right singular vector,
  $\rightsingularveck{\reducedspaceindex}$, describes three different temporal
  behaviors of a left singular basis vector
  $\leftsingularveck{\reducedspaceindex}$, i.e., three different temporal
  behaviors of a spatial mode. Each temporal behavior is denoted as
  $\rightsingularveck{\reducedspaceindex}^1$,
  $\rightsingularveck{\reducedspaceindex}^2$, and
  $\rightsingularveck{\reducedspaceindex}^3$.  }
\label{fig:SVDexplanation}
\end{figure}

We follow the method of snapshots described by Sirovich
\cite{sirovich1987turbulence}. First, let
$\{\paramk{1},\ldots,\paramk{\nparam}\}$ be a set of parameter samples, where we
run full order model simulations. Let $\stateSnapshotk{p} \equiv
\bmat{\statek{1}(\paramk{p}) & \cdots &
\statek{\ntime}(\paramk{p})}\in\RR{\nspace\times\ntime}$, $p\in\nat{\nparam}$,
be a full order model solution matrix for a sample parameter,
$\paramk{p}\in\paramspace$. Then concatenating all the solution matrix defines a
snapshot matrix, $\stateSnapshot\in\RR{\nspace\times\nparam\ntime}$, i.e., 
\begin{equation}\label{eq:snapshotmat}
\stateSnapshot \equiv \bmat{\stateSnapshotk{1} & \cdots &
\stateSnapshotk{\nparam}}.
\end{equation}
We use Proper Orthogonal Decomposition (POD) to construct the spatial basis,
$\spatialbasis$. POD \cite{berkooz1993proper} obtains $\spatialbasis$ by
choosing the leading $\nreducedspace$ columns of the left singular matrix,
$\leftsingularmat$, of the following Singular Value Decomposition (SVD) with
$\ell \equiv \min(\nspace,\nparam\ntime)$ and $\nreducedspace < \nparam\ntime$:
\begin{align}
\stateSnapshot &= \leftsingularmat\singularvalmat\rightsingularmat^T
\label{eq:SVD}\\
              &= \sum_{\reducedspaceindex=1}^{\ell}\singularvaluesArg{i}
              \leftsingularveck{\reducedspaceindex}
              \rightsingularveck{\reducedspaceindex}^T,\label{eq:SVDsummation}
\end{align} 
where $\leftsingularmat\in\RR{\nspace\times\ell}$ and
$\rightsingularmat\in\RR{\nparam\ntime\times\ell}$ are orthogonal matrices and
$\singularvalmat\in\RR{\ell\times\ell}$ is a diagonal matrix with singular
values on its diagonal. The spatial POD basis, $\spatialbasis$ minimizes 
\begin{equation}\label{eq:POD}
\left \|\stateSnapshot -
  \spatialbasis\spatialbasis^T{\stateSnapshot} \right \|_F^2
\end{equation}
over  all $\spatialbasis\in\RR{\nspace\times\nreducedspace}$ with orthonormal
columns, where $\|\cdot\|_F$ denotes the Frobenius norm. The POD procedure seeks
the $\nreducedspace$-dimensional subspace that optimally represents the solution
snapshot, $\stateSnapshot$.
The equivalent summation form is written in \eqref{eq:SVDsummation}, where
$\singularvaluesArg{\reducedspaceindex}\in\RR{}$ is $\reducedspaceindex$th
singular value, $\leftsingularveck{\reducedspaceindex}$ and
$\rightsingularveck{i}$ are $i$th left and right singular vectors, respectively.
Note that $\rightsingularveck{\reducedspaceindex}$ describes $\nparam$ different
temporal behavior of $\leftsingularveck{\reducedspaceindex}$. For example,
Figure~\ref{fig:SVDexplanation} illustrates the case of $\nparam=3$, where
$\rightsingularveck{\reducedspaceindex}^1$,
$\rightsingularveck{\reducedspaceindex}^2$, and
$\rightsingularveck{\reducedspaceindex}^3$ describe three different temporal
behavior of a specific spatial basis vector, i.e.,
$\leftsingularveck{\reducedspaceindex}$. For general $\nparam$, we note that
$\rightsingularveck{\reducedspaceindex}$ describes $\nparam$ different temporal
behavior of $\reducedspaceindex$th spatial basis vector, i.e.,
$\spatialbasisveck{\reducedspaceindex} = \leftsingularveck{\reducedspaceindex}$.
We set $\temporalSnapshotk{\reducedspaceindex} =
\bmat{\rightsingularveck{\reducedspaceindex}^1 & \cdots &
\rightsingularveck{\reducedspaceindex}^{\nparam}   }$ to be
$\reducedspaceindex$th temporal snapshot matrix, where
$\rightsingularveck{\reducedspaceindex}^k \equiv [
  \rightsingularveck{\reducedspaceindex}(1+(k-1)\ntime),
  \rightsingularveck{\reducedspaceindex}(2+(k-1)\ntime), \cdots,
  \rightsingularveck{\reducedspaceindex}(k\ntime)]^T \in \RR{\ntime}$ for
$k\in\nat{\nparam}$, where
$\rightsingularveck{\reducedspaceindex}(j), j\in\nat{\nparam\ntime}$ is the
$j$th component of the vector.  The SVD of
$\temporalSnapshotk{\reducedspaceindex}$ is
\begin{equation}\label{eq:temporalSVD}
\temporalSnapshotk{\reducedspaceindex} =
  \leftsingularmattemporalk{\reducedspaceindex}
  \singularvalmatk{\reducedspaceindex}
  \rightsingularmattemporalk{\reducedspaceindex}^T.
\end{equation}
Then, choosing the leading $\nreducedtime$ vectors of
$\leftsingularmattemporalk{\reducedspaceindex}$ yields the temporal basis,
$\temporalbasisi{\reducedspaceindex}$ for $\reducedspaceindex$th spatial basis
vector. Finally, we can construct a space--time basis vector,
$\spacetimebasisveck{\reducedspaceindex+\nreducedspace(\reducedtimeindex-1)}
\in\RR{\nspace\ntime}$,
in Equation~\eqref{eq:spacetimebasis} as
\begin{equation}\label{eq:spacetimebasisvec}
\spacetimebasisveck{\reducedspaceindex+\nreducedspace(\reducedtimeindex-1)} =
  \temporalbasisveck{\reducedspaceindex}{\reducedtimeindex} \otimes
  \spatialbasisveck{\reducedspaceindex},
\end{equation}
where $\otimes$ denotes Kronecker product,
$\spatialbasisveck{\reducedspaceindex} \in \RR{\nspace}$ is
$\reducedspaceindex$th vector of the spatial basis, $\spatialbasis$, and
$\temporalbasisveck{\reducedspaceindex}{\reducedtimeindex} \in\RR{\ntime}$ is
$\reducedtimeindex$th vector of the temporal basis,
$\temporalbasisi{\reducedspaceindex}$ that describes a temporal behavior of
$\spatialbasisveck{\reducedspaceindex}$.

\section{Space-time reduced order models in block structure}
\label{sec:BlockStructure}
We avoid building the space--time basis vector defined in
Equation~\eqref{eq:spacetimebasisvec} because it requires much memory for
storage. Thus, we can exploit the block structure of the matrices to save
computational cost and storage of the matrices in memory.
Section~\ref{sec:BlockG} introduces such block structures for the space--time
Galerkin projection, while Section~\ref{sec:BlockPG} shows block structures for
the space--time Petrov--Galerkin projection. First, we introduce common block
structures that appear both the Galerkin and Petrov--Galerkin projections in
Section~\ref{sec:BlockSpaceTimeBasis}.

\subsection{Block structures of space--time
basis}\label{sec:BlockSpaceTimeBasis}
Following \cite{choi2020space}'s notation, we define the block structure of the
space--time basis to be:
\begin{equation}\label{eq:Phist}
 \basismatspace_{st} =
 \begin{pmatrix}
  \basismatspace_{s} \bds{D}_1^1 & \cdots & \cdots & \cdots  &
   \basismatspace_{s} \bds{D}_1^{n_t} \\
  \vdots & \ddots & \vdots & \iddots  & \vdots \\
  \vdots  & \cdots  & \basismatspace_{s} \bds{D}_{k}^j  & \cdots  & \vdots  \\
  \vdots  & \iddots  & \vdots & \ddots  & \vdots  \\
  \basismatspace_{s} \bds{D}_{N_t}^1 & \cdots & \cdots & \cdots  &
   \basismatspace_{s} \bds{D}_{N_t}^{n_t}
 \end{pmatrix} \in \mathbb R^{N_s N_t \times n_s n_t}
  \end{equation}
where the $k$th time step of the temporal basis matrix is a diagonal matrix
defined as 

\begin{equation}\label{eq:D}
 \bds{D}_k^j =
 \begin{pmatrix}
  \phi_{1j,k}^t  & 0 & \cdots & \cdots &0  \\
  0 & \ddots &0 & \cdots & \vdots   \\
  \vdots  & 0  & \phi_{ij,k}^t  & \cdots & \vdots    \\
   \vdots  & \vdots  & \ddots & \ddots & 0    \\
  0& 0 & \cdots & 0 &  \phi_{n_s j,k}^t
 \end{pmatrix} \in \mathbb R^{n_s \times n_s }
 \end{equation}
where $\phi_{ij,k}^t  \in \mathbb R$ is the $k$th element of
$\boldsymbol{\phi}_{ij}^t \in \mathbb R^{N_t}$.  

\subsection{Block structures of Galerkin projection}\label{sec:BlockG}
As shown in Table~\ref{tab:table_formulations}, the reduced space--time Galerkin
system matrix, $\AhatGalerkin $ is:
\begin{equation}\label{eq:Astg}
\begin{aligned}
       \AhatGalerkin= &\spacetimebasis^{T} \Astparameter \spacetimebasis  \\
    \end{aligned}
 \end{equation}
Now, We define the block structure of this matrix as: 
\begin{equation}\label{eq:AstgBlock}
\begin{aligned}
   \AhatGalerkin= \begin{pmatrix}
\hat{\bds{A}}_{(1,1)}^{st,g}(\mubold) & \cdots & \cdots & \cdots  &
     \hat{\bds{A}}_{(1,n_t)}^{st,g}(\mubold)  \\
\vdots & \ddots & \vdots & \iddots  & \vdots \\
\vdots  & \cdots  &  \AhatPetrovGalerkinBlock  & \cdots  & \vdots  \\
\vdots  & \iddots  & \vdots & \ddots  & \vdots  \\
\hat{\bds{A}}_{(n_t,1)}^{st,g}(\mubold)  & \cdots & \cdots & \cdots  &
     \hat{\bds{A}}_{(n_t,n_t)}^{st,g}(\mubold) 
\end{pmatrix}
\end{aligned}
\end{equation}
so that we can exploit the block structure of these matrices such that we do not
need to form the entire matrix. We derive that $\AhatGalerkinBlock \in  \mathbb
R^{n_s \times n_s}$ where the $\jprimej$th block matrix is:
\begin{equation}\label{eq:AstgBlockSum}
  \AhatGalerkinBlock =  \sum\limits_{k=1}^{N_t} \Big(
  \temporalBasisMatrixk^{j^\prime}\temporalBasisMatrixk^{j} - \timestepk{k}
  \temporalBasisMatrixk^{j^\prime} \basismatspace_s^T \Amatrix\basismatspace_s
  \temporalBasisMatrixk^{j} \Big) - \sum\limits_{k=1}^{N_t-1}
  \temporalBasisMatrixkpone^{j^\prime}\temporalBasisMatrixk^j 
 \end{equation}
 
The reduced space--time Galerkin input vector $\fhatGalerkin \in \mathbb R^{n_s
n_t}$ is
\begin{equation}\label{eq:fstg}
   \fhatGalerkin= \spacetimebasis^{T} \stinputvec(\param). 
\end{equation}
Again, utilizing the block structure of matrices, we compute $j$th block vector
$\fhatGalerkin_{(j)} \in \mathbb R^{n_t}$ to be:
\begin{equation}\label{eq:fstgSum}
   \fhatGalerkin_{(j)} = \sum\limits_{k=1}^{N_t} \timestepk{k}
   \temporalBasisMatrixk^j \basismatspace_s^T \Bhatmatrix
   \inputveck{k}(\param).
\end{equation}
 
Finally, the space--time Galerkin initial vector, $\uhatGalerkin \in \mathbb
R^{n_s n_t}$, can be computed as:
\begin{equation}\label{eq:ustg}
   \uhatGalerkin = \spacetimebasis^{T} \stinitial(\param)
\end{equation}
where the $j$th block vector, $\uhatGalerkin_{(j)} \in \mathbb R^{n_s}$, is:
\begin{equation}\label{eq:ustgSum}
   \uhatGalerkin_{(j)}  = \bds{D}_{1}^{j} \basismatspace_s^T \uinit.
\end{equation}

\subsection{Block structures of least-squares Petrov--Galerkin projection}
\label{sec:BlockPG}
As shown in Table~\ref{tab:table_formulations}, the reduced space--time
Petrov--Galerkin system matrix, $\AhatPetrovGalerkin $ is:
\begin{equation}\label{eq:Astpg}
\begin{aligned}
   \AhatPetrovGalerkin= &\basismatspace_{st}^{T}  \Astparameter^{T}
   \Astparameter \basismatspace_{st}  \\
\end{aligned}
\end{equation}
Now, We define the block structure of this matrix as: 
\begin{equation}\label{eq:AstpgBlock}
\begin{aligned}
   \AhatPetrovGalerkin= \begin{pmatrix}
\hat{\bds{A}}_{(1,1)}^{st,pg}(\mubold) & \cdots & \cdots & \cdots  &
     \hat{\bds{A}}_{(1,n_t)}^{st,pg}(\mubold)  \\
\vdots & \ddots & \vdots & \iddots  & \vdots \\
\vdots  & \cdots  &  \AhatPetrovGalerkinBlock  & \cdots  & \vdots  \\
\vdots  & \iddots  & \vdots & \ddots  & \vdots  \\
\hat{\bds{A}}_{(n_t,1)}^{st,pg}(\mubold)  & \cdots & \cdots & \cdots  &
     \hat{\bds{A}}_{(n_t,n_t)}^{st,pg}(\mubold) 
\end{pmatrix}
\end{aligned}
\end{equation}
so that we can exploit the block structure of these matrices such that we do not
need to form the entire matrix. We derive that $\AhatPetrovGalerkinBlock \in
\mathbb R^{n_s \times n_s}$ where the $\jprimej$th block matrix is:
\begin{equation}\label{eq:AstpgBlockSum}
	\begin{aligned}
      \AhatPetrovGalerkinBlock =  & \sum\limits_{k=1}^{N_t} \Big[
        \temporalBasisMatrixk^{j^\prime}   \basismatspace_s^T \Big(\IdentityNs -
        \Delta t^{(k)} \Amatrix^T    \Big)  \Big(\IdentityNs - \Delta t^{(k)}
        \Amatrix    \Big) \basismatspace_s  \temporalBasisMatrixk^{j} \Big]  \\
      & + \sum\limits_{k=1}^{N_t-1} \Big[     \temporalBasisMatrixk^{j^\prime}
      \temporalBasisMatrixk^{j}  -   \temporalBasisMatrixk^{j^\prime}
      \basismatspace_s^T \Big(\IdentityNs - \Delta t^{(k+1)} \Amatrix \Big)
      \basismatspace_s \temporalBasisMatrixkpone^j  \\
     & - \temporalBasisMatrixkpone^{j^{\prime}}  \basismatspace_s^T
     \Big(\IdentityNs - \Delta t^{(k+1)} \Amatrix^T    \Big) \basismatspace_s
     \temporalBasisMatrixk^{j} \Big] 
      \end{aligned}
 \end{equation}
 
The reduced space--time Petrov--Galerkin input vector $\fhatPetrovGalerkin \in
\mathbb R^{n_s n_t}$ is
\begin{equation}\label{eq:fstpg}
   \fhatPetrovGalerkin= \basismatspace_{st}^{T}  \Astparameter^{T}
   \stinputvec(\param). 
\end{equation}
Again, utilizing the block structure of matrices, we compute $j$th block vector
$\fhatPetrovGalerkin_{(j)} \in \mathbb R^{n_t}$ to be:
\begin{equation}\label{eq:fstpgSum}
 \begin{aligned}
   \fhatPetrovGalerkin_{(j)} = & \sum\limits_{k=1}^{N_t} \Big[
     \temporalBasisMatrixk^j   \basismatspace_s^T   \Big(\IdentityNs -
     \timestepk{k} \Amatrix^T    \Big)  \timestepk{k}    \Bhatmatrix
     \inputveck{k}(\param)   \Big] \\
&   -  \sum\limits_{k=1}^{N_t-1} \Big[   \temporalBasisMatrixk^j
   \basismatspace_s^T \timestepk{k+1}  \Bhatmatrix  \inputveck{k+1}(\param)
   \Big].
   \end{aligned}
\end{equation}
 
Finally, the space--time Petrov--Galerkin initial vector, $\uhatPetrovGalerkin \in
\mathbb R^{n_s n_t}$, can be computed as:
\begin{equation}\label{eq:ustpg}
   \uhatPetrovGalerkin = \basismatspace_{st}^{T} \Astparameter^T
   \stinitial(\param)
\end{equation}
where the $j$th block vector, $\uhatPetrovGalerkin_{(j)} \in \mathbb R^{n_s}$,
is:
\begin{equation}\label{eq:ustpgSum}
   \uhatPetrovGalerkin_{(j)}  = \bds{D}_{1}^{j}    \basismatspace_s^T
   \Big(\IdentityNs - \timestepk{1} \Amatrix^T    \Big) \uinit.
\end{equation}
 
\subsection{Comparison of Galerkin and Petrov--Galerkin block structures}
\label{sec:BlockGvsBlcokPG}
The block structures of space--time reduced order model operators are summarized
in Table~\ref{tab:table_blockstructures}.
\begin{table}[H]
\caption{Comparison of Galerkin and Petrov--Galerkin block structures}
\label{tab:table_blockstructures}
\resizebox{\textwidth}{!}{\begin{tabular}{ |c|c| } 
 \hline
 \bf{Galerkin} & \bf{Petrov--Galerkin}  \\ 
  \hline
 $\AhatGalerkin_{\jprimej} = $ & $\AhatPetrovGalerkin_{\jprimej}=  $  \\ 
 $\sum\limits_{k=1}^{N_t} \Big( \temporalBasisMatrixk^{j^\prime}
  \temporalBasisMatrixk^{j}  - \timestepk{k} \temporalBasisMatrixk^{j^\prime}
  \basismatspace_s^T  \Amatrix \basismatspace_s \temporalBasisMatrixk^{j} \Big)$
  & $\sum\limits_{k=1}^{N_t} \Big[     \temporalBasisMatrixk^{j^\prime}
  \basismatspace_s^T \Big(\IdentityNs - \timestepk{k} \Amatrix^T    \Big)
  \Big(\IdentityNs - \timestepk{k} \Amatrix    \Big) \basismatspace_s
  \temporalBasisMatrixk^{j} \Big]  $ \\
 $- \sum\limits_{k=1}^{N_t-1} \temporalBasisMatrixkpone^{j^{\prime}}
  \temporalBasisMatrixk^{j} $&$+ \sum\limits_{k=1}^{N_t-1} \Big[
    \temporalBasisMatrixk^{j^\prime}   \temporalBasisMatrixk^{j}  -
    \temporalBasisMatrixk^{j^\prime}  \basismatspace_s^T \Big(\IdentityNs -
    \timestepk{k+1} \Amatrix \Big) \basismatspace_s \temporalBasisMatrixkpone^j
    $\\
 &$- \temporalBasisMatrixkpone^{j^{\prime}}  \basismatspace_s^T \Big(\IdentityNs
  - \timestepk{k+1} \Amatrix^T    \Big) \basismatspace_s
  \temporalBasisMatrixk^{j} \Big] $\\
  \hline
 $\fhatGalerkin_{(j)} = $ & $\fhatPetrovGalerkin_{(j)}= \sum\limits_{k=1}^{N_t}
  \Big[     \temporalBasisMatrixk^j   \basismatspace_s^T   \Big(\IdentityNs -
  \timestepk{k} \Amatrix^T    \Big)  \timestepk{k} \Bhatmatrix
  \inputveck{k}(\param)   \Big]  $  \\ 
 $\sum\limits_{k=1}^{N_t}  \temporalBasisMatrixk^j  \timestepk{k}
  \basismatspace_s^T  \Bhatmatrix  \inputveck{k}(\param) $ & $ -
  \sum\limits_{k=1}^{N_t-1} \Big[   \temporalBasisMatrixk^j   \basismatspace_s^T
  \timestepk{k+1}  \Bhatmatrix  \inputveck{k+1}(\param)  \Big]$\\
 \hline
 $\uhatGalerkin_{(j)} = \bds{D}_{1}^{j} \basismatspace_s^T \uinit $ &
  $\uhatPetrovGalerkin_{(j)} = \bds{D}_{1}^{j}  \basismatspace_s^T
  \Big(\IdentityNs - \timestepk{1} \Amatrix^T    \Big) \uinit$ \\
 \hline
\end{tabular}}
\end{table}

\subsection{Computational complexity of forming space--time ROM
operators}\label{sec:ComputationalCost}
To compute computational complexity of forming the reduced space--time system
matrices, input vectors, and initial state vectors for Galerkin and
Petrov--Galerkin projections, we assume that $\Aparameter \in
\RR{\nspace\times\nspace}$ is a band matrix with the bandwidth, $b$ and
$\Bparameter$ is a identity matrix, $\IdentityNs$. The band structure of
$\Aparameter$ is often seen in mathematical models because of local
approximations to derivative terms. Then, the bandwidth of $\Astparameter\in
\RR{\nspace\ntime\times\nspace\ntime}$ formed with backward Euler scheme is
$\nspace$. We also assume that the spatial basis vectors
$\spatialbasisvec_{i}\in\RR{\nspace}, i\in\nat{\nreducedspace}$ and temporal
basis vectors $\temporalbasisvec_{ij}\in\RR{\ntime}, i\in\nat{\nreducedspace}
\text{ and } j\in\nat{\nreducedtime}$ are given.

Let us start to compute the computational cost without use of block structures.
Constructing space--time basis costs
$\mathcal{O}(\nspace\ntime\nreducedspace\nreducedtime)$. For Galerkin
projections, computing the reduced space--time system matrix, input vectors, and
initial state vector costs
$\mathcal{O}(2\nspace^2\ntime\nreducedspace\nreducedtime)+\mathcal{O}(\nspace\ntime(\nreducedspace\nreducedtime)^2)$,
$\mathcal{O}(\nspace\ntime\nreducedspace\nreducedtime)$, and
$\mathcal{O}(\nspace\nreducedspace\nreducedtime)$, respectively. Thus, keeping
the dominant terms and taking off coefficient $2$ lead to
$\mathcal{O}(\nspace^2\ntime\nreducedspace\nreducedtime)+\mathcal{O}(\nspace\ntime(\nreducedspace\nreducedtime)^2)$.
With the assumption of $\nreducedspace\nreducedtime \ll \nspace$, we have
$\mathcal{O}(\nspace^2\ntime\nreducedspace\nreducedtime)$. For Petrov--Galerkin
projections, we first compute $\spacetimebasis^T\Astparameter$, resulting in
$\mathcal{O}(2\nspace^2\ntime\nreducedspace\nreducedtime)$. Then, computing the
reduced space--time system matrix, input vectors, and initial state vector costs
$\mathcal{O}(\nspace\ntime(\nreducedspace\nreducedtime)^2)$,
$\mathcal{O}(\nspace\ntime\nreducedspace\nreducedtime)$, and
$\mathcal{O}(\nspace\nreducedspace\nreducedtime)$, respectively. Thus, keeping
the dominant terms and taking off coefficient $2$ lead to
$\mathcal{O}(\nspace^2\ntime\nreducedspace\nreducedtime)+\mathcal{O}(\nspace\ntime(\nreducedspace\nreducedtime)^2)$.
With the assumption of $\nreducedspace\nreducedtime \ll \nspace$, we have
$\mathcal{O}(\nspace^2\ntime\nreducedspace\nreducedtime)$.

Now, let us compute the computational complexity with use of block structures.
For Galerkin projection, we first compute
$\basismatspace_s^T\Aparameter\basismatspace_s$ which will be reused, resulting
in $\mathcal{O}(2b\nspace\nreducedspace)+\mathcal{O}(\nspace\nreducedspace^2)$.
Then, we compute $\nreducedtime^2$ blocks for the reduced space--time system
matrix. It takes $\mathcal{O}(\ntime(2\nreducedspace^2+2\nreducedspace))$ to
compute each block. Thus, it costs
$\mathcal{O}(\ntime\nreducedtime^2(2\nreducedspace^2+2\nreducedspace))$ to
compute the reduced space--time system matrix. For the reduced space--time input
vector, $\nreducedtime$ blocks are needed and each block costs
$\mathcal{O}(\ntime(\nspace\nreducedspace+\nreducedspace))$, resulting in
$\mathcal{O}(\ntime\nreducedtime(\nspace\nreducedspace+\nreducedspace))$. The
computing the reduced space--time initial vector costs
$\mathcal{O}(\nspace\nreducedspace+\nreducedspace)$. Thus, keeping the dominant
terms and taking off coefficient $2$ lead to
$\mathcal{O}(b\nspace\nreducedspace)+\mathcal{O}(\nspace\nreducedspace^2)+\mathcal{O}(\ntime(\nreducedspace\nreducedtime)^2)+\mathcal{O}(\nspace\ntime\nreducedspace\nreducedtime)$.
With the assumptions of $b \ll \ntime\nreducedtime$, $\nreducedspace \ll
\ntime\nreducedtime$, and $\nreducedspace\nreducedtime \ll \nspace$, we have
$\mathcal{O}(\nspace\ntime\nreducedspace\nreducedtime)$. For Petrov--Galerkin
projection, we compute $(\IdentityNs-\timestep\Amatrix)\basismatspace_s$,
resulting in $\mathcal{O}(2b\nspace\nreducedspace)$. Then, we compute
$\basismatspace_s^T(\IdentityNs-\timestep\Amatrix^T)(\IdentityNs-\timestep\Amatrix)\basismatspace_s$
and $\basismatspace_s^T(\IdentityNs-\timestep\Amatrix)\basismatspace_s$ for
re-use. Each of them costs $\mathcal{O}(\nspace\nreducedspace^2)$. Now, we
compute $\nreducedtime^2$ blocks for the reduced space--time system matrix. It
takes $\mathcal{O}(\ntime(6\nreducedspace^2+\nreducedspace))$ to compute each
block. Thus, it costs
$\mathcal{O}(\ntime\nreducedtime^2(6\nreducedspace^2+\nreducedspace))$ to
compute the reduced space--time system matrix. For the reduced space--time input
vector, $\nreducedtime$ blocks are needed and each block costs
$\mathcal{O}(\ntime(\nspace\nreducedspace+\nreducedspace))$, resulting in
$\mathcal{O}(\ntime\nreducedtime(\nspace\nreducedspace+\nreducedspace))$. The
computing the reduced space--time initial vector costs
$\mathcal{O}(\nspace\nreducedspace+\nreducedspace)$. Thus, keeping the dominant
terms and taking off coefficients $2$ and $6$ lead to
$\mathcal{O}(b\nspace\nreducedspace)+\mathcal{O}(\nspace\nreducedspace^2)+\mathcal{O}(\ntime(\nreducedspace\nreducedtime)^2)+\mathcal{O}(\nspace\ntime\nreducedspace\nreducedtime)$.
With the assumptions of $b \ll \ntime\nreducedtime$, $\nreducedspace \ll
\ntime\nreducedtime$, and $\nreducedspace\nreducedtime \ll \nspace$, we have
$\mathcal{O}(\nspace\ntime\nreducedspace\nreducedtime)$.

In summary, the computational complexities of forming space--time ROM operators
in training phase for Galerkin and Petrov--Galerkin projections are presented in
Table~\ref{tab:table_costs}. We observe that a lot of computational costs are
reduced by making use of block structures for forming space--time reduced order
models.
\begin{table}[H]
\caption{Comparison of Galerkin and Petrov--Galerkin computational complexities}
\label{tab:table_costs}
\centering
  \begin{tabular}{ |c|c|c| } 
\hline
 & \bf{Galerkin} & \bf{Petrov--Galerkin}  \\ 
\hline
Not using block structures &
  $\mathcal{O}(\nspace^2\ntime\nreducedspace\nreducedtime)$ &
  $\mathcal{O}(\nspace^2\ntime\nreducedspace\nreducedtime)$ \\ 
\hline
Using block structures & $\mathcal{O}(\nspace\ntime\nreducedspace\nreducedtime)$
  & $\mathcal{O}(\nspace\ntime\nreducedspace\nreducedtime)$ \\
\hline
\end{tabular}
\end{table}

\section{Error analysis}\label{sec:ErrorAnalysis}
We present error analysis of the space--time ROM method. The error analysis is
based on \cite{choi2020space}. \textit{A posteriori} error bound is derived in
this section. Here, we drop the parameter dependence for notational simplicity.

{\theorem\label{theorem:STerrorbound}
We define the error at $k$th time step as $\errork{k} \equiv \statek{k} -
\approxstatek{k} \in \RR{\nspace}$ where $\statek{k} \in \RR{\nspace}$ denotes
FOM solution, $\approxstatek{k} \in \RR{\nspace}$ denotes approximate solution,
and $k \in \nat{\ntime}$. Let $\ststatemat \in \RR{\nspace\ntime \times
\nspace\ntime}$ be the space--time system matrix, $\resk{k} \in \RR{\nspace}$ be
the residual computed using FOM solution at $k$th time step, and $\approxresk{k}
\in \RR{\nspace}$ be the residual computed using approximate solution at $k$th
time step. For example, $\resk{k}$ and $\approxresk{k}$ after applying the
backward Euler scheme with the uniform time step become
\begin{align}
    \resk{k}(\statek{k},\statek{k-1}) &= \timestep\inputveck{k} + \statek{k-1} -
    (\identity-\timestep\statemat)\statek{k} = \zero \\
    \approxresk{k}(\approxstatek{k},\approxstatek{k-1}) &=
    \timestep\inputveck{k} + \approxstatek{k-1} -
    (\identity-\timestep\statemat)\approxstatek{k}
\end{align}
with $\approxstatek{0}=\initialstate$. Then, the error bound is given by
\begin{equation}\label{eq:STerrorbound}
  \max_{k\in\nat{\ntime}}\|\errork{k}\|_2 \leq \eta
  \max_{k\in\nat{\ntime}}\|\approxresk{k}\|_2
\end{equation}
where $\eta\equiv\sqrt{\ntime}\|(\ststatemat)^{-1}\|_2$ denotes the stability
constant.}

\begin{proof}
Let us define the space--time residual as
\begin{equation}
    \stres: \ststatedummy \mapsto
    \stinputvec+\stinitial-\ststatemat\ststatedummy
\end{equation}
with $\stres:\RR{\nspace\ntime}\rightarrow\RR{\nspace\ntime}$. Then, we have
\begin{align}
    \stres(\ststate)&=\stinputvec+\stinitial-\ststatemat\ststate=\zero
    \label{eq:stres}\\
    \stres(\approxststate)&=\stinputvec+\stinitial-\ststatemat\approxststate
    \label{eq:approxstres}
\end{align}
where $\ststate\in\RR{\nspace\ntime}$ is the space--time FOM solution and $\approxststate\in\RR{\nspace\ntime}$ is the approximate space--time solution.
Subtracting Equation~\eqref{eq:approxstres} from Equation~\eqref{eq:stres} gives
\begin{equation}
    \stres(\approxststate)=\ststatemat\sterror
\end{equation}
where $\sterror \equiv \ststate - \approxststate \in \RR{\nspace\ntime}$.
Inverting $\ststatemat$ yields
\begin{equation}
    \sterror = (\ststatemat)^{-1}\stres(\approxststate).
\end{equation}
Taking $\ell_2$ norm and H{\"o}lders' inequality gives
\begin{equation}
    \|\sterror\|_2 \leq \|(\ststatemat)^{-1}\|_2\|\stres(\approxststate)\|_2.
\end{equation}
We can re-write this in the following form
\begin{equation}
    \sqrt{\sum_{k=1}^{\ntime}\|\errork{k}\|_2^2} \leq \|(\ststatemat)^{-1}\|_2 \sqrt{\sum_{k=1}^{\ntime}\|\approxresk{k}\|_2^2}.
\end{equation}
Using the relations
\begin{equation}
  \max_{k\in\nat{\ntime}}\|\errork{k}\|_2^2 \leq
  \sum_{k=1}^{\ntime}\|\errork{k}\|_2^2
\end{equation}
and
\begin{equation}
    \sum_{k=1}^{\ntime}\|\approxresk{k}\|_2^2 \leq \ntime\max_{k\in\nat{\ntime}}\|\approxresk{k}\|_2^2,
\end{equation}
we have
\begin{equation}
    \max_{k\in\nat{\ntime}}\|\errork{k}\|_2 \leq
    \sqrt{\ntime}\|(\ststatemat)^{-1}\|_2
    \max_{k\in\nat{\ntime}}\|\approxresk{k}\|_2,
\end{equation}
which is equivalent to the error bound in \eqref{eq:STerrorbound}.
\end{proof}

\begin{figure}[!htbp]
    \centering
    \begin{subfigure}[b]{0.49\textwidth}
        \centering
        \includegraphics[width=\textwidth]{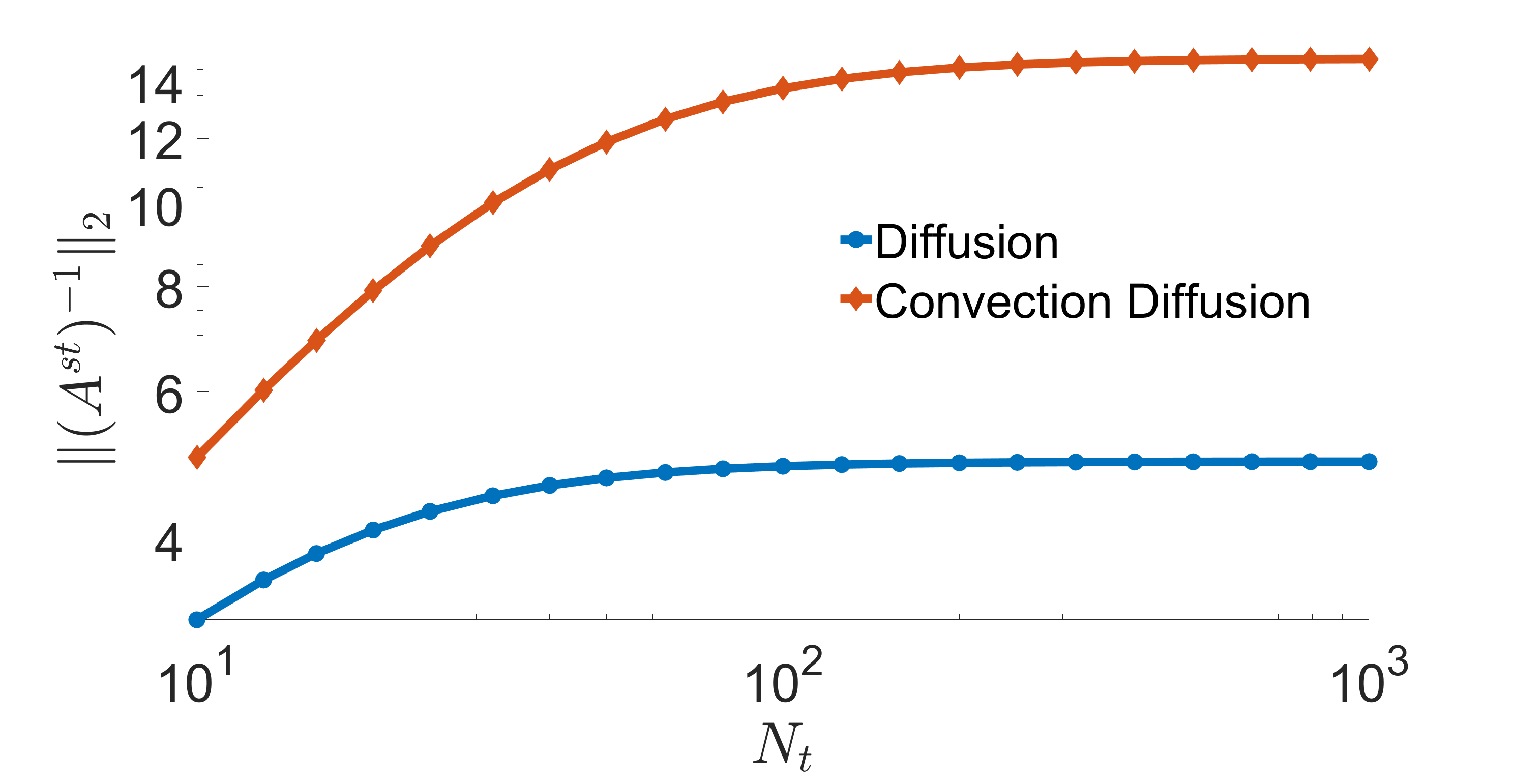}
        \caption{$\|(\ststatemat)^{-1}\|_2$ in inequality
        \eqref{eq:STerrorbound}}
    \end{subfigure}
    \begin{subfigure}[b]{0.49\textwidth}
        \centering
        \includegraphics[width=\textwidth]{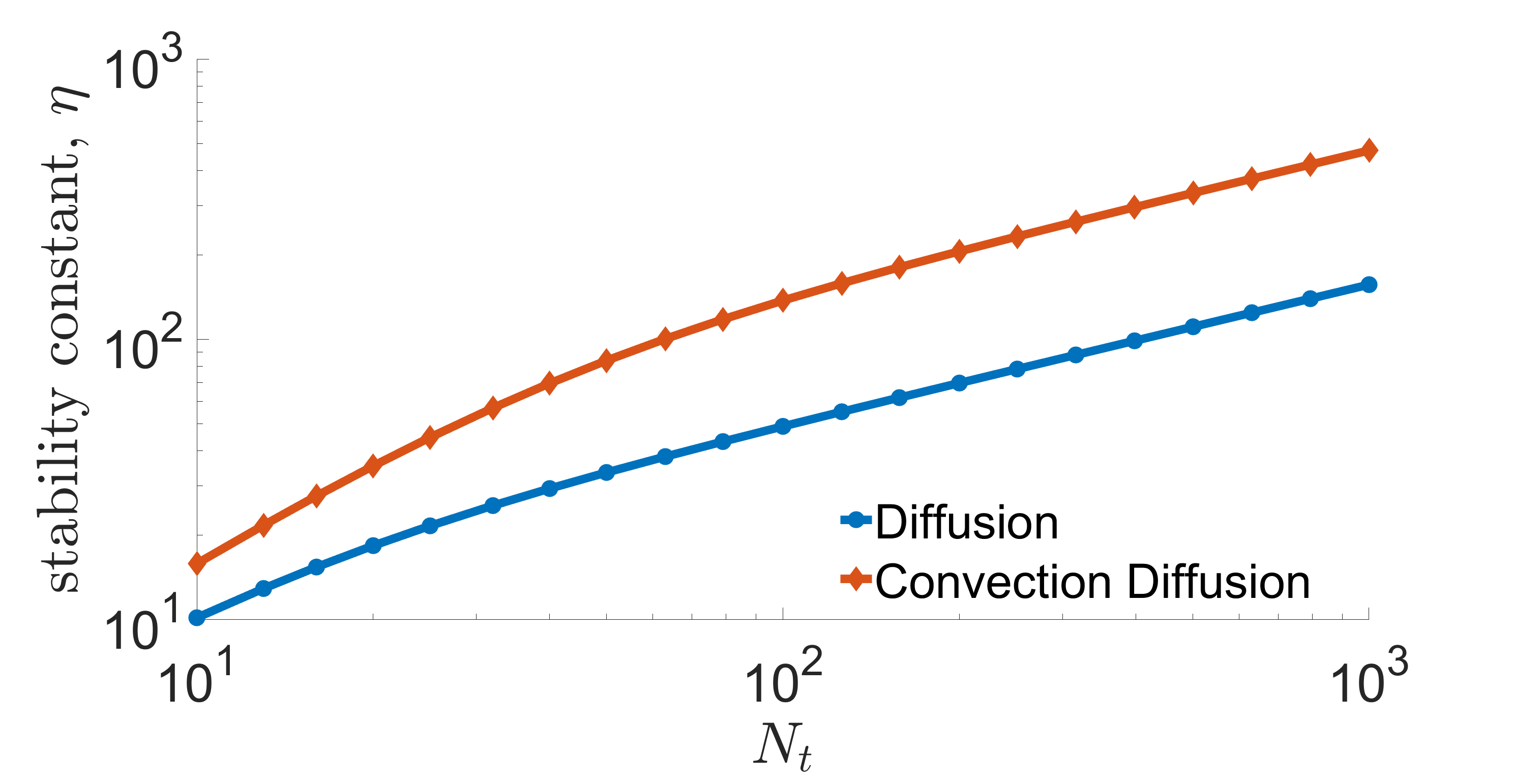}
        \caption{Stability constant, $\eta$ in inequality
        \eqref{eq:STerrorbound}}
    \end{subfigure}
    \caption{Growth rate of stability constant in
    Theorem~\ref{theorem:STerrorbound}. Backward Euler time stepping scheme with
    uniform time step size, $\timestep=10^{-2}$ is used.}
    \label{fig:fix_dt1e-2}
\end{figure}

A numerical demonstration with space--time system matrices, $\ststatemat$ that
have the same structure as the ones used in Section~\ref{sec:2Ddiffusion} and
Section~\ref{sec:2DconvectiondiffusionWOsource} shows the magnitude of
$\|(\ststatemat)^{-1}\|_2$ increases linearly for small $\ntime$, while it
becomes eventually flattened for large $\ntime$ as shown in
Fig.~\ref{fig:fix_dt1e-2}(a) for the backward Euler time integrator with uniform
time step size. Combined with $\sqrt{\ntime}$, the stability constant $\eta$
growth rate is shown in Fig.~\ref{fig:fix_dt1e-2}(b). These error bound shows
much improvement against the ones for the spatial Galerkin and Petrov--Galerkin
ROMs, which grows exponentially in time \cite{choi2020space}.

\section{Numerical results}\label{sec:NumericalResults}
In this section, we apply both the space--time Galerkin and Petrov--Galerkin ROMs
to two model problems: (i) a 2D linear diffusion equation in
Section~\ref{sec:2Ddiffusion} and (ii) a 2D linear convection-diffusion equation 
in Section~\ref{sec:2Dconvectiondiffusion}. We demonstrate their accuracy and
speed-up. The space--time ROMs are trained with solution snapshots associated with
parameters in a chosen domain and used to predict the solution of a
parameter that is not included in the trained parameter domain. We refer to this as 
the predictive case. The accuracy of space--time ROM solution
$\approxststate(\param)$ is assessed from its relative error by:
\begin{equation}\label{eq:relError}
    \relError=\frac{\left\Vert\approxststate(\param)-\ststate(\param)\right\Vert_2}{\left\Vert\ststate(\param)\right\Vert_2} 
\end{equation}
and the $\ell_2$ norm of space--time residual:
\begin{equation}\label{eq:stRes}
    \left\Vert \stres(\approxststate(\param)) \right\Vert_2.
\end{equation}
The computational cost is measured in terms of CPU wall-clock time. The online
speed-up is evaluated by dividing the wall-clock time of the FOM by the online
phase of the ROM. For the multi-query problems, total speed-up is evaluated by
dividing the time of all FOMs by the time of all ROMs including training time.
All calculations are performed on an Intel(R) Core(TM) i9-10900T CPU @ 1.90GHz
and DDR4 Memory @ 2933MHz.

\subsection{2D linear diffusion equation}\label{sec:2Ddiffusion}
We consider a parameterized 2D linear diffusion equation with a source term
\begin{equation} \label{eqn:linear_diff}
\begin{split}
\frac{\partial u}{\partial t} = \left[ \frac{\partial^2 u}{\partial x^2} +
  \frac{\partial^2 u}{\partial y^2} \right] &- \left[\frac{1}{\sqrt{(x-\mu_1)^2
  + (y-\mu_2)^2}} \right] u \\
&+ \frac{\sin{2 \pi t}}{\sqrt{(x-\mu_1)^2 + (y-\mu_2)^2}}
\end{split}
\end{equation}
where $(x,y) \in [0,1]\times [0,1]$, $t \in [0,2]$ and $(\mu_1,\mu_2) \in
[-1.7,-0.2]\times[-1.7,-0.2]$. The boundary condition is
\begin{equation}\label{eqn:linear_diff_BC} 
\begin{aligned}
u(x=0,y,t)  & = 0                     \\
u(x=1,y,t) & = 0     \\
u(x,y=0,t) & = 0   \\
u(x,y=1,t) &  = 0
\end{aligned}
\end{equation}
and the initial condition is 
\begin{equation}\label{eqn:linear_diff_IC} 
u(x,y,t=0) = 0.
\end{equation}

The backward Euler with uniform time step size $\frac{2}{\ntime}$ is employed,
where we set $\ntime=50$. For spatial differentiation, the second order central
difference scheme is implemented for the diffusion terms. Discretizing the space
domain into $N_x=70$ and $N_y=70$ uniform meshes in $x$ and $y$ directions,
respectively, gives $\nspace = (N_x-1)\times (N_y-1) =4,761$ grid points,
excluding boundary grid points. As a result, there are $238,050$ free degrees of
freedom in space--time.

For training phase, we collect solution snapshots associated with the following
parameters:
\[ (\mu_1,\mu_2)\in\{(-0.9,-0.9),(-0.9,-0.5),(-0.5,-0.9),(-0.5,-0.5)\} \] 
at which the FOM is solved.

\begin{figure}[!htbp]
  \centering
  \begin{subfigure}[b]{0.49\textwidth}
      \centering
      \includegraphics[width=\textwidth]{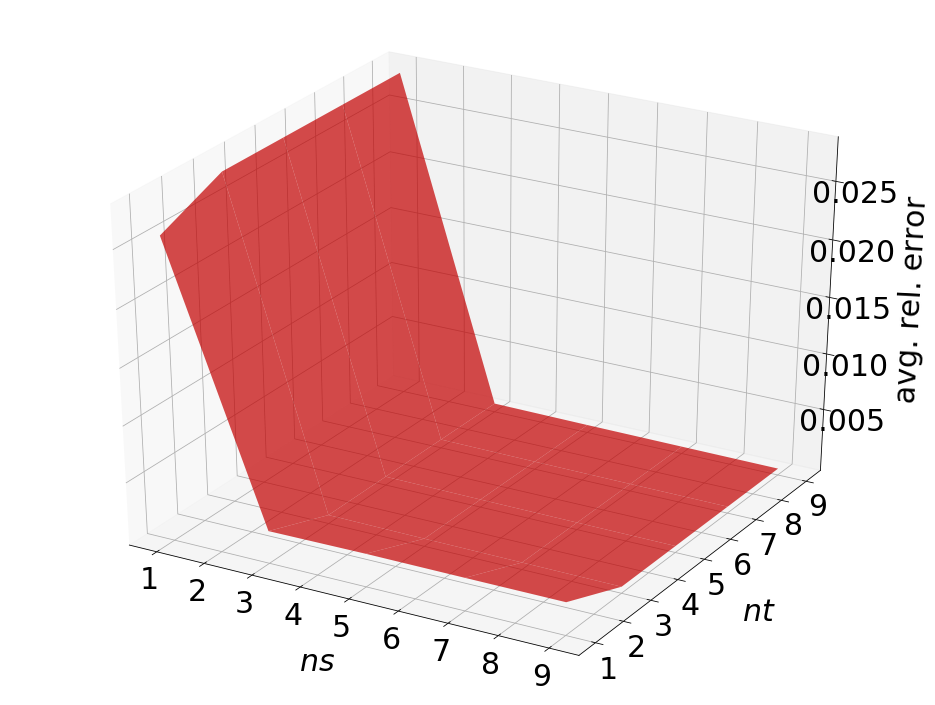}
      \caption{Relative errors vs reduced dimensions for Galerkin projection}
  \end{subfigure}
  \begin{subfigure}[b]{0.49\textwidth}
      \centering
      \includegraphics[width=\textwidth]{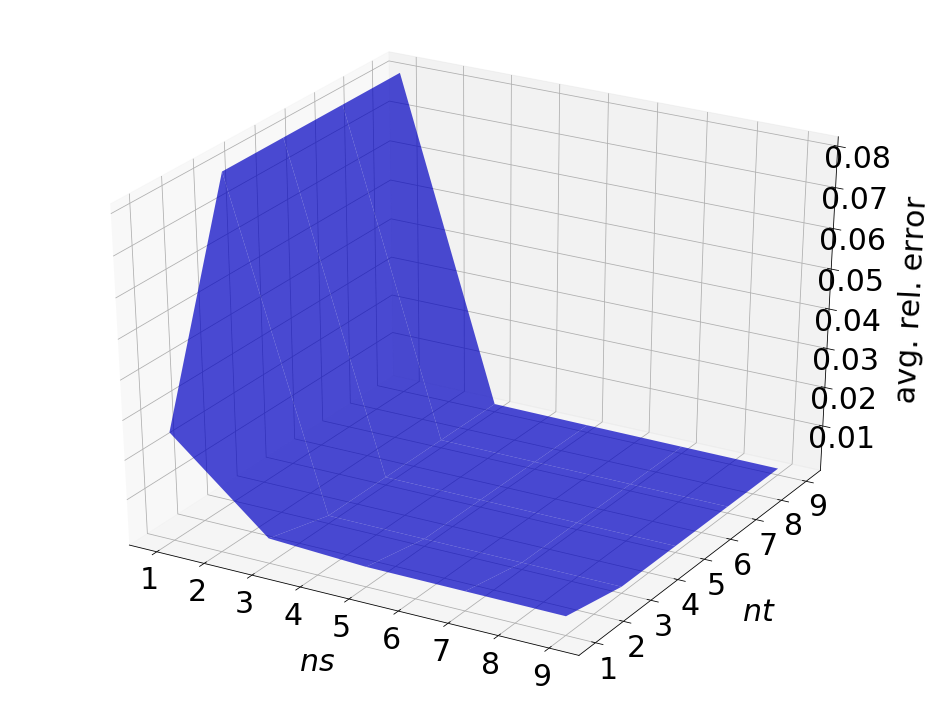}
      \caption{Relative errors vs reduced dimensions for Petrov--Galerkin
      projection}
  \end{subfigure}
  \caption{2D linear diffusion equation. Relative errors vs reduced dimensions.}
  \label{fig:DiffErrVSRed}
\end{figure}

\begin{figure}[!htbp]
  \centering
  \begin{subfigure}[b]{0.49\textwidth}
      \centering
      \includegraphics[width=\textwidth]{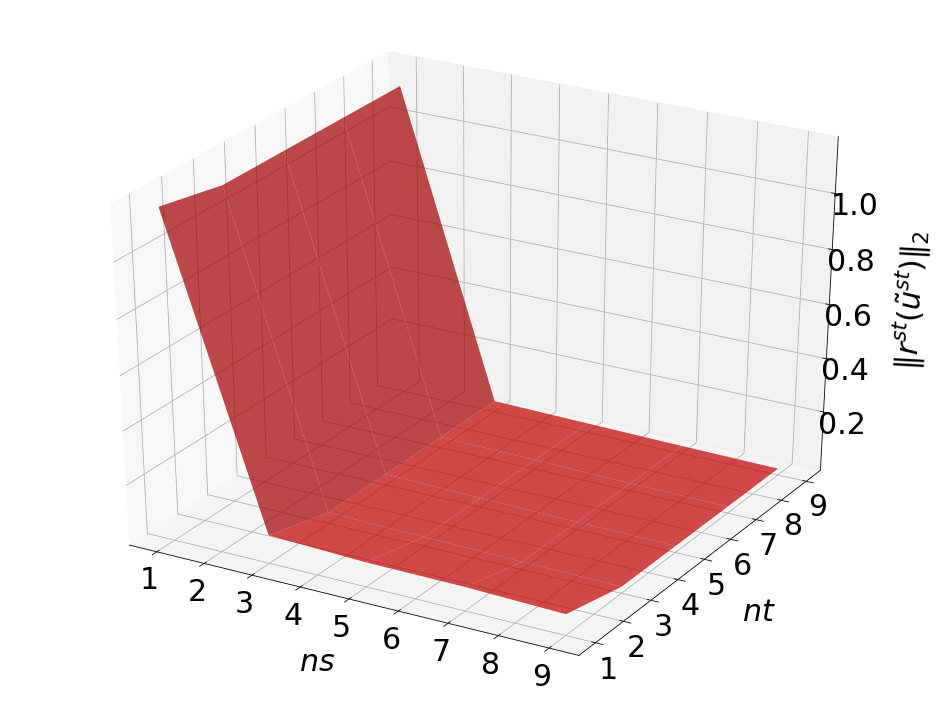}
      \caption{Space-time residuals vs reduced dimensions for Galerkin
      projection}
  \end{subfigure}
  \begin{subfigure}[b]{0.49\textwidth}
      \centering
      \includegraphics[width=\textwidth]{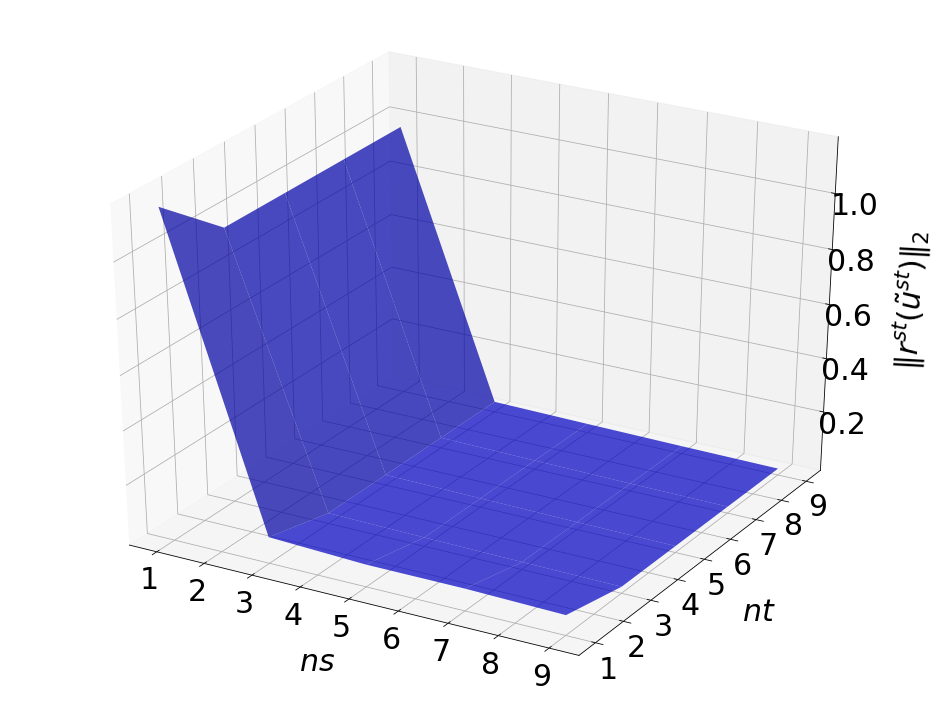}
      \caption{Space-time residuals vs reduced dimensions for Petrov--Galerkin
      projection}
  \end{subfigure}
  \caption{2D linear diffusion equation. Space-time residuals vs reduced
  dimensions.}
  \label{fig:DiffSTResVSRed}
\end{figure}

\begin{figure}[!htbp]
    \centering
    \begin{subfigure}[b]{0.49\textwidth}
        \centering
        \includegraphics[width=\textwidth]{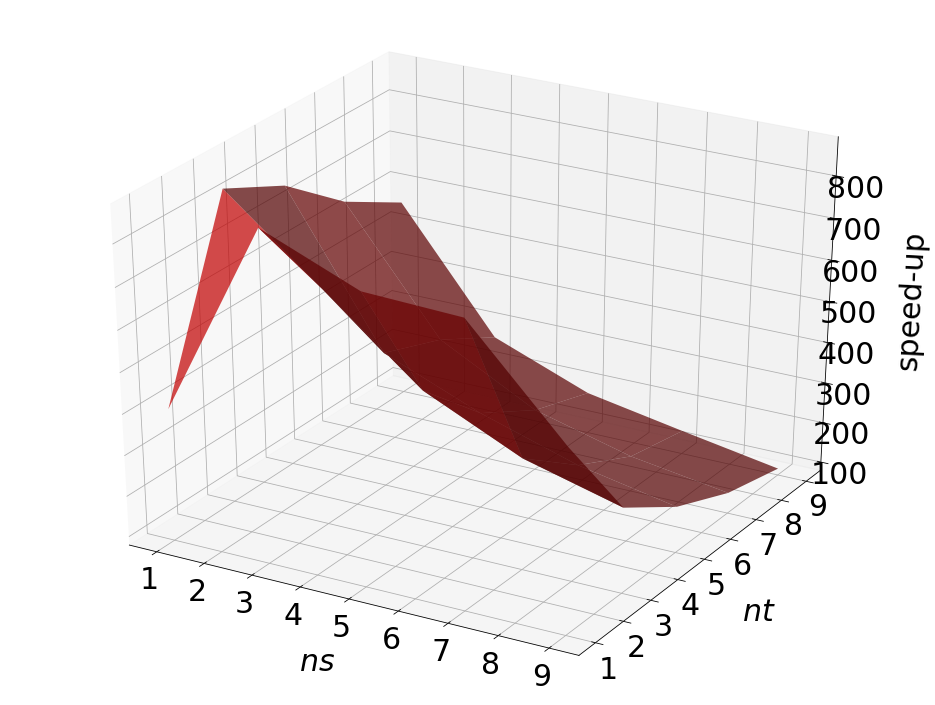}
        \caption{Speedups vs reduced dimensions for Galerkin projection}
    \end{subfigure}
    \begin{subfigure}[b]{0.49\textwidth}
        \centering
        \includegraphics[width=\textwidth]{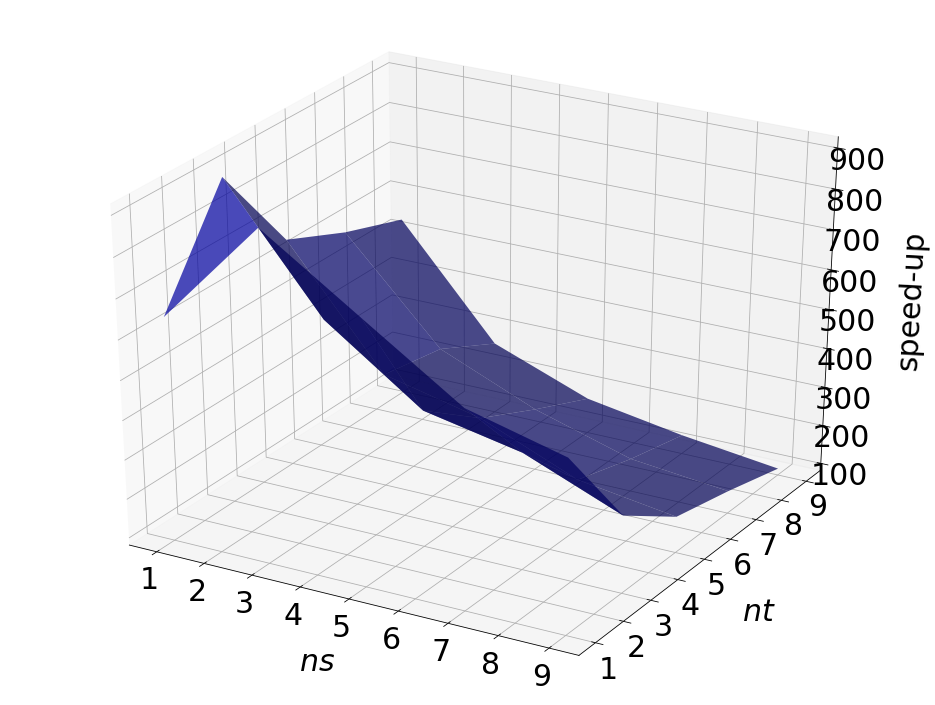}
        \caption{Speedups vs reduced dimensions for Petrov--Galerkin projection}
    \end{subfigure}
    \caption{2D linear diffusion equation. Speedups vs reduced dimensions.}
    \label{fig:DiffSpeedupVSRed}
\end{figure}

The Galerkin and Petrov--Galerkin space--time ROMs solve the
Equation~\eqref{eqn:linear_diff} with the target parameter
$(\mu_1,\mu_2)=(-0.7,-0.7))$. Fig.~\ref{fig:DiffErrVSRed},
~\ref{fig:DiffSTResVSRed}, and~\ref{fig:DiffSpeedupVSRed} show the relative
errors, the space--time residuals, and the online speed-ups as a function of the
reduced dimension $\nreducedspace$ and $\nreducedtime$. We observe that both
Galerkin and Petrov--Galerkin ROMs with $\nreducedspace=5$ and $\nreducedtime=3$
achieve a good accuracy (i.e., relative errors of 0.012\% and 0.026\%,
respectively) and speed-up (i.e., 350.31 and 376.04, respectively). We also
observe that the relative errors of Galerkin projection is smaller but the
space--time residual is larger than Petrov--Galerkin projection. This is because
Petrov--Galerkin space--time ROM solution minimizes the space--time residual.

\begin{figure}[!htbp]
    \centering
    \begin{subfigure}[b]{0.32\textwidth}
        \centering
        \includegraphics[width=\textwidth]{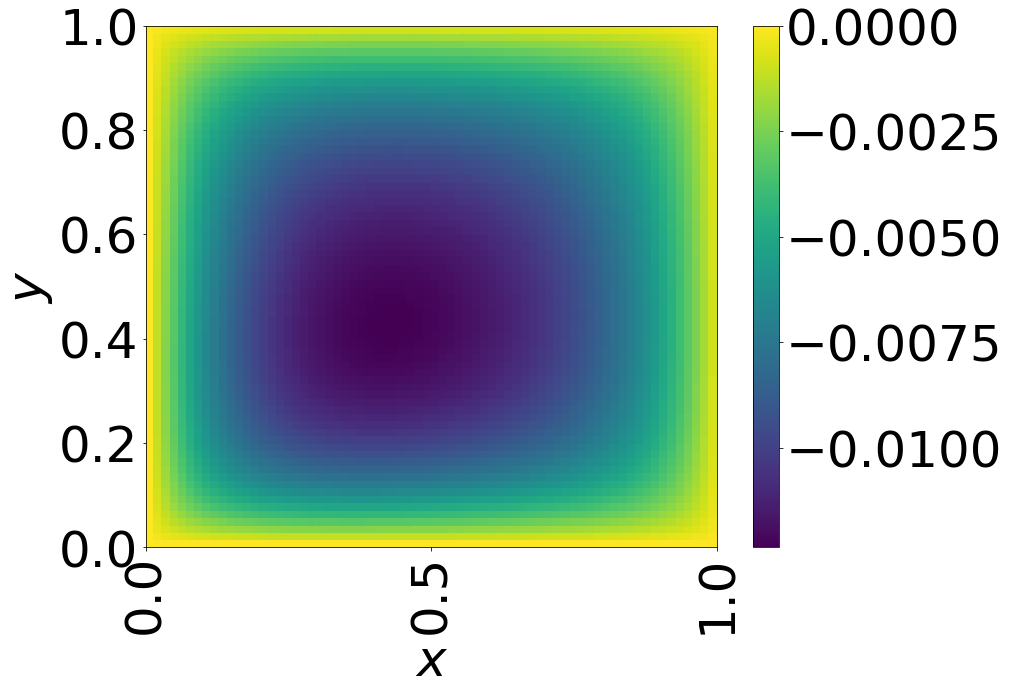}
        \caption{FOM}
    \end{subfigure}    
    \begin{subfigure}[b]{0.32\textwidth}
        \centering
        \includegraphics[width=\textwidth]{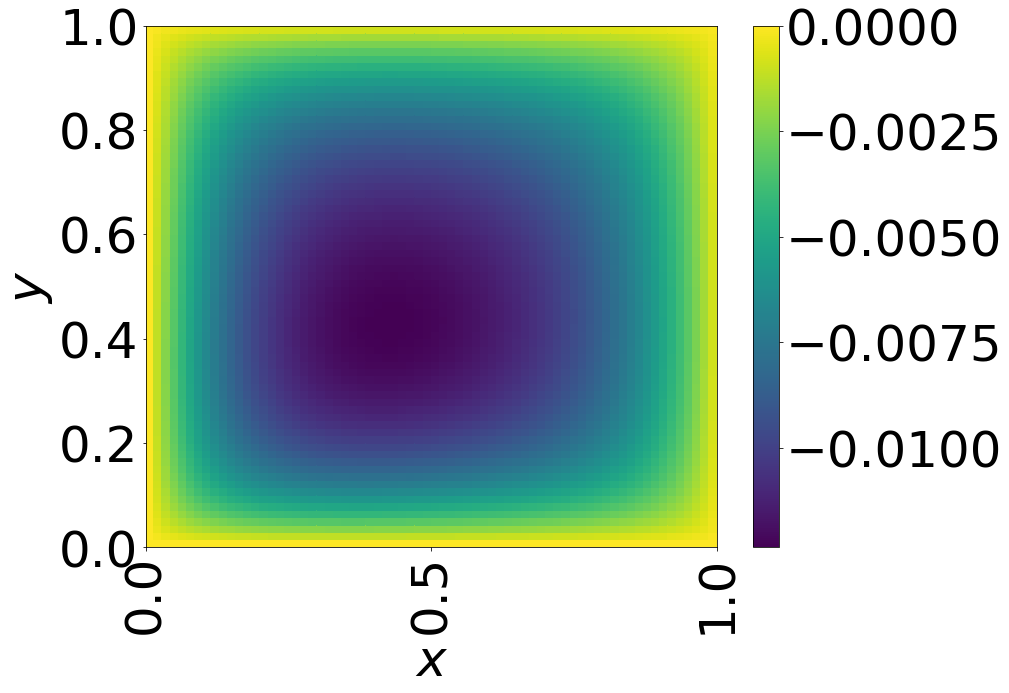}
        \caption{Galerkin ROM}
    \end{subfigure}
    \begin{subfigure}[b]{0.32\textwidth}
        \centering
        \includegraphics[width=\textwidth]{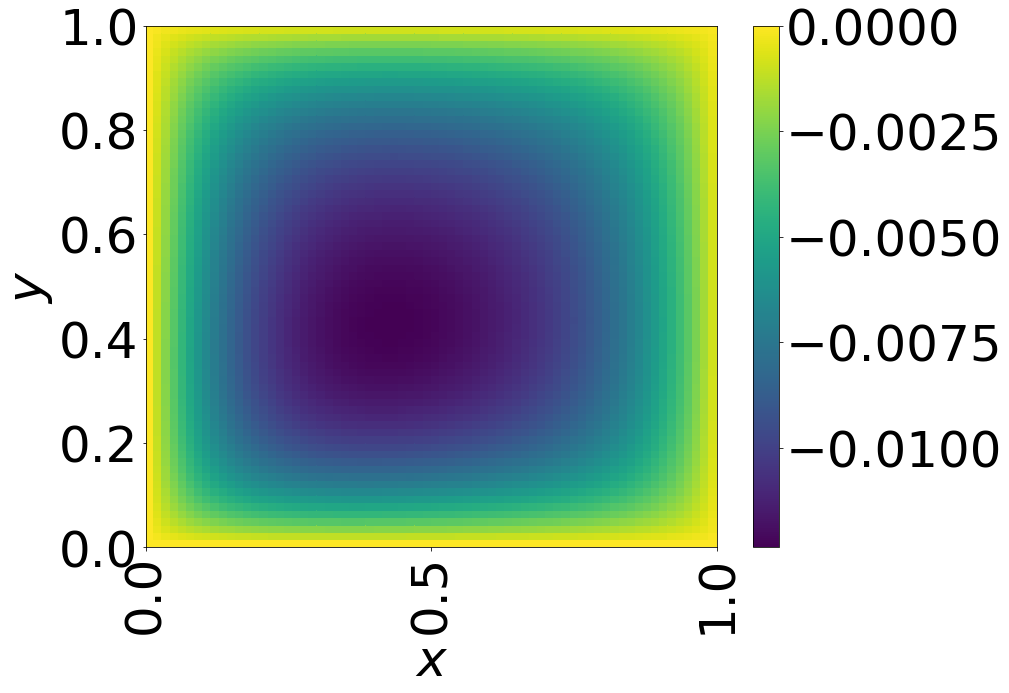}
        \caption{Petrov--Galerkin ROM}
    \end{subfigure}
    \caption{Solution snapshots of FOM, Galerkin ROM, and Petrov--Galerkin ROM at $t=2$.}
    \label{fig:DiffSol}
\end{figure}

The final time snapshots of FOM, Galerkin space--time ROM, and Petrov--Galerkin
space--time ROM are seen in Fig.~\ref{fig:DiffSol}. Both ROMs have a basis size
of $\nreducedspace=5$ and $\nreducedtime=3$, resulting in a reduction factor of
$(\nspace\ntime)\small/(\nreducedspace\nreducedtime)=15,870$. For the Galerkin
method, the FOM and space--time ROM simulation with $\nreducedspace=5$ and
$\nreducedtime=3$ takes an average time of $6.1816\times 10^{-1}$ and
$1.7646\times 10^{-3}$ seconds, respectively, resulting in speed-up of $350.31$.
For the Petrov--Galerkin method, the FOM and space--time ROM simulation with
$\nreducedspace=5$ and $\nreducedtime=3$ takes an average time of $6.0809\times
10^{-1}$ and $1.6171\times 10^{-3}$ seconds, respectively, resulting in speed-up
of $376.04$. For accuracy, the Galerkin method results in $1.210\times 10^{-2}$
\% relative error and $1.249\times 10^{-2}$ space--time residual norm while the
Petrov--Galerkin results in $2.626\times 10^{-2}$ \% relative error and
$1.029\times 10^{-2}$ space--time residual norm. 

\begin{figure}[!htbp]
    \centering
    \begin{subfigure}[b]{0.49\textwidth}
        \centering
        \includegraphics[width=\textwidth]{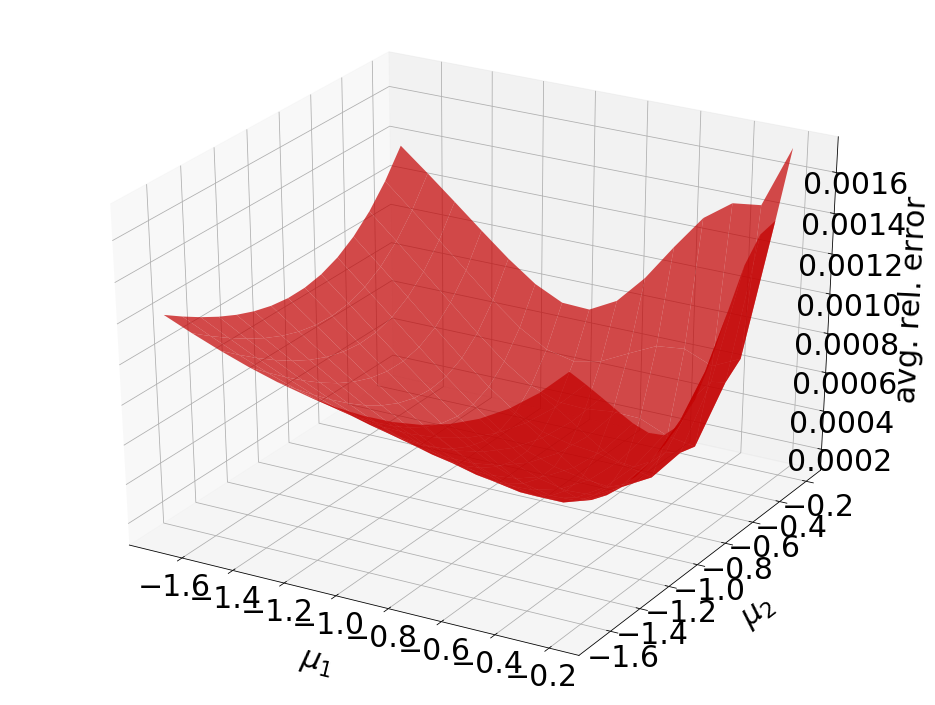}
        \caption{Glaerkin}
    \end{subfigure}    
    \begin{subfigure}[b]{0.49\textwidth}
        \centering
        \includegraphics[width=\textwidth]{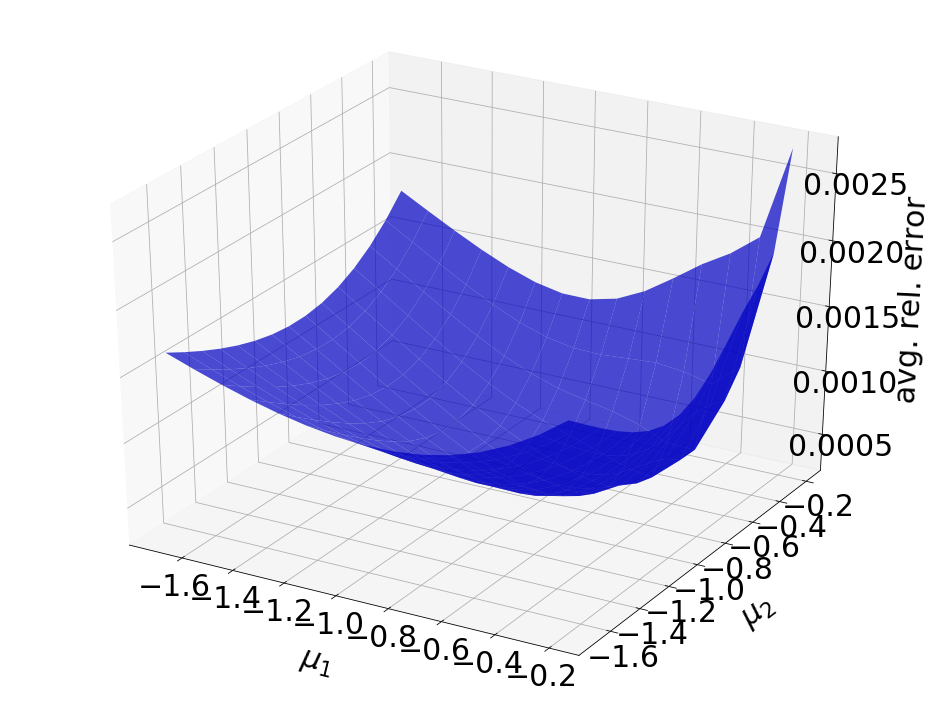}
        \caption{Petrov--Galerkin}
    \end{subfigure}
    \caption{The comparison of the Galerkin and Petrov--Galerkin ROMs for
    predictive cases}
    \label{fig:DiffErrVSParam}
\end{figure}

We investigate the numerical tests to see the generalization capability of both
Galerkin and Petrov--Galerkin ROMs. The train parameter set,
$(\mu_1,\mu_2)\in\{(-0.9,-0.9),(-0.9,-0.5),(-0.5,-0.9),(-0.5,-0.5)\}$ is used to
train a space--time ROMs with a basis of $\nreducedspace=5$ and
$\nreducedtime=3$. Then trained ROMs solve predictive cases with the test
parameter set,
$(\mu_1,\mu_2)\in\{\mu_1|\mu_1=-1.7+1.5\small/14i,i=0,1,\cdots,14\}\times
\{\mu_2|\mu_2=-1.7+1.5\small/14j,j=0,1,\cdots,14\}$.
Fig.~\ref{fig:DiffErrVSParam} shows the relative errors over the test parameter
set. The Galerkin and Petrov--Galerkin ROMs are the most accurate within the
range of the train parameter points, i.e., $[-0.9,-0.5]\times [-0.9,-0.5]$.  As
the parameter points go beyond the train parameter domain, the accuracy of the
Galerkin and Petrov--Galerkin ROMs start to deteriorate gradually. This implies
that the Galerkin and Petrov--Galerkin ROMs have a trust region. Its trust
region should be determined by the application space. For Galerkin ROM, online speed-up
is about $389$ in average and total time for ROM and FOM are $107.14$ and
$132.66$ seconds, respectively, resulting in total speed-up of $1.24$. For
Petrov--Galerkin ROM, online speed-up is about $386$ in average and total time
for ROM and FOM are $117.96$ and $132.42$ seconds, respectively, resulting in
total speed-up of $1.12$. Since the training time doesn't depend on the number
of test cases, we expect more speed-up for a larger number of test cases.

\subsection{2D linear convection diffusion equation}\label{sec:2Dconvectiondiffusion}
\subsubsection{Without source term}\label{sec:2DconvectiondiffusionWOsource} 
\begin{figure}[H]
\centering
\includegraphics[width=0.5\textwidth]{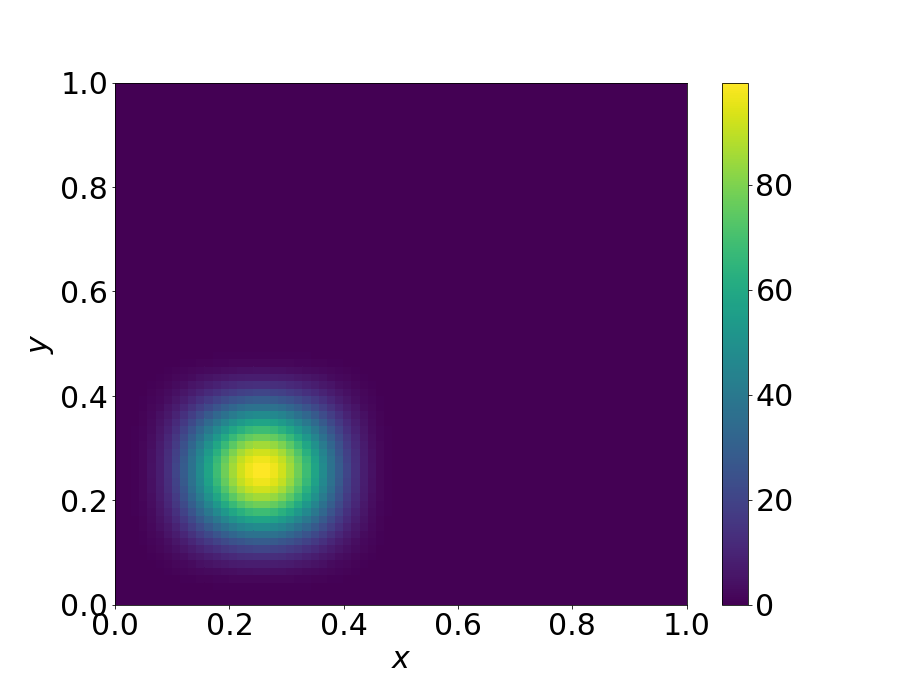}
\caption{Plot of Equation~\eqref{eqn:linear_conv_diff_IC} }
\label{fig:conv_diff_IC.png}
\end{figure}

We
consider a parameterized 2D linear convection diffusion equation
\begin{equation} \label{eqn:linear_conv_diff}
\frac{\partial u}{\partial t} = - \mu_1 \left[ \frac{\partial u}{\partial x}+
  \frac{\partial u}{\partial y} \right] + \mu_2 \left[ \frac{\partial^2
  u}{\partial x^2} + \frac{\partial^2 u}{\partial y^2} \right]
\end{equation}
where $(x,y) \in [0,1]\times [0,1]$, $t \in [0,1]$ and $(\mu_1,\mu_2) \in
[0.01,0.07]\times[0.31,0.37]$. The boundary condition is given by
\begin{equation}  \label{eqn:linear_conv_diff_BC}
\begin{aligned}
u(x=0,y,t)  & = 0                     \\
u(x=1,y,t) & = 0     \\
u(x,y=0,t) & = 0   \\
u(x,y=1,t) &  = 0.
\end{aligned}
\end{equation}
The initial condition is given by 
\begin{align}  \label{eqn:linear_conv_diff_IC}
u(x,y,t=0) = \left\{ 
\begin{array}{ll}
100\sin{(2\pi x)}^3 \cdot \sin{(2\pi y)}^3 \quad &\text{if } (x,y) \in
  [0,0.5]\times [0,0.5] \\
 0 \quad &\text{otherwise}
\end{array} \right.
\end{align}
and shown in Fig.~\ref{fig:conv_diff_IC.png}.

The backward Euler with uniform time step size $\frac{1}{\ntime}$ is employed
where we set $\ntime=50$. For spatial differentiation, a second order central
difference scheme for the diffusion terms and a first order backward difference
scheme for the convection terms are implemented. Discretizing the space domain
into $N_x=70$ and $N_y=70$ uniform meshes in $x$ and $y$ directions,
respectively, gives $\nspace = (N_x-1)\times (N_y-1) =4,761$ grid points,
excluding boundary grid points. As a result, there are $238,050$ free degrees of
freedom in space--time.

For training phase, we collect solution snapshots associated with the following
parameters:
\[(\mu_1,\mu_2) \in \{(0.03,0.33), (0.03,0.35), (0.05,0.33), (0.05,0.35)\},\] 
at which the FOM is solved.

\begin{figure}[!htbp]
    \centering
    \begin{subfigure}[b]{0.49\textwidth}
        \centering
        \includegraphics[width=\textwidth]{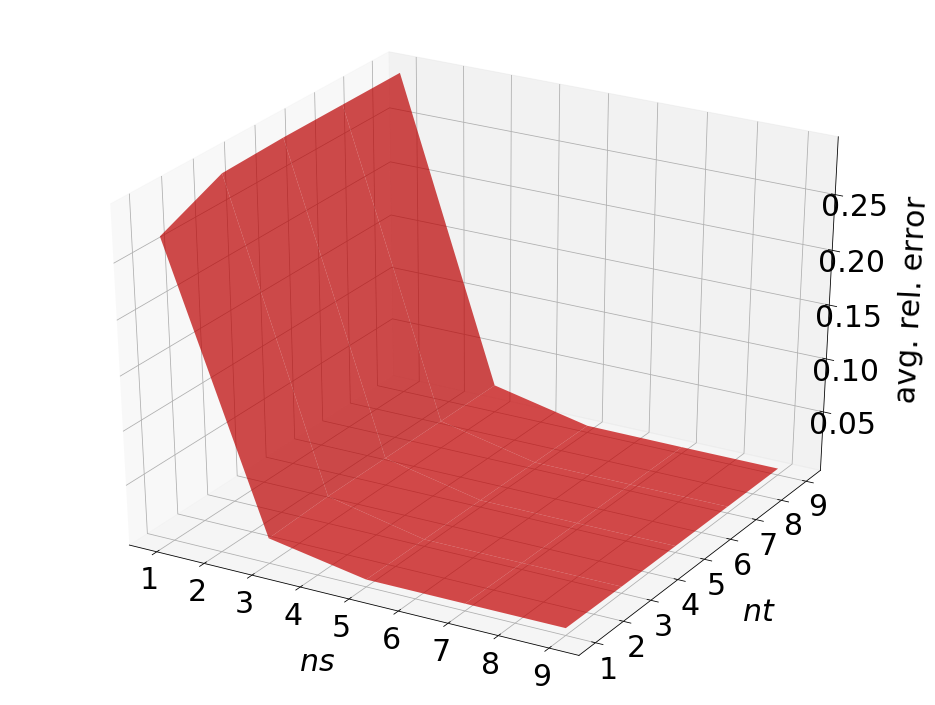}
        \caption{Relative errors vs reduced dimensions for Galerkin projection}
    \end{subfigure}
    \begin{subfigure}[b]{0.49\textwidth}
        \centering
        \includegraphics[width=\textwidth]{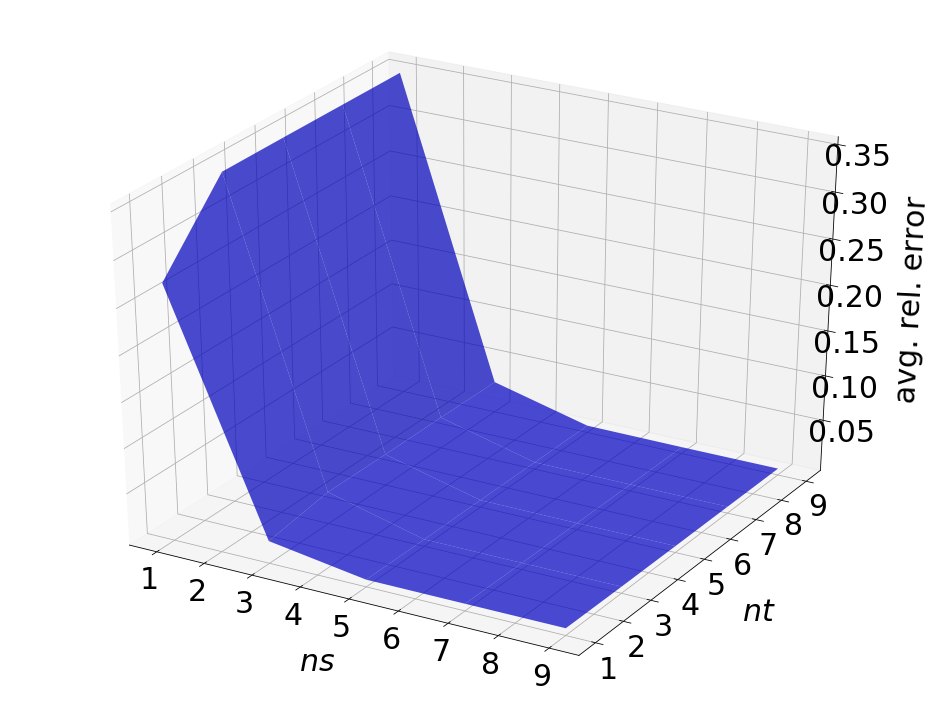}
        \caption{Relative errors vs reduced dimensions for Petrov--Galerkin
        projection}
    \end{subfigure}
    \caption{2D linear convection diffusion equation. Relative errors vs reduced
    dimensions.}
    \label{fig:ConvDiffErrVSRed}
\end{figure}

\begin{figure}[!htbp]
    \centering
    \begin{subfigure}[b]{0.49\textwidth}
        \centering
        \includegraphics[width=\textwidth]{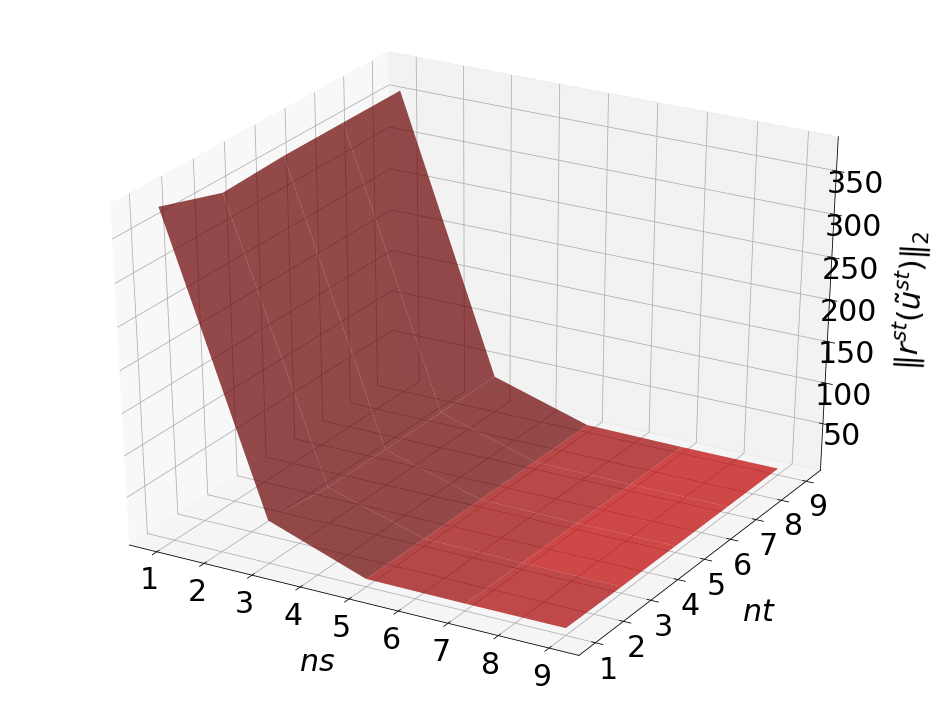}
        \caption{Space-time residuals vs reduced dimensions for Galerkin
        projection}
    \end{subfigure}
    \begin{subfigure}[b]{0.49\textwidth}
        \centering
        \includegraphics[width=\textwidth]{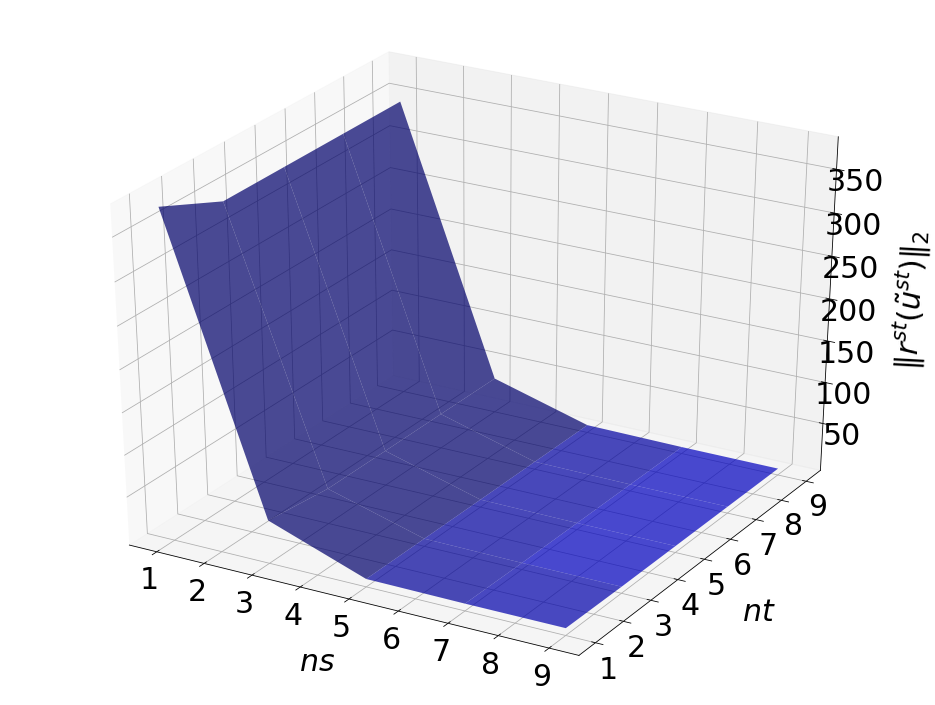}
        \caption{Space-time residuals vs reduced dimensions for Petrov--Galerkin
        projection}
    \end{subfigure}
    \caption{2D linear convection diffusion equation. Space-time residuals vs
    reduced dimensions.}
    \label{fig:ConvDiffSTResVSRed}
\end{figure}

\begin{figure}[!htbp]
    \centering
    \begin{subfigure}[b]{0.49\textwidth}
        \centering
        \includegraphics[width=\textwidth]{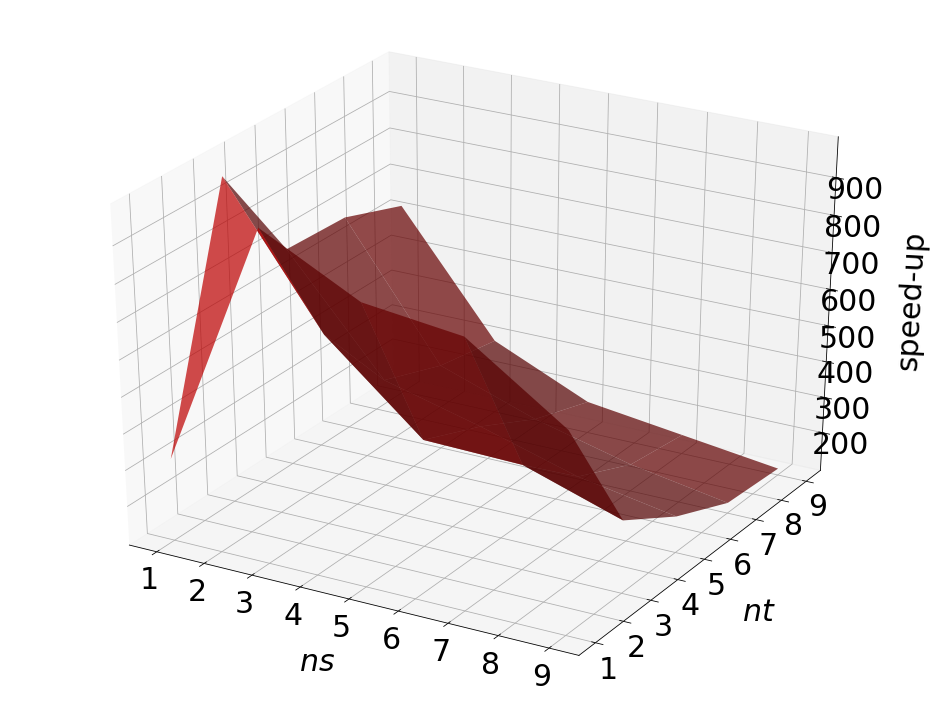}
        \caption{Speedups vs reduced dimensions for Galerkin projection}
    \end{subfigure}
    \begin{subfigure}[b]{0.49\textwidth}
        \centering
        \includegraphics[width=\textwidth]{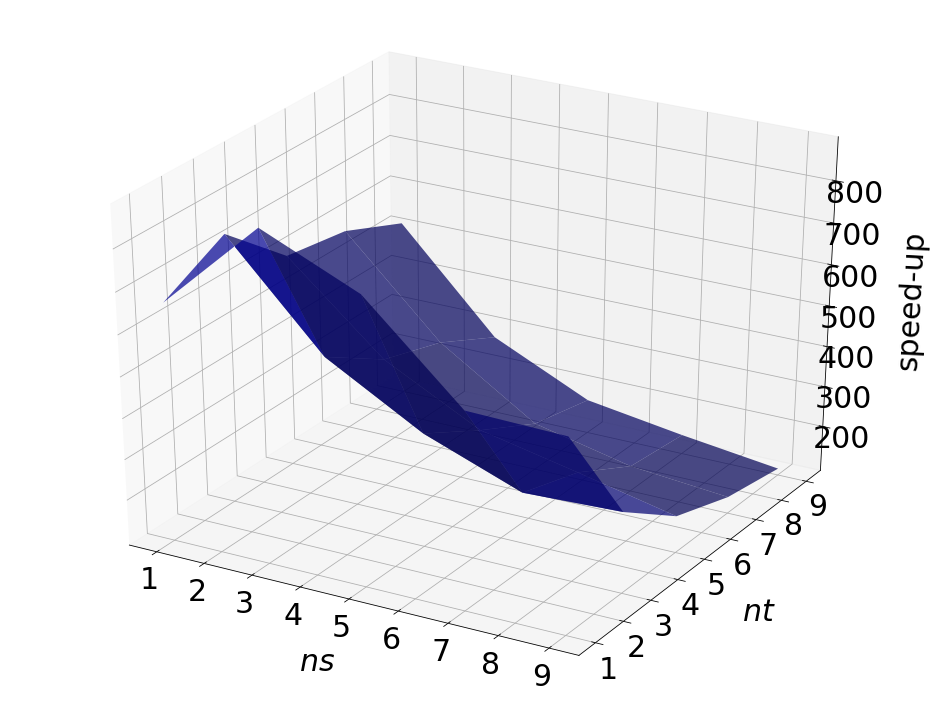}
        \caption{Speedups vs reduced dimensions for Petrov--Galerkin projection}
    \end{subfigure}
    \caption{2D linear convection diffusion equation. Speedups vs reduced
    dimensions.}
    \label{fig:ConvDiffSpeedupVSRed}
\end{figure}

The Galerkin and Petrov--Galerkin space--time ROMs solve the
Equation~\eqref{eqn:linear_conv_diff} with the target parameter
$(\mu_1,\mu_2)=(0.04,0.34)$. Fig.~\ref{fig:ConvDiffErrVSRed},
~\ref{fig:ConvDiffSTResVSRed}, and~\ref{fig:ConvDiffSpeedupVSRed} show the
relative errors, the space--time residuals, and the online speed-ups as a
function of the reduced dimension $\nreducedspace$ and $\nreducedtime$. We
observe that both Galerkin and Petrov--Galerkin ROMs with $\nreducedspace=5$ and
$\nreducedtime=3$ achieve a good accuracy (i.e., relative errors of 0.049\% and
0.059\%, respectively) and speed-up (i.e., 451.17 and 370.74, respectively). We
also observe that the relative errors of Galerkin projection is smaller but the
space--time residual is larger than Petrov--Galerkin projection. This is because
Petrov--Galerkin space--time ROM solution minimizes the space--time residual.

\begin{figure}[!htbp]
    \centering
    \begin{subfigure}[b]{0.32\textwidth}
        \centering
        \includegraphics[width=\textwidth]{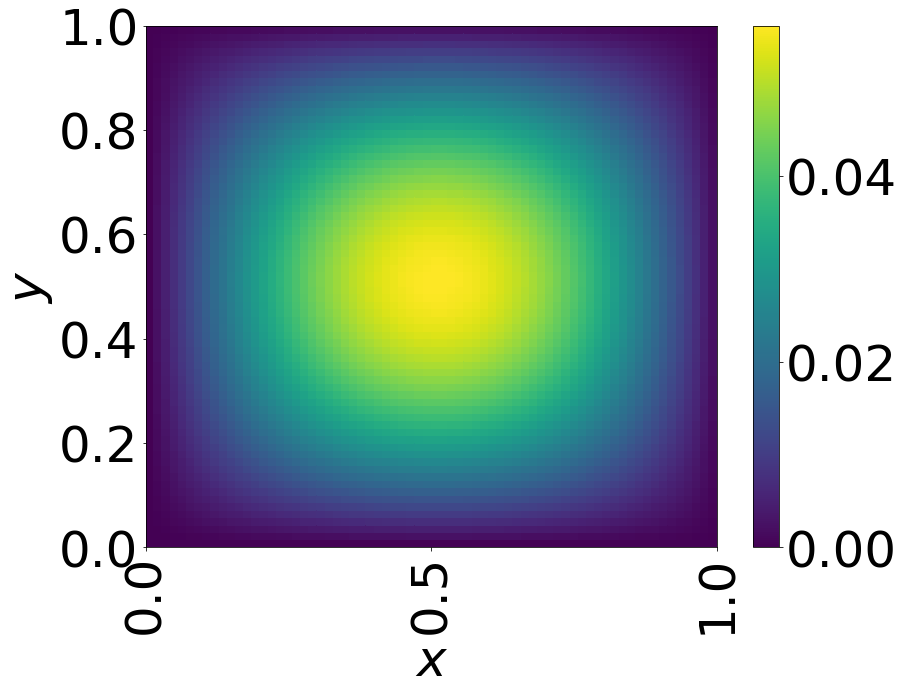}
        \caption{FOM}
    \end{subfigure}    
    \begin{subfigure}[b]{0.32\textwidth}
        \centering
        \includegraphics[width=\textwidth]{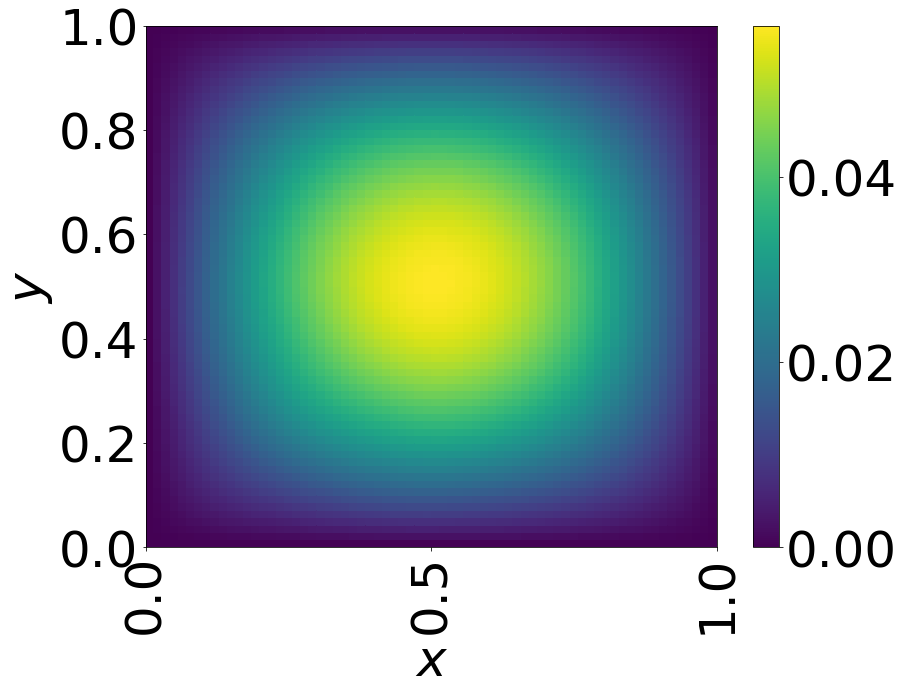}
        \caption{Galerkin ROM}
    \end{subfigure}
    \begin{subfigure}[b]{0.32\textwidth}
        \centering
        \includegraphics[width=\textwidth]{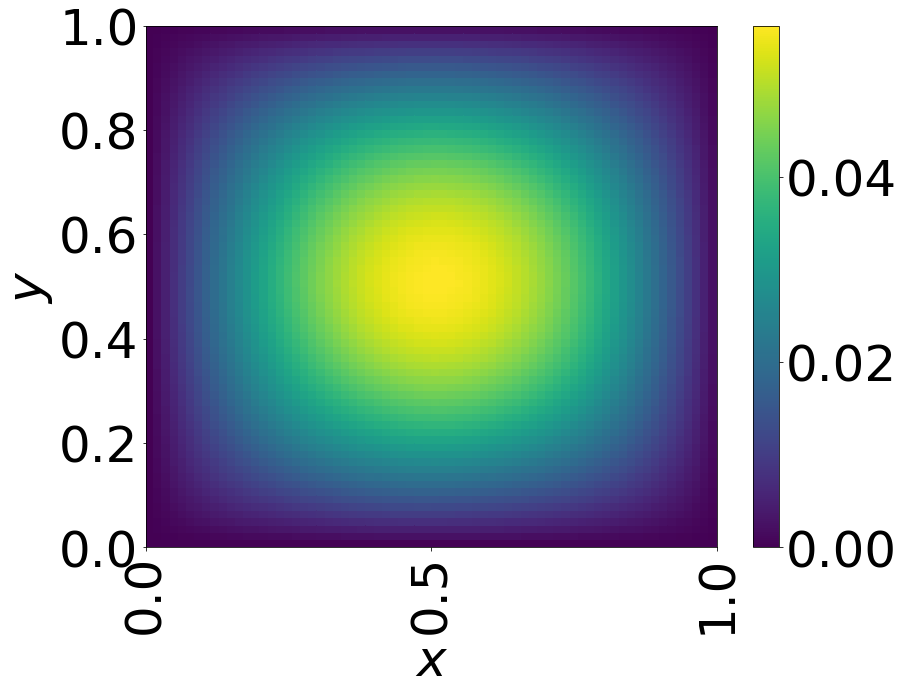}
        \caption{Petrov--Galerkin ROM}
    \end{subfigure}
    \caption{Solution snapshots of FOM, Galerkin ROM, and Petrov--Galerkin ROM
    at $t=1$.}
    \label{fig:ConvDiffSol}
\end{figure}

The final time snapshots of FOM, Galerkin space--time ROM, and Petrov--Galerkin
space--time ROM are seen in Fig.~\ref{fig:ConvDiffSol}. Both ROMs have a basis
size of $\nreducedspace=5$ and $\nreducedtime=3$, resulting in a reduction
factor of $(\nspace\ntime)\small/(\nreducedspace\nreducedtime)=15,870$. For the
Galerkin method, the FOM and space--time ROM simulation with $\nreducedspace=5$
and $\nreducedtime=3$ takes an average time of $6.1562\times 10^{-1}$ and
$1.3645\times 10^{-3}$ seconds, respectively, resulting in speed-up of $451.17$.
For the Petrov--Galerkin method, the FOM and space--time ROM simulation with
$\nreducedspace=5$ and $\nreducedtime=3$ takes an average time of
$5.7617\times 10^{-1}$ and $1.5541\times 10^{-3}$ seconds, respectively, resulting
in speed-up of $370.74$. For accuracy, the Galerkin method results in
$4.898\times 10^{-2}$ \% relative error and $1.503$ space--time residual norm while
the Petrov--Galerkin results in $5.878\times 10^{-2}$ \% relative error and
$1.459$ space--time residual norm. 

\begin{figure}[!htbp]
    \centering
    \begin{subfigure}[b]{0.49\textwidth}
        \centering
        \includegraphics[width=\textwidth]{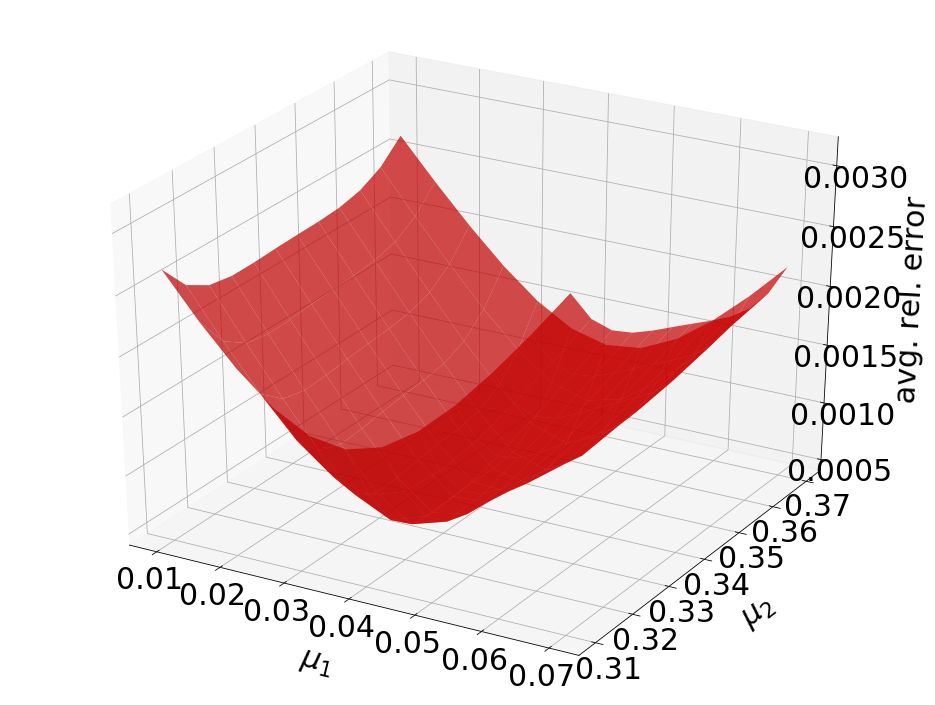}
        \caption{Glaerkin}
    \end{subfigure}    
    \begin{subfigure}[b]{0.49\textwidth}
        \centering
        \includegraphics[width=\textwidth]{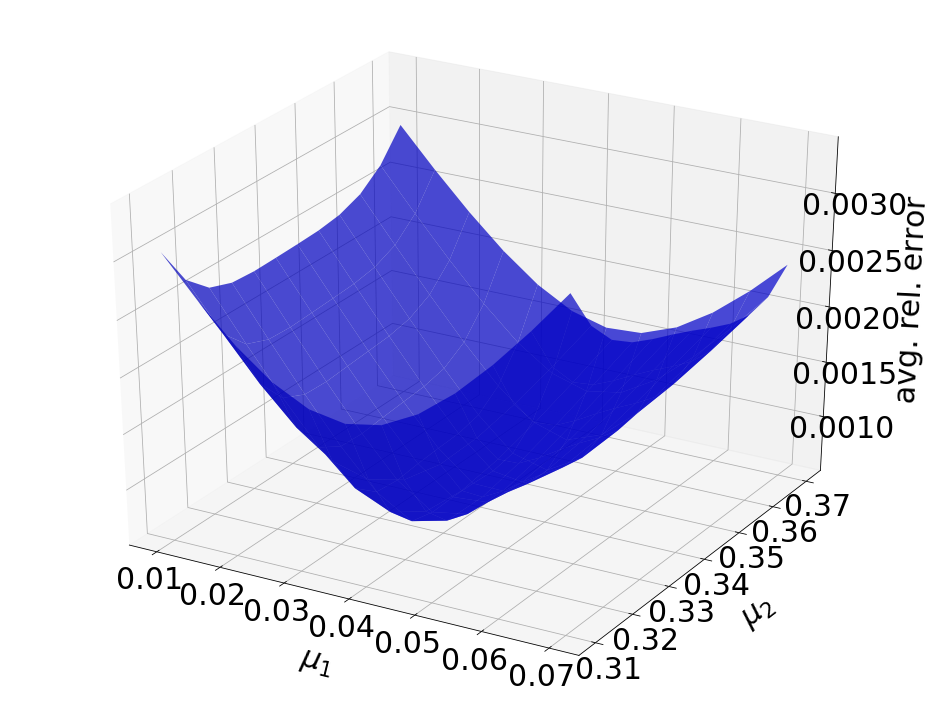}
        \caption{Petrov--Galerkin}
    \end{subfigure}
    \caption{The comparison of the Galerkin and Petrov--Galerkin ROMs for
    predictive cases}
    \label{fig:ConvDiffErrVSParam}
\end{figure}

We investigate the numerical tests to see the generalization capability of both
Galerkin and Petrov--Galerkin ROMs. The train parameter set,
$(\mu_1,\mu_2)\in\{(0.03,0.33),(0.03,0.35),(0.05,0.33),(0.05,0.35)\}$ is used to
train a space--time ROMs with a basis of $\nreducedspace=5$ and
$\nreducedtime=3$. Then trained ROMs solve predictive cases with the test
parameter set,
$(\mu_1,\mu_2)\in\{\mu_1|\mu_1=0.01+0.06\small/11i,i=0,1,\cdots,11\}\times
\{\mu_2|\mu_2=0.31+0.06\small/11j,j=0,1,\cdots,11\}$.
Fig.~\ref{fig:ConvDiffErrVSParam} shows the relative errors over the test
parameter set. The Galerkin and Petrov--Galerkin ROMs are the
most accurate within the range of the train parameter points, i.e.,
$[0.03,0.33]\times [0.05,0.35]$.  As the parameter points go beyond the train
parameter domain, the accuracy of the Galerkin and Petrov--Galerkin ROMs start
to deteriorate gradually. This implies that the Galerkin and Petrov--Galerkin
ROMs have a trust region. Its trust region should be determined by an
application. For Galerkin ROM, online speed-up is about $387$ in average and
total time for ROM and FOM are $65.03$ and $83.89$ seconds, respectively,
resulting in total speed-up of $1.29$. For Petrov--Galerkin ROM, online speed-up
is about $385$ in average and total time for ROM and FOM are $70.55$ and $83.34$
seconds, respectively, resulting in total speed-up of $1.18$. Since the training
time doesn't depend on the number of test cases, we expect more speed-up for the
larger number of test cases.

\subsubsection{With source term}
We consider a parameterized 2D linear convection diffusion equation
\begin{equation} \label{eqn:linear_conv_diff_source}
\frac{\partial u}{\partial t} = - \mu_1 \left[0.1 \frac{\partial u}{\partial x}+
  \frac{\partial u}{\partial y} \right] + \mu_2 \left[ \frac{\partial^2
  u}{\partial x^2} + \frac{\partial^2 u}{\partial y^2} \right] + f(x,y,t)
\end{equation}
with the source term $f(x,y,t)$ which is given by
\begin{equation} \label{eqn:sourceterm_linear_conv_diff}
f(x,y,t) = 10^{5}\exp(-( (\frac{x-0.5+0.2\sin{(2\pi
  t})}{0.1})^2+(\frac{y}{0.05})^2 ) )
\end{equation}
where $(x,y) \in [0,1]\times [0,1]$, $t \in [0,2]$ and $(\mu_1,\mu_2) \in
[0.195,0.205]\times[0.018,0.022]$. The boundary condition is given by
\begin{equation}  \label{eqn:linear_conv_diff_source_BC}
\begin{aligned}
u(x=0,y,t)  & = 0                     \\
u(x=1,y,t) & = 0     \\
u(x,y=0,t) & = 0   \\
u(x,y=1,t) &  = 0
\end{aligned}
\end{equation}
and the initial condition is given by 
\begin{equation}  \label{eqn:linear_conv_diff_source_IC}
u(x,y,t=0) = 0.
\end{equation}

The backward Euler with uniform time step size $\frac{2}{\ntime}$ is employed
where we set $\ntime=50$. For spatial differentiation, a second order central
difference scheme for the diffusion terms and a first order backward difference
scheme for the convection terms are implemented. Discretizing the space domain
into $N_x=70$ and $N_y=70$ uniform meshes in $x$ and $y$ directions,
respectively, gives $\nspace = (N_x-1)\times (N_y-1) = 4,761$ grid points,
excluding boundary grid points. As a result, there are $238,050$ free degrees of
freedom in space--time.

For training phase, we collect solution snapshots associated with the following
parameters:
\[ (\mu_1,\mu_2) \in \{(0.195,0.018), (0.195,0.022), (0.205,0.018),
(0.205,0.022)\},\] 
at which the FOM is solved.

\begin{figure}[!htbp]
  \centering
  \begin{subfigure}[b]{0.49\textwidth}
      \centering
      \includegraphics[width=\textwidth]{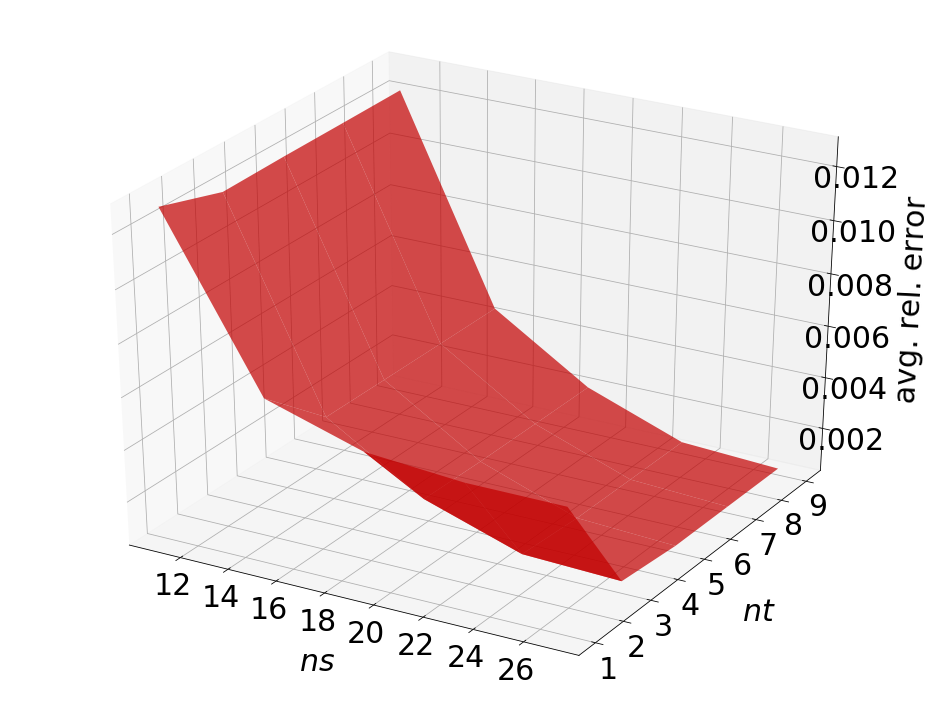}
      \caption{Relative errors vs reduced dimensions for Galerkin projection}
  \end{subfigure}
  \begin{subfigure}[b]{0.49\textwidth}
      \centering
      \includegraphics[width=\textwidth]{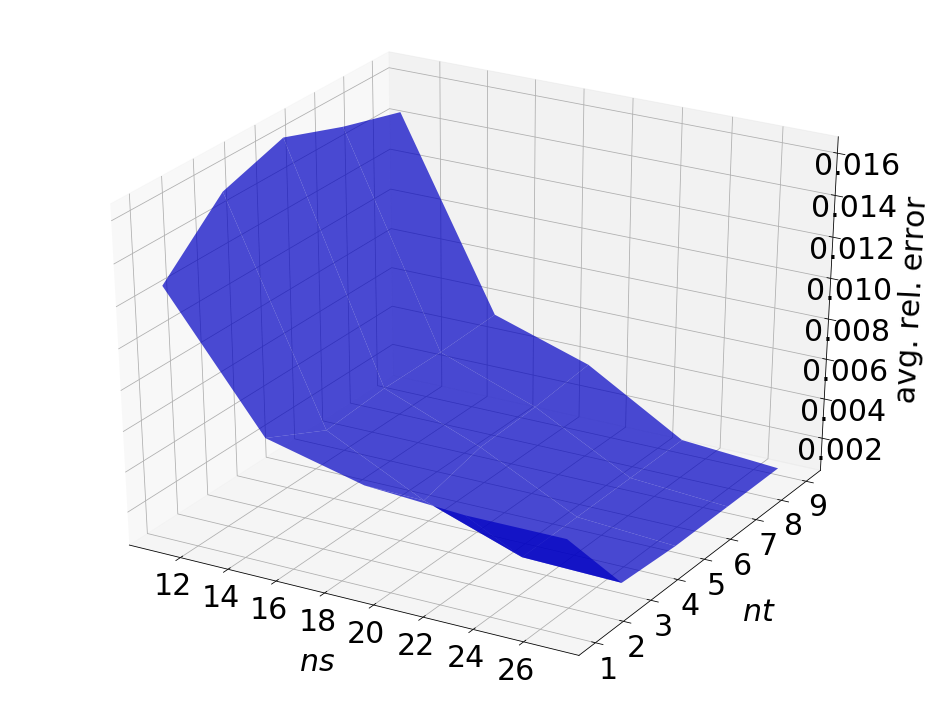}
      \caption{Relative errors vs reduced dimensions for Petrov--Galerkin
      projection}
  \end{subfigure}
  \caption{2D linear convection diffusion equation with source term. Relative
  errors vs reduced dimensions.}
  \label{fig:ConvDiffSourceErrVSRed}
\end{figure}

\begin{figure}[!htbp]
  \centering
  \begin{subfigure}[b]{0.49\textwidth}
      \centering
      \includegraphics[width=\textwidth]{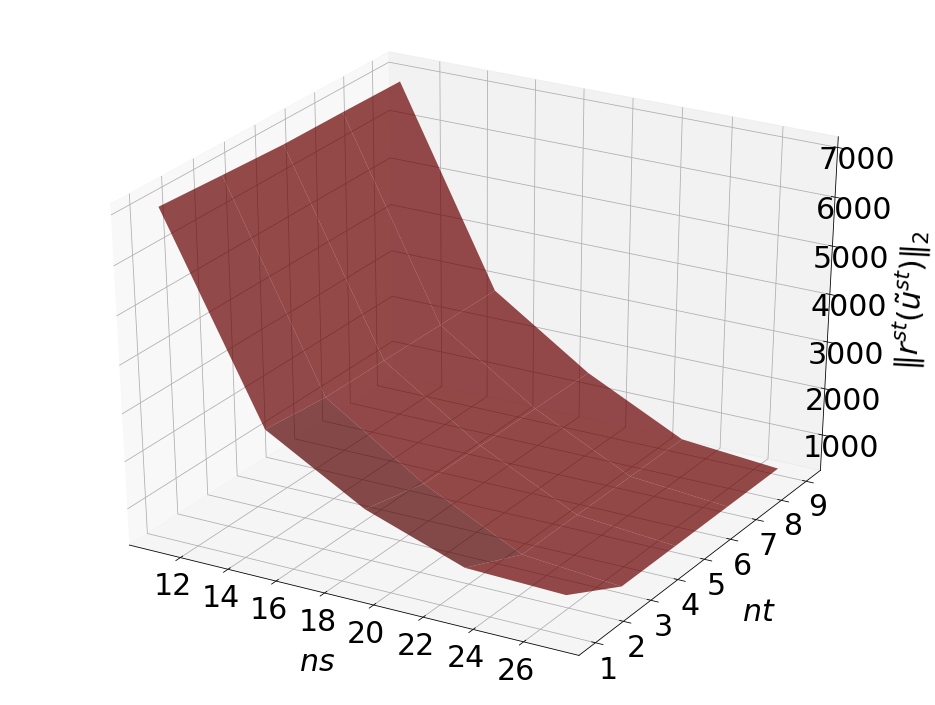}
      \caption{Space-time residuals vs reduced dimensions for Galerkin
      projection}
  \end{subfigure}
  \begin{subfigure}[b]{0.49\textwidth}
      \centering
      \includegraphics[width=\textwidth]{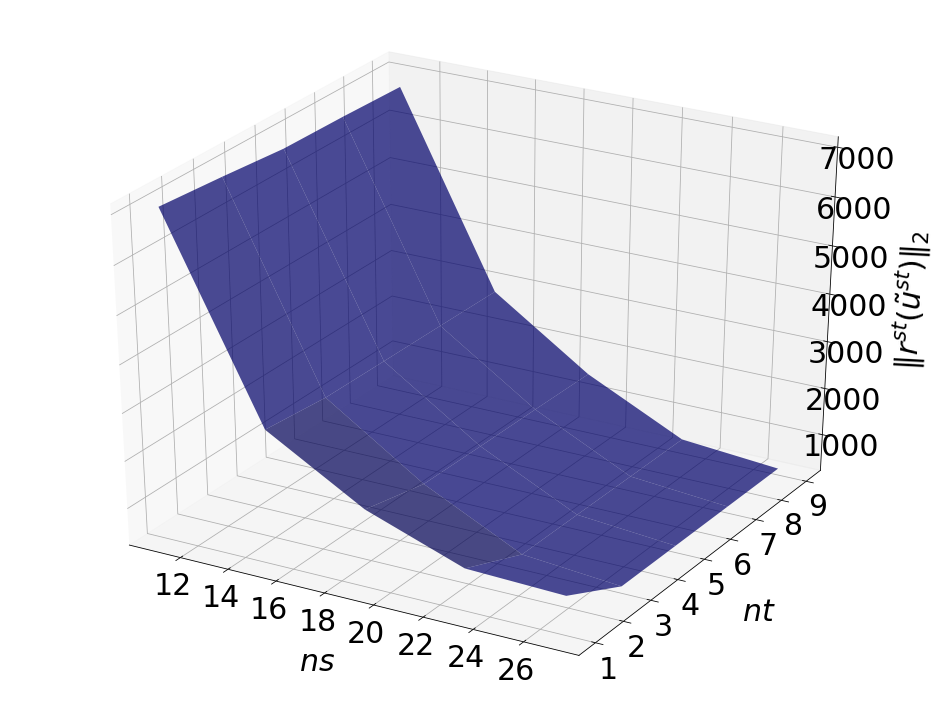}
      \caption{Space-time residuals vs reduced dimensions for Petrov--Galerkin
      projection}
  \end{subfigure}
  \caption{2D linear convection diffusion equation with source term. Space-time
  residuals vs reduced dimensions.}
  \label{fig:ConvDiffSourceSTResVSRed}
\end{figure}

\begin{figure}[!htbp]
  \centering
  \begin{subfigure}[b]{0.49\textwidth}
    \centering
    \includegraphics[width=\textwidth]{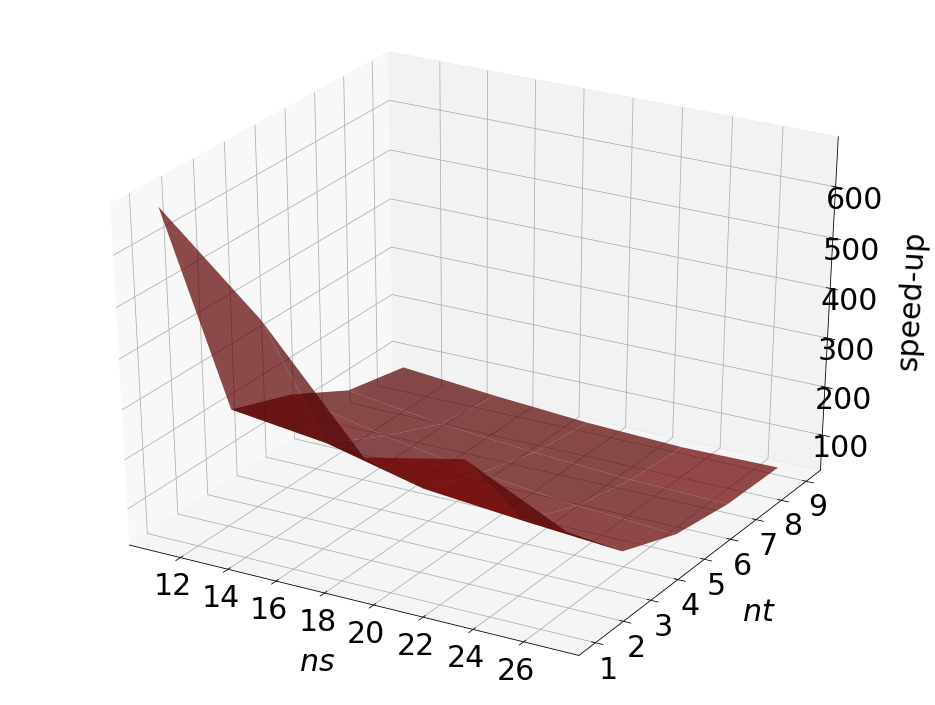}
    \caption{Speedups vs reduced dimensions for Galerkin projection}
  \end{subfigure}
  \begin{subfigure}[b]{0.49\textwidth}
    \centering
    \includegraphics[width=\textwidth]{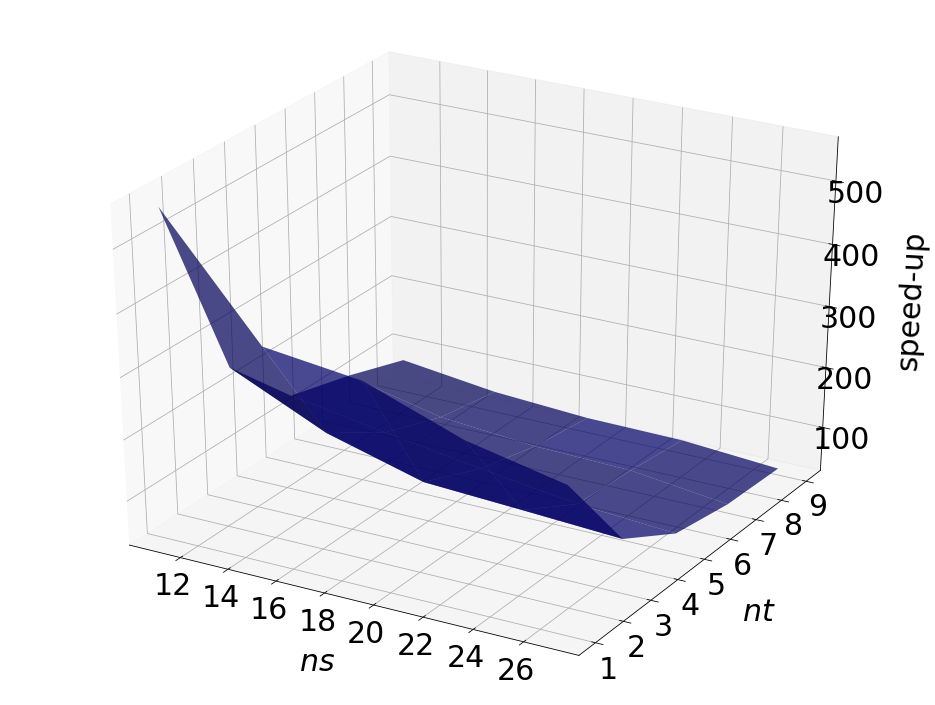}
    \caption{Speedups vs reduced dimensions for Petrov--Galerkin projection}
  \end{subfigure}
  \caption{2D linear convection diffusion equation with source term. Speedups vs
  reduced dimensions.}
  \label{fig:ConvDiffSourceSpeedupVSRed}
\end{figure}

The Galerkin and Petrov--Galerkin space--time ROMs solve the
Equation~\eqref{eqn:linear_conv_diff_source} with the target parameter
$(\mu_1,\mu_2)=(0.2,0.02))$. Fig.~\ref{fig:ConvDiffSourceErrVSRed},
~\ref{fig:ConvDiffSourceSTResVSRed}, and~\ref{fig:ConvDiffSourceSpeedupVSRed}
show the relative errors, the space--time residuals, and the online speed-ups as
a function of the reduced dimension $\nreducedspace$ and $\nreducedtime$. We
observe that both Galerkin and Petrov--Galerkin ROMs with $\nreducedspace=19$
and $\nreducedtime=3$ achieve a good accuracy (i.e., relative errors of 0.217\%
and 0.265\%, respectively) and speed-up (i.e., 153.87 and 139.88, respectively).
We also observe that the relative errors of Galerkin projection is smaller but
the space--time residual is larger than Petrov--Galerkin projection. This is
because Petrov--Galerkin space--time ROM solution minimizes the space--time
residual.

\begin{figure}[!htbp]
    \centering
    \begin{subfigure}[b]{0.32\textwidth}
        \centering
        \includegraphics[width=\textwidth]{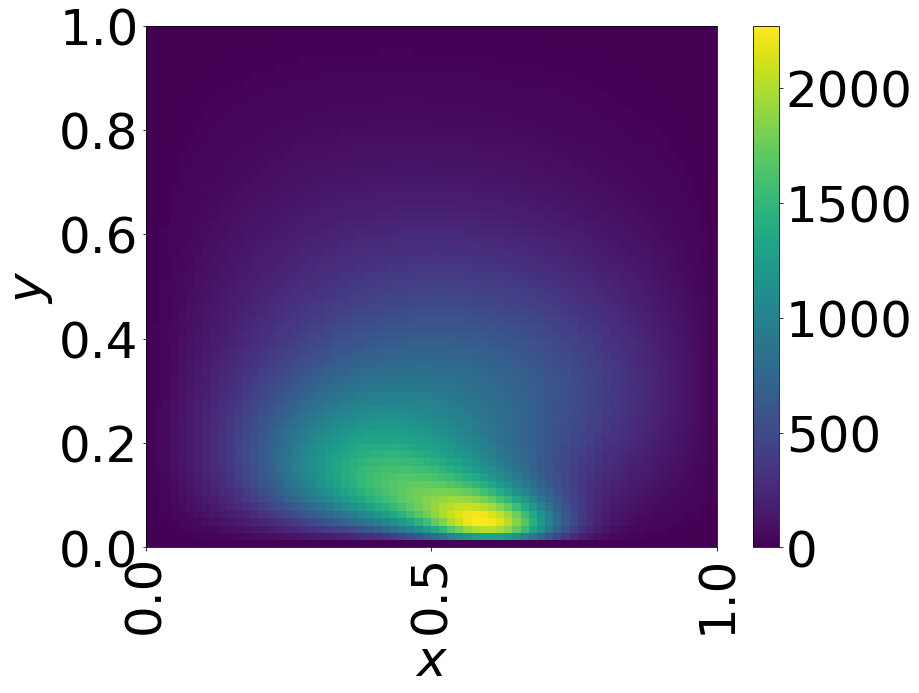}
        \caption{FOM}
    \end{subfigure}    
    \begin{subfigure}[b]{0.32\textwidth}
        \centering
        \includegraphics[width=\textwidth]{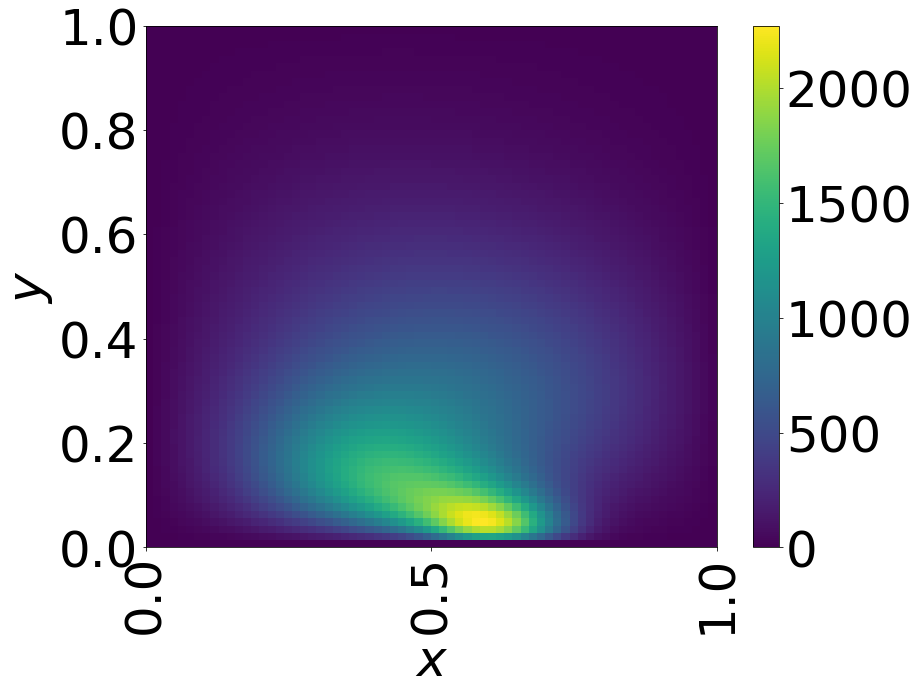}
        \caption{Galerkin ROM}
    \end{subfigure}
    \begin{subfigure}[b]{0.32\textwidth}
        \centering
        \includegraphics[width=\textwidth]{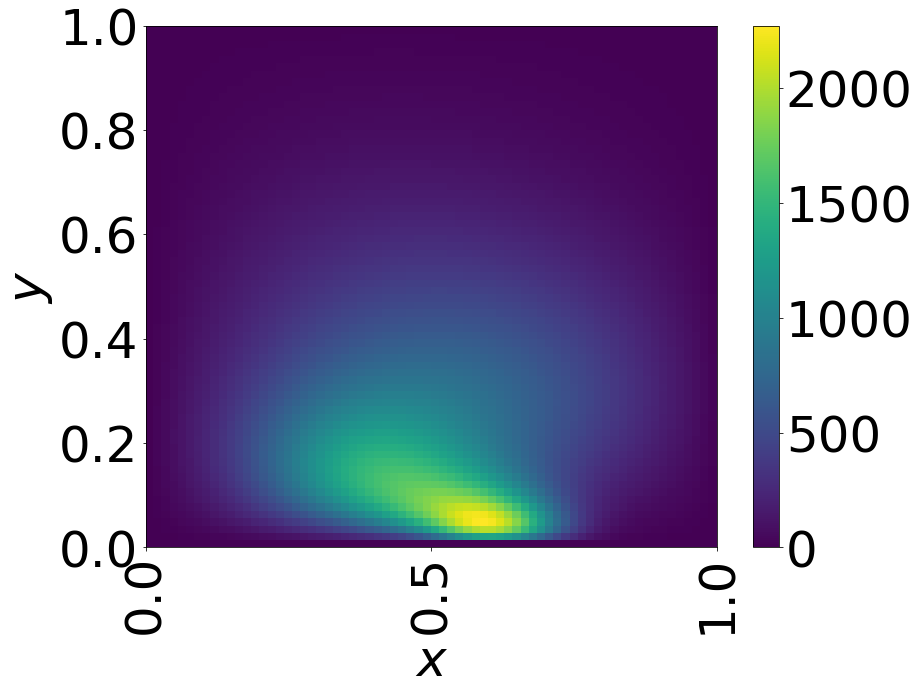}
        \caption{Petrov--Galerkin ROM}
    \end{subfigure}
    \caption{Solution snapshots of FOM, Galerkin ROM, and Petrov--Galerkin ROM at $t=2$.}
    \label{fig:ConvDiffSourceSol}
\end{figure}

The final time snapshots of FOM, Galerkin space--time ROM, and Petrov--Galerkin
space--time ROM are seen in Fig.~\ref{fig:ConvDiffSourceSol}. Both ROMs have a
basis size of $\nreducedspace=19$ and $\nreducedtime=3$, resulting in a
reduction factor of $(\nspace\ntime)\small/(\nreducedspace\nreducedtime)=4,176$.
For the Galerkin method, the FOM and space--time ROM simulation with
$\nreducedspace=19$ and $\nreducedtime=3$ takes an average time of
$6.1209\times 10^{-1}$ and $3.9780\times 10^{-3}$ seconds, respectively, resulting
in speed-up of $153.87$. For the Petrov--Galerkin method, the FOM and
space--time ROM simulation with $\nreducedspace=19$ and $\nreducedtime=3$ takes
an average time of $5.8780\times 10^{-1}$ and $4.2020\times 10^{-3}$ seconds,
respectively, resulting in speed-up of $139.89$. For accuracy, the Galerkin
method results in $2.174\times 10^{-1}$ \% relative error and $1.564\times 10^{3}$
space--time residual norm while the Petrov--Galerkin results in
$2.652\times 10^{-1}$ \% relative error and $1.550\times 10^{3}$ space--time
residual norm. 

\begin{figure}[!htbp]
  \centering
  \begin{subfigure}[b]{0.49\textwidth}
    \centering
    \includegraphics[width=\textwidth]{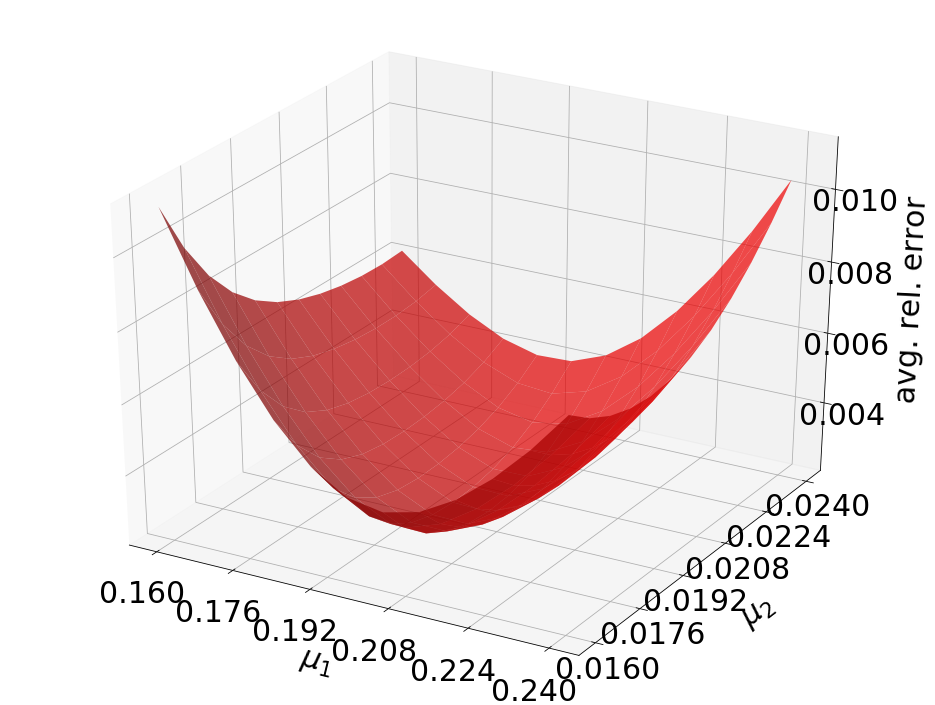}
    \caption{Glaerkin}
  \end{subfigure}    
  \begin{subfigure}[b]{0.49\textwidth}
    \centering
    \includegraphics[width=\textwidth]{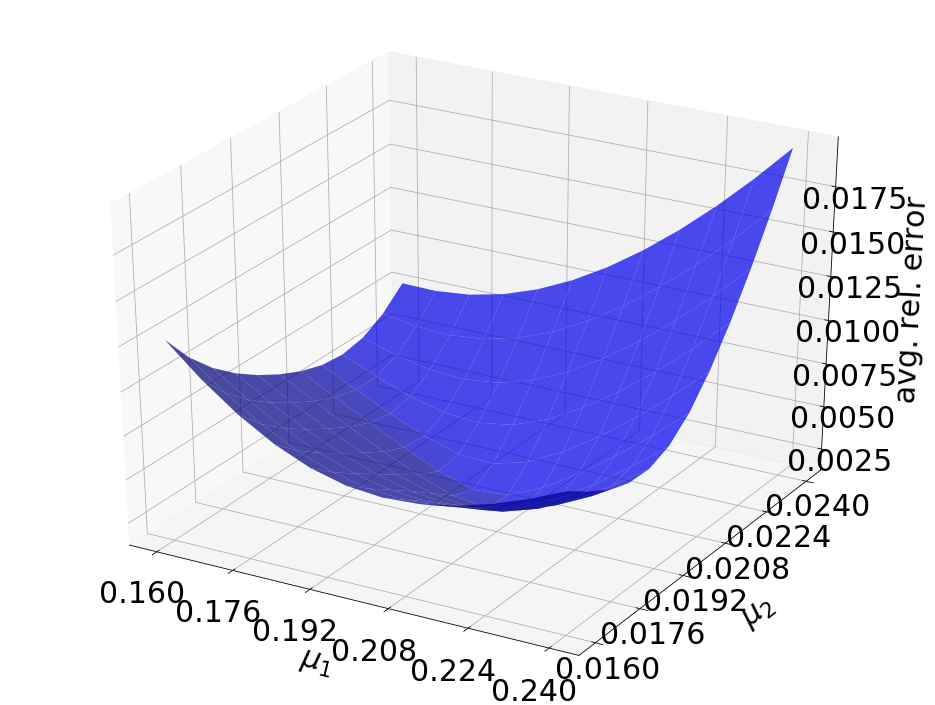}
    \caption{Petrov--Galerkin}
  \end{subfigure}
  \caption{The comparison of the Galerkin and Petrov--Galerkin ROMs for
  predictive cases}
  \label{fig:ConvDiffSourceErrVSParam}
\end{figure}

We investigate the numerical tests to see the generalization capability of both
Galerkin and Petrov--Galerkin ROMs. The train parameter set,
$(\mu_1,\mu_2)\in\{(0.195,0.018),(0.195,0.022),(0.205,0.018),(0.205,0.022)\}$ is
used to train a space--time ROMs with a basis of $\nreducedspace=19$ and
$\nreducedtime=3$. Then trained ROMs solve predictive cases with the test
parameter set,
$(\mu_1,\mu_2)\in\{\mu_1|\mu_1=0.160+0.08\small/11i,i=0,1,\cdots,11\}\times
\{\mu_2|\mu_2=0.016+0.008\small/11j,j=0,1,\cdots,11\}$.
Fig.~\ref{fig:ConvDiffSourceErrVSParam} shows the relative errors over the test
parameter set. The Galerkin and Petrov--Galerkin ROMs are the most accurate
within the range of the train parameter points, i.e., $[0.195,0.205]\times
[0.018,0.022]$. As the parameter points go beyond the train parameter domain,
the accuracy of the Galerkin and Petrov--Galerkin ROMs start to deteriorate
gradually. This implies that the Galerkin and Petrov--Galerkin ROMs have a trust
region. Its trust region should be determined by an application. For Galerkin
ROM, online speed-up is about $133$ in average and total time for ROM and FOM
are $68.49$ and $82.35$ seconds, respectively, resulting in total speed-up of
$1.20$. For Petrov--Galerkin ROM, online speed-up is about $138$ in average and
total time for ROM and FOM are $75.75$ and $86.04$ seconds, respectively,
resulting in total speed-up of $1.14$. Since the training time doesn't depend on
the number of test cases, we expect more speed-up for the larger number of test
cases.

\section{Conclusion}\label{sec:Conclusion}
In this work, we have formulated Galerkin and Petrov--Galerkin space--time ROMs
using block structures which enable us to implement the space--time ROM
operators efficiently. We also presented an \textit{a posteriori} error bound for
both Galerkin and Petrov--Galerkin space--time ROMs. We demonstrated that both
Galerkin and Petrov--Galerkin space--time ROMs solves 2D linear diffusion problems
and 2D linear convection diffusion problems accurately and efficiently. Both
space--time reduced order models were able to achieve $\mathcal{O}(10^{-3})$ to
$\mathcal{O}(10^{-4})$ relative errors with $\mathcal{O}(10^{2})$ speed-ups. We
also presented our Python codes used for the numerical examples in
Appendix~\ref{sec:codes} so that readers can easily reproduce our numerical
results. Furthermore, each Python code is less than 120 lines, demonstrating the ease of implementing our space--time ROMs.

We used a linear subspace based ROM which is suitable for accelerating physical
simulations whose solution space has a small Kolmogorov $n$-width. However, the
linear subspace based ROM is not able to represent advection-dominated or sharp
gradient solutions with a small number of bases. To address this challenge, a
nonlinear manifold based ROM can be used, and recently, a nonlinear manifold based
ROM has been developed for spatial ROMs \cite{lee2020model,kim2020fast}. In future
work, we aim to develop a nonlinear manifold based space--time ROM. 

\section*{Acknowledgments}
This work was performed at Lawrence Livermore National Laboratory and was
supported by LDRD program  (project 20-FS-007).  Youngkyu was also supported for
this work through generous funding from DTRA. Lawrence Livermore National
Laboratory is operated by Lawrence Livermore National Security, LLC, for the
U.S. Department of Energy, National Nuclear Security Administration under
Contract DE-AC52-07NA27344 and LLNL-JRNL-816093.

\appendix
\section{Python codes in less than 120 lines of code for all numerical models described in Section \ref{sec:NumericalResults}}\label{sec:codes}

The Python code used for the numerical examples described in this paper are
included in the following pages of the appendix are listed below. The total
number of lines in each of the files are denoted in the parentheses. Note that
we removed print statements of the results.
\begin{enumerate}
\setlength\itemsep{0.03125em}
\item All input code for the Galerkin Reduced Order Model for 2D Implicit Linear
  Diffusion Equation with Source Term (111 lines)
\item All input code for the Petrov--Galerkin Reduced Order Model for 2D
  Implicit Linear Diffusion Equation with Source Term (117 lines)
\item All input code for the Galerkin Reduced Order Model for 2D Implicit Linear
  Convection Diffusion Equation (114 lines)
\item All input code for the Petrov--Galerkin Reduced Order Model for 2D
  Implicit Linear Convection Diffusion Equation (116 lines)
\item All input code for the Galerkin Reduced Order Model for 2D Implicit Linear
  Convection Diffusion Equation with Source Term (115 lines)
\item All input code for the Petrov--Galerkin Reduced Order Model for 2D
  Implicit Linear Convection Diffusion Equation with Source Term (119 lines)
\end{enumerate}

\subsection{Galerkin Reduced Order Model for 2D Implicit Linear Diffusion
Equation with Source Term}
\fvset{frame=single,numbers=left,fontsize=\relscale{0.7}}
\begin{Verbatim}
import numpy as np; import scipy as sp
from scipy.sparse.linalg import spsolve

# train parameter space
mu1,mu2=np.meshgrid(np.linspace(-0.9,-0.5,2),np.linspace(-0.9,-0.5,2),indexing='ij')
mu1=mu1.flatten(); mu2=mu2.flatten(); no_para=np.size(mu1,axis=0)

# space and time domain
N=70; Ns=(N-1)**2; Nt=50; h=1/N; k=2/Nt
x1,x2=np.meshgrid(np.linspace(0,1,N+1)[1:-1],np.linspace(0,1,N+1)[1:-1],indexing='ij')
x1=x1.flatten(); x2=x2.flatten(); t=np.linspace(0,2,Nt+1)

# basic operators
e=np.ones(N-1); I=sp.sparse.eye(Ns,format='csc')
A1D=1/h**2*sp.sparse.spdiags(np.vstack((e,-2*e,e)),[-1,0,1],N-1,N-1,format='csc')
A2D=sp.sparse.kron(A1D,sp.sparse.eye(N-1,format='csc'),format='csc')\
    +sp.sparse.kron(sp.sparse.eye(N-1,format='csc'),A1D,format='csc')

# snapshot, free DoFs, No IC included
U=np.zeros((Ns,no_para*Nt)) 
for p in range(no_para):
    A=A2D-sp.sparse.diags(1/np.sqrt((x1-mu1[p])**2+(x2-mu2[p])**2),format='csc')
    f=(1/np.sqrt((x1-mu1[p])**2+(x2-mu2[p])**2).reshape(-1,1))@(np.sin(2*np.pi*t).reshape(1,-1))
    u=np.zeros((Ns,Nt+1)) # IC is zero
    for n in range(Nt):
        u[:,n+1]=spsolve(I-k*A,u[:,n]+k*f[:,n+1])      
    U[:,np.arange(p*Nt,(p+1)*Nt)]=u[:,1:]
    
W,S,VT=np.linalg.svd(U) # POD of solution snapshot

# test parameter space
nparam1=15; nparam2=15; ParamD1,ParamD2=np.meshgrid(np.linspace(-1.7,-0.2,nparam1),np.linspace(-1.7,-0.2,nparam2))
ParamD1=ParamD1.flatten(); ParamD2=ParamD2.flatten(); num_test=np.prod(ParamD1.shape)

# variables to store (also store ParamD1 and ParamD2)
UROM_g2=np.zeros((num_test,Nt+1,N+1,N+1)); UFOM_2=np.zeros((num_test,Nt+1,N+1,N+1))
AvgRelErr_g2=np.zeros((num_test)); STResROM_g2=np.zeros((num_test))

# Construct Phi_s
ns=5; PHIs=W[:,:ns]; PHIsT=PHIs.T

# construct Djk
nt=3; D=np.zeros((ns,nt*Nt))
for i in range(ns):
    Ri=VT[i,:]; Ri=Ri.reshape(-1,no_para,order='F')
    Wi,Si,ViT=np.linalg.svd(Ri); PHIti=Wi[:,:nt]
    for j in range(nt):
        D[i,Nt*j:Nt*(j+1)]=PHIti[:,j]

# construct PHIst for post processing, i.e., reconstruction
PHIst=np.zeros((Ns*Nt,ns*nt))
for i in range(Nt):
    for j in range(nt):
        PHIstij=np.zeros((Ns,ns))
        Dij=sp.sparse.diags(D[:,j*Nt+i],format='csc')
        PHIst[Ns*i:Ns*(i+1),ns*j:ns*(j+1)]=PHIs@Dij  
PHIstT=PHIst.T

for test_count in range(num_test): # predictive cases
    param1,param2=ParamD1[test_count],ParamD2[test_count] # set parameter

    # set model
    A=A2D-sp.sparse.diags(1/np.sqrt((x1-param1)**2+(x2-param2)**2),format='csc')
    f=(1/np.sqrt((x1-param1)**2+(x2-param2)**2).reshape(-1,1))@(np.sin(2*np.pi*t).reshape(1,-1))

    # construct Ast and fst using block structure (ust0 is zero)
    As=PHIsT@A@PHIs; fs=PHIsT@f
    Ast=np.zeros((ns*nt,ns*nt))
    for i in range(nt):
        for j in range(nt):
            Astij=np.zeros((ns,ns))
            for kk in range(Nt):
                Dik=sp.sparse.diags(D[:,i*Nt+kk],format='csc')
                Djk=sp.sparse.diags(D[:,j*Nt+kk],format='csc')
                Astij+=Dik@Djk-k*Dik@As@Djk
                if kk != Nt-1:
                    Dik_next=sp.sparse.diags(D[:,i*Nt+kk+1],format='csc')
                    Astij-=Dik_next@Djk   
            Ast[ns*i:ns*(i+1),ns*j:ns*(j+1)]=Astij  
    fst=np.zeros(ns*nt)
    for j in range(nt):
        for kk in range(Nt):
            Djk=sp.sparse.diags(D[:,j*Nt+kk],format='csc')
            fst[ns*j:ns*(j+1)]+=k*Djk@fs[:,kk+1]

    # Space-Time ROM (online phase)
    UstROM=PHIst@np.linalg.solve(Ast,fst)

    # FOM
    u=np.zeros((Ns,Nt+1))
    for n in range(Nt):
        u[:,n+1]=spsolve(I-k*A,u[:,n]+k*f[:,n+1])
    UstFOM=u[:,1:].flatten(order='F')
    
    AvgRelErr_g2[test_count]=np.linalg.norm(UstROM-UstFOM)/np.linalg.norm(UstFOM) # Avg. rel. error
    
    rstROM=np.zeros(Ns*Nt) # ST residual
    for n in range(1,Nt+1):
        if n==1:
            rstROM[(n-1)*Ns:n*Ns]=k*f[:,n]-(I-k*A).dot(UstROM[(n-1)*Ns:n*Ns])
        else:
            rstROM[(n-1)*Ns:n*Ns]=k*f[:,n]+UstROM[(n-2)*Ns:(n-1)*Ns]-(I-k*A).dot(UstROM[(n-1)*Ns:n*Ns])
    STResROM_g2[test_count]=np.linalg.norm(rstROM)
    
    # store solutions for each param
    urom=np.zeros((Nt+1,N+1,N+1)); ufom=np.zeros((Nt+1,N+1,N+1))
    urom[0]=np.zeros((N+1,N+1)); ufom[0]=np.zeros((N+1,N+1))
    for n in range(1,Nt+1):
        urom[n,1:-1,1:-1]=UstROM[(n-1)*Ns:n*Ns].reshape((N-1,N-1))
        ufom[n,1:-1,1:-1]=UstFOM[(n-1)*Ns:n*Ns].reshape((N-1,N-1))
    UROM_g2[test_count]=urom; UFOM_2[test_count]=ufom
\end{Verbatim}

\subsection{Petrov--Galerkin Reduced Order Model for 2D Implicit Linear
Diffusion Equation with Source Term}
\fvset{frame=single,numbers=left,fontsize=\relscale{0.7}}
\begin{Verbatim}
import numpy as np; import scipy as sp
from scipy.sparse.linalg import spsolve

# train parameter space
mu1,mu2=np.meshgrid(np.linspace(-0.9,-0.5,2),np.linspace(-0.9,-0.5,2),indexing='ij')
mu1=mu1.flatten(); mu2=mu2.flatten(); no_para=np.size(mu1,axis=0)

# space and time domain
N=70; Ns=(N-1)**2; Nt=50; h=1/N; k=2/Nt
x1,x2=np.meshgrid(np.linspace(0,1,N+1)[1:-1],np.linspace(0,1,N+1)[1:-1],indexing='ij')
x1=x1.flatten(); x2=x2.flatten(); t=np.linspace(0,2,Nt+1)

# basic operators
e=np.ones(N-1); I=sp.sparse.eye(Ns,format='csc')
A1D=1/h**2*sp.sparse.spdiags(np.vstack((e,-2*e,e)),[-1,0,1],N-1,N-1,format='csc')
A2D=sp.sparse.kron(A1D,sp.sparse.eye(N-1,format='csc'),format='csc')\
    +sp.sparse.kron(sp.sparse.eye(N-1,format='csc'),A1D,format='csc')

# snapshot, free DoFs, No IC included
U=np.zeros((Ns,no_para*Nt)) 
for p in range(no_para):
    A=A2D-sp.sparse.diags(1/np.sqrt((x1-mu1[p])**2+(x2-mu2[p])**2),format='csc')
    f=(1/np.sqrt((x1-mu1[p])**2+(x2-mu2[p])**2).reshape(-1,1))@(np.sin(2*np.pi*t).reshape(1,-1))
    u=np.zeros((Ns,Nt+1)) # IC is zero
    for n in range(Nt):
        u[:,n+1]=spsolve(I-k*A,u[:,n]+k*f[:,n+1])      
    U[:,np.arange(p*Nt,(p+1)*Nt)]=u[:,1:]
    
W,S,VT=np.linalg.svd(U) # POD of solution snapshot

# test parameter space
nparam1=15; nparam2=15; ParamD1,ParamD2=np.meshgrid(np.linspace(-1.7,-0.2,nparam1),np.linspace(-1.7,-0.2,nparam2))
ParamD1=ParamD1.flatten(); ParamD2=ParamD2.flatten(); num_test=np.prod(ParamD1.shape)

# variables to store (also store ParamD1 and ParamD2)
UROM_pg2=np.zeros((num_test,Nt+1,N+1,N+1)); UFOM_2=np.zeros((num_test,Nt+1,N+1,N+1))
AvgRelErr_pg2=np.zeros((num_test)); STResROM_pg2=np.zeros((num_test))

# Construct Phi_s
ns=5; PHIs=W[:,:ns]; PHIsT=PHIs.T

# construct Djk
nt=3;  D=np.zeros((ns,nt*Nt))
for i in range(ns):
    Ri=VT[i,:]; Ri=Ri.reshape(-1,no_para,order='F')
    Wi,Si,ViT=np.linalg.svd(Ri); PHIti=Wi[:,:nt]
    for j in range(nt):
        D[i,Nt*j:Nt*(j+1)]=PHIti[:,j]

# construct PHIst for post processing, i.e., reconstruction
PHIst=np.zeros((Ns*Nt,ns*nt))
for i in range(Nt):
    for j in range(nt):
        PHIstij=np.zeros((Ns,ns))
        Dij=sp.sparse.diags(D[:,j*Nt+i],format='csc')
        PHIst[Ns*i:Ns*(i+1),ns*j:ns*(j+1)]=PHIs@Dij  
PHIstT=PHIst.T

# predictive cases
for test_count in range(num_test):
    param1,param2=ParamD1[test_count],ParamD2[test_count] # set parameter

    # set model
    A=A2D-sp.sparse.diags(1/np.sqrt((x1-param1)**2+(x2-param2)**2),format='csc')
    f=(1/np.sqrt((x1-param1)**2+(x2-param2)**2).reshape(-1,1))@(np.sin(2*np.pi*t).reshape(1,-1))

    # construct Ast and fst using block structure (ust0 is zero)
    As=PHIsT@A@PHIs; fs=PHIsT@f
    ATAs=PHIsT@(I-k*A.T)@(I-k*A)@PHIs; AIs=PHIsT@(I-k*A)@PHIs; ATIs=PHIsT@(I-k*A.T)@PHIs; ATIfs=PHIsT@(I-k*A.T)@f
    Ast=np.zeros((ns*nt,ns*nt))
    for i in range(nt):
        for j in range(nt):
            Astij=np.zeros((ns,ns))
            for kk in range(Nt):
                Dik=sp.sparse.diags(D[:,i*Nt+kk],format='csc')
                Djk=sp.sparse.diags(D[:,j*Nt+kk],format='csc')
                Astij += Dik@ATAs@Djk
                if kk != Nt-1:
                    Dik_next=sp.sparse.diags(D[:,i*Nt+kk+1],format='csc')
                    Djk_next=sp.sparse.diags(D[:,j*Nt+kk+1],format='csc')
                    Astij += Dik@Djk - Dik@AIs@Djk_next - Dik_next@ATIs@Djk
            Ast[ns*i:ns*(i+1),ns*j:ns*(j+1)]=Astij  
    fst=np.zeros(ns*nt)
    for j in range(nt):
        for kk in range(Nt):
            Djk=sp.sparse.diags(D[:,j*Nt+kk],format='csc')
            fst[ns*j:ns*(j+1)]+= k*Djk@ATIfs[:,kk+1]
            if kk != Nt-1:
                fs_next = fs[:,kk+1+1]
                fst[ns*j:ns*(j+1)]-=k*Djk@fs_next 

    # Space-Time ROM (online phase)
    UstROM=PHIst@np.linalg.solve(Ast,fst)
    
    # FOM
    u=np.zeros((Ns,Nt+1))
    for n in range(Nt):
        u[:,n+1]=spsolve(I-k*A,u[:,n]+k*f[:,n+1])
    UstFOM=u[:,1:].flatten(order='F')

    AvgRelErr_pg2[test_count]=np.linalg.norm(UstROM-UstFOM)/np.linalg.norm(UstFOM) # Avg. rel. error

    rstROM=np.zeros(Ns*Nt) # ST residual
    for n in range(1,Nt+1):
        if n==1:
            rstROM[(n-1)*Ns:n*Ns]=k*f[:,n]-(I-k*A).dot(UstROM[(n-1)*Ns:n*Ns])
        else:
            rstROM[(n-1)*Ns:n*Ns]=k*f[:,n]+UstROM[(n-2)*Ns:(n-1)*Ns]-(I-k*A).dot(UstROM[(n-1)*Ns:n*Ns])
    STResROM_pg2[test_count]=np.linalg.norm(rstROM)

    # store solutions for each param
    urom=np.zeros((Nt+1,N+1,N+1)); ufom=np.zeros((Nt+1,N+1,N+1))
    urom[0]=np.zeros((N+1,N+1)); ufom[0]=np.zeros((N+1,N+1))
    for n in range(1,Nt+1):
        urom[n,1:-1,1:-1]=UstROM[(n-1)*Ns:n*Ns].reshape((N-1,N-1))
        ufom[n,1:-1,1:-1]=UstFOM[(n-1)*Ns:n*Ns].reshape((N-1,N-1))
    UROM_pg2[test_count]=urom; UFOM_2[test_count]=ufom
\end{Verbatim}

\subsection{Galerkin Reduced Order Model for 2D Implicit Linear Convection
Diffusion Equation}
\fvset{frame=single,numbers=left,fontsize=\relscale{0.7}}
\begin{Verbatim}
import numpy as np; import scipy as sp
from scipy.sparse.linalg import spsolve

# train parameter space
mu1,mu2=np.meshgrid(np.linspace(0.03,0.05,2),np.linspace(0.33,0.35,2),indexing='ij')
mu1=mu1.flatten(); mu2=mu2.flatten(); no_para=np.size(mu1,axis=0) 

# space and time domain
N=70; Ns=(N-1)**2; Nt=50; h=1./N; Tfinal = 1.; k=Tfinal/Nt
x1,x2=np.meshgrid(np.linspace(0,1,N+1)[1:-1],np.linspace(0,1,N+1)[1:-1],indexing='ij')

# IC
u0=100*np.sin(2*np.pi*x1)**3*np.sin(2*np.pi*x2)**3; u0[np.nonzero(x1>0.5)]=0.0; u0[np.nonzero(x2>0.5)]=0.0; u0=u0.flatten()

# basic operators
e=np.ones(N-1); I=sp.sparse.eye(Ns,format='csc')
A1D_diff=1/h**2*sp.sparse.spdiags(np.vstack((e,-2*e,e)),[-1,0,1],N-1,N-1,format='csc') 
A2D_diff=sp.sparse.kron(A1D_diff,sp.sparse.eye(N-1,format='csc'),format='csc')\
    +sp.sparse.kron(sp.sparse.eye(N-1,format='csc'),A1D_diff,format='csc')
A1D_conv=1/(h)*sp.sparse.spdiags(np.vstack((-1*e,1*e,0*e)),[-1,0,1],N-1,N-1,format='csc') 
A2D_conv=sp.sparse.kron(A1D_conv,sp.sparse.eye(N-1,format='csc'),format='csc')\
    +sp.sparse.kron(sp.sparse.eye(N-1,format='csc'),A1D_conv,format='csc')

# snapshot, free DoFs, No IC included
U=np.zeros((Ns,no_para*Nt))
for p in range(no_para):
    A=-mu1[p]*A2D_conv + mu2[p]*A2D_diff 
    u=np.zeros((Ns,Nt+1)); u[:,0]=u0
    for n in range(Nt):
        u[:,n+1]=spsolve(I-k*A,u[:,n])     
    U[:,np.arange(p*Nt,(p+1)*Nt)]=u[:,1:]

# POD of solution snapshot
W,S,VT=np.linalg.svd(U)

# test parameter space
nparam1=12; nparam2=12; ParamD1,ParamD2=np.meshgrid(np.linspace(0.01,0.07,nparam1),np.linspace(0.31,0.37,nparam2))
ParamD1=ParamD1.flatten(); ParamD2=ParamD2.flatten(); num_test=np.prod(ParamD1.shape)

# variables to store (also store ParamD1 and ParamD2)
UROM_g2=np.zeros((num_test,Nt+1,N+1,N+1)); UFOM_2=np.zeros((num_test,Nt+1,N+1,N+1))
AvgRelErr_g2=np.zeros((num_test)); STResROM_g2=np.zeros((num_test))

# Construct Phi_s
ns=5; PHIs=W[:,:ns]; PHIsT=PHIs.T

# construct Djk
nt=3; D=np.zeros((ns,nt*Nt))
for i in range(ns):
    Ri=VT[i,:]; Ri=Ri.reshape(-1,no_para,order='F')
    Wi,Si,ViT=np.linalg.svd(Ri); PHIti=Wi[:,:nt]
    for j in range(nt):
        D[i,Nt*j:Nt*(j+1)]=PHIti[:,j]

# construct PHIst for post processing, i.e., reconstruction
PHIst=np.zeros((Ns*Nt,ns*nt))
for i in range(Nt):
    for j in range(nt):
        PHIstij=np.zeros((Ns,ns))
        Dij=sp.sparse.diags(D[:,j*Nt+i],format='csc')
        PHIst[Ns*i:Ns*(i+1),ns*j:ns*(j+1)]=PHIs@Dij  
PHIstT=PHIst.T 

# predictive cases
for test_count in range(num_test):
    param1,param2=ParamD1[test_count],ParamD2[test_count] # set parameter mu

    A = -param1*A2D_conv + param2*A2D_diff # set model

    # construct Ast, ust (fst is zero)
    As=PHIsT@A@PHIs; us0=PHIsT@u0
    Ast=np.zeros((ns*nt,ns*nt))
    for i in range(nt):
        for j in range(nt):
            Astij=np.zeros((ns,ns))
            for kk in range(Nt):
                Dik=sp.sparse.diags(D[:,i*Nt+kk],format='csc')
                Djk=sp.sparse.diags(D[:,j*Nt+kk],format='csc')
                Astij+=Dik@Djk-k*Dik@As@Djk
                if kk != Nt-1:
                    Dik_next=sp.sparse.diags(D[:,i*Nt+kk+1],format='csc')
                    Astij-=Dik_next@Djk   
            Ast[ns*i:ns*(i+1),ns*j:ns*(j+1)]=Astij  
    ust=np.zeros(ns*nt)
    for j in range(nt):
        Djk=sp.sparse.diags(D[:,j*Nt],format='csc')
        ust[ns*j:ns*(j+1)] += Djk@us0

    # Space-Time ROM (online phase)
    UstROM=PHIst@np.linalg.solve(Ast,ust)

    # FOM
    u=np.zeros((Ns,Nt+1)); u[:,0] = u0
    for n in range(Nt):
        u[:,n+1]=spsolve(I-k*A,u[:,n])
    UstFOM=u[:,1:].flatten(order='F')

    AvgRelErr_g2[test_count]=np.linalg.norm(UstROM-UstFOM)/np.linalg.norm(UstFOM) # Avg. rel. error

    rstROM=np.zeros(Ns*Nt) # ST residual
    for n in range(1,Nt+1):
        if n==1:
            rstROM[(n-1)*Ns:n*Ns]=u0-(I-k*A).dot(UstROM[(n-1)*Ns:n*Ns])
        else:
            rstROM[(n-1)*Ns:n*Ns]=UstROM[(n-2)*Ns:(n-1)*Ns]-(I-k*A).dot(UstROM[(n-1)*Ns:n*Ns])
    STResROM_g2[test_count]=np.linalg.norm(rstROM)

    # store solutions for each param
    urom=np.zeros((Nt+1,N+1,N+1)); ufom=np.zeros((Nt+1,N+1,N+1))
    urom[0,1:-1,1:-1]=u0.reshape((N-1,N-1)); ufom[0,1:-1,1:-1]=u0.reshape((N-1,N-1))
    for n in range(1,Nt+1):
        urom[n,1:-1,1:-1]=UstROM[(n-1)*Ns:n*Ns].reshape((N-1,N-1))
        ufom[n,1:-1,1:-1]=UstFOM[(n-1)*Ns:n*Ns].reshape((N-1,N-1))
    UROM_g2[test_count]=urom; UFOM_2[test_count]=ufom
\end{Verbatim}

\subsection{Petrov--Galerkin Reduced Order Model for 2D Implicit Linear
Convection Diffusion Equation}
\fvset{frame=single,numbers=left,fontsize=\relscale{0.7}}
\begin{Verbatim}
import numpy as np; import scipy as sp
from scipy.sparse.linalg import spsolve

# train parameter space
mu1,mu2=np.meshgrid(np.linspace(0.03,0.05,2),np.linspace(0.33,0.35,2),indexing='ij')
mu1=mu1.flatten(); mu2=mu2.flatten(); no_para=np.size(mu1,axis=0) 

# space and time domain
N=70; Ns=(N-1)**2; Nt=50; h=1./N; Tfinal = 1.; k=Tfinal/Nt
x1,x2=np.meshgrid(np.linspace(0,1,N+1)[1:-1],np.linspace(0,1,N+1)[1:-1],indexing='ij')

# IC
u0=100*np.sin(2*np.pi*x1)**3*np.sin(2*np.pi*x2)**3; u0[np.nonzero(x1>0.5)]=0.0; u0[np.nonzero(x2>0.5)]=0.0; u0=u0.flatten()

# basic operators
e=np.ones(N-1); I=sp.sparse.eye(Ns,format='csc')
A1D_diff=1/h**2*sp.sparse.spdiags(np.vstack((e,-2*e,e)),[-1,0,1],N-1,N-1,format='csc') 
A2D_diff=sp.sparse.kron(A1D_diff,sp.sparse.eye(N-1,format='csc'),format='csc')\
    +sp.sparse.kron(sp.sparse.eye(N-1,format='csc'),A1D_diff,format='csc')
A1D_conv=1/(h)*sp.sparse.spdiags(np.vstack((-1*e,1*e,0*e)),[-1,0,1],N-1,N-1,format='csc') 
A2D_conv=sp.sparse.kron(A1D_conv,sp.sparse.eye(N-1,format='csc'),format='csc')\
    +sp.sparse.kron(sp.sparse.eye(N-1,format='csc'),A1D_conv,format='csc')

# snapshot, free DoFs, No IC included
U=np.zeros((Ns,no_para*Nt))
for p in range(no_para):
    A=-mu1[p]*A2D_conv + mu2[p]*A2D_diff 
    u=np.zeros((Ns,Nt+1)); u[:,0]=u0
    for n in range(Nt):
        u[:,n+1]=spsolve(I-k*A,u[:,n])     
    U[:,np.arange(p*Nt,(p+1)*Nt)]=u[:,1:]

# POD of solution snapshot
W,S,VT=np.linalg.svd(U)

# test parameter space
nparam1=12; nparam2=12; ParamD1,ParamD2=np.meshgrid(np.linspace(0.01,0.07,nparam1),np.linspace(0.31,0.37,nparam2))
ParamD1=ParamD1.flatten(); ParamD2=ParamD2.flatten(); num_test=np.prod(ParamD1.shape)

# variables to store (also store ParamD1 and ParamD2)
UROM_pg2=np.zeros((num_test,Nt+1,N+1,N+1)); UFOM_2=np.zeros((num_test,Nt+1,N+1,N+1))
AvgRelErr_pg2=np.zeros((num_test)); STResROM_pg2=np.zeros((num_test))

# Construct Phi_s
ns=5; PHIs=W[:,:ns]; PHIsT=PHIs.T

# construct Djk
nt=3; D=np.zeros((ns,nt*Nt))
for i in range(ns):
    Ri=VT[i,:]; Ri=Ri.reshape(-1,no_para,order='F')
    Wi,Si,ViT=np.linalg.svd(Ri); PHIti=Wi[:,:nt]
    for j in range(nt):
        D[i,Nt*j:Nt*(j+1)]=PHIti[:,j]

# construct PHIst for post processing, i.e., reconstruction
PHIst=np.zeros((Ns*Nt,ns*nt))
for i in range(Nt):
    for j in range(nt):
        PHIstij=np.zeros((Ns,ns))
        Dij=sp.sparse.diags(D[:,j*Nt+i],format='csc')
        PHIst[Ns*i:Ns*(i+1),ns*j:ns*(j+1)]=PHIs@Dij  
PHIstT=PHIst.T  

# predictive cases
for test_count in range(num_test):
    param1,param2=ParamD1[test_count],ParamD2[test_count] # set parameter mu

    A= -param1*A2D_conv + param2*A2D_diff # set model

    # construct Ast, ust (fst is zero)
    As=PHIsT@A@PHIs; us0=PHIsT@u0
    ATAs=PHIsT@(I-k*A.T)@(I-k*A)@PHIs; AIs=PHIsT@(I-k*A)@PHIs; ATIs=PHIsT@(I-k*A.T)@PHIs; ATIu0=PHIsT@(I-k*A.T)@u0
    Ast=np.zeros((ns*nt,ns*nt))
    for i in range(nt):
        for j in range(nt):
            Astij=np.zeros((ns,ns))
            for kk in range(Nt):
                Dik=sp.sparse.diags(D[:,i*Nt+kk],format='csc')
                Djk=sp.sparse.diags(D[:,j*Nt+kk],format='csc')
                Astij += Dik@ATAs@Djk
                if kk != Nt-1:
                    Dik_next=sp.sparse.diags(D[:,i*Nt+kk+1],format='csc')
                    Djk_next=sp.sparse.diags(D[:,j*Nt+kk+1],format='csc')
                    Astij += Dik@Djk - Dik@AIs@Djk_next - Dik_next@ATIs@Djk
            Ast[ns*i:ns*(i+1),ns*j:ns*(j+1)]=Astij  
    ust=np.zeros(ns*nt)
    for j in range(nt):
        Djk=sp.sparse.diags(D[:,j*Nt],format='csc')
        ust[ns*j:ns*(j+1)] += Djk@ATIu0

    # Space-Time ROM (online phase)
    UstROM=PHIst@np.linalg.solve(Ast,ust)

    # FOM
    u=np.zeros((Ns,Nt+1)); u[:,0] = u0
    for n in range(Nt):
        u[:,n+1]=spsolve(I-k*A,u[:,n])
    UstFOM=u[:,1:].flatten(order='F')
   
    AvgRelErr_pg2[test_count]=np.linalg.norm(UstROM-UstFOM)/np.linalg.norm(UstFOM)  # Avg. rel. error

    rstROM=np.zeros(Ns*Nt) # ST residual
    for n in range(1,Nt+1):
        if n==1:
            rstROM[(n-1)*Ns:n*Ns]=u0-(I-k*A).dot(UstROM[(n-1)*Ns:n*Ns])
        else:
            rstROM[(n-1)*Ns:n*Ns]=UstROM[(n-2)*Ns:(n-1)*Ns]-(I-k*A).dot(UstROM[(n-1)*Ns:n*Ns])
    STResROM_pg2[test_count]=np.linalg.norm(rstROM)

    # store solutions for each param
    urom=np.zeros((Nt+1,N+1,N+1)); ufom=np.zeros((Nt+1,N+1,N+1))
    urom[0,1:-1,1:-1]=u0.reshape((N-1,N-1)); ufom[0,1:-1,1:-1]=u0.reshape((N-1,N-1))
    for n in range(1,Nt+1):
        urom[n,1:-1,1:-1]=UstROM[(n-1)*Ns:n*Ns].reshape((N-1,N-1))
        ufom[n,1:-1,1:-1]=UstFOM[(n-1)*Ns:n*Ns].reshape((N-1,N-1))
    UROM_pg2[test_count]=urom; UFOM_2[test_count]=ufom
\end{Verbatim}

\subsection{Galerkin Reduced Order Model for 2D Implicit Linear Convection
Diffusion Equation with Source Term}
\fvset{frame=single,numbers=left,fontsize=\relscale{0.7}}
\begin{Verbatim}
import numpy as np; import scipy as sp
from scipy.sparse.linalg import spsolve

# train parameter space
mu1,mu2=np.meshgrid(np.linspace(0.195,0.205,2),np.linspace(0.018,0.022,2),indexing='xy')
mu1=mu1.flatten(); mu2=mu2.flatten(); no_para=np.size(mu1,axis=0) 

# space and time domain
N=70; Ns=(N-1)**2; Nt=50; h=1./N; Tfinal = 2.; k=Tfinal/Nt; t=np.linspace(0,Tfinal,Nt+1)
x1,x2=np.meshgrid(np.linspace(0,1,N+1)[1:-1],np.linspace(0,1,N+1)[1:-1],indexing='xy')
x1=x1.flatten(); x2=x2.flatten()

# basic operators
e=np.ones(N-1); I=sp.sparse.eye(Ns,format='csc')
A1D_diff=1/h**2*sp.sparse.spdiags(np.vstack((e,-2*e,e)),[-1,0,1],N-1,N-1,format='csc') 
A2D_diff=sp.sparse.kron(A1D_diff,sp.sparse.eye(N-1,format='csc'),format='csc')+\
    sp.sparse.kron(sp.sparse.eye(N-1,format='csc'),A1D_diff,format='csc')
A1D_conv=1/(h)*sp.sparse.spdiags(np.vstack((-1*e,1*e,0*e)),[-1,0,1],N-1,N-1,format='csc') 
A2D_conv=sp.sparse.kron(A1D_conv,sp.sparse.eye(N-1,format='csc'),format='csc') + \
         0.1*sp.sparse.kron(sp.sparse.eye(N-1,format='csc'),A1D_conv,format='csc')

# snapshot, free DoFs, No IC included
f=np.zeros((Ns,Nt+1))
for i in range(Nt+1):
    f[:,i] = 1e5*np.exp(-(((x1-0.5+(0.2*np.sin(2*np.pi*t[i])))/0.1)**2 + ((x2-0)/0.05)**2))
U=np.zeros((Ns,no_para*Nt))
for p in range(no_para):
    A=-mu1[p]*A2D_conv + mu2[p]*A2D_diff 
    u=np.zeros((Ns,Nt+1))
    for n in range(Nt):
        u[:,n+1]=spsolve(I-k*A,u[:,n]+k*f[:,n+1])     
    U[:,np.arange(p*Nt,(p+1)*Nt)]=u[:,1:]

# POD of solution snapshot
W,S,VT=np.linalg.svd(U)

# test parameter space
nparam1=12; nparam2=12; ParamD1,ParamD2=np.meshgrid(np.linspace(0.160,0.240,nparam1),np.linspace(0.016,0.024,nparam2))
ParamD1=ParamD1.flatten(); ParamD2=ParamD2.flatten(); num_test=np.prod(ParamD1.shape)

# variables to store (also store ParamD1 and ParamD2)
UROM_g2=np.zeros((num_test,Nt+1,N+1,N+1)); UFOM_2=np.zeros((num_test,Nt+1,N+1,N+1))
AvgRelErr_g2=np.zeros((num_test)); STResROM_g2=np.zeros((num_test))

# Construct Phi_s
ns=19; PHIs=W[:,:ns]; PHIsT=PHIs.T

# construct Djk
nt=3; D=np.zeros((ns,nt*Nt))
for i in range(ns):
    Ri=VT[i,:]; Ri=Ri.reshape(-1,no_para,order='F')
    Wi,Si,ViT=np.linalg.svd(Ri); PHIti=Wi[:,:nt]
    for j in range(nt):
        D[i,Nt*j:Nt*(j+1)]=PHIti[:,j]

# construct PHIst for post processing, i.e., reconstruction
PHIst=np.zeros((Ns*Nt,ns*nt))
for i in range(Nt):
    for j in range(nt):
        PHIstij=np.zeros((Ns,ns))
        Dij=sp.sparse.diags(D[:,j*Nt+i],format='csc')
        PHIst[Ns*i:Ns*(i+1),ns*j:ns*(j+1)]=PHIs@Dij  
PHIstT=PHIst.T

# predictive cases
for test_count in range(num_test):
    param1,param2=ParamD1[test_count],ParamD2[test_count] # set parameter mu

    A= -param1*A2D_conv + param2*A2D_diff # set model

    # construct Ast, fst (ust is zero)
    As=PHIsT@A@PHIs; fs=PHIsT@f
    Ast=np.zeros((ns*nt,ns*nt))
    for i in range(nt):
        for j in range(nt):
            Astij=np.zeros((ns,ns))
            for kk in range(Nt):
                Dik=sp.sparse.diags(D[:,i*Nt+kk],format='csc')
                Djk=sp.sparse.diags(D[:,j*Nt+kk],format='csc')
                Astij+=Dik@Djk-k*Dik@As@Djk
                if kk != Nt-1:
                    Dik_next=sp.sparse.diags(D[:,i*Nt+kk+1],format='csc')
                    Astij-=Dik_next@Djk   
            Ast[ns*i:ns*(i+1),ns*j:ns*(j+1)]=Astij  
    fst=np.zeros(ns*nt)
    for j in range(nt):
        for kk in range(Nt):
            Djk=sp.sparse.diags(D[:,j*Nt+kk],format='csc')
            fst[ns*j:ns*(j+1)]+=k*Djk@fs[:,kk+1]

    # Space-Time ROM (online phase)
    UstROM=PHIst@np.linalg.solve(Ast,fst)

    # FOM
    u=np.zeros((Ns,Nt+1))
    for n in range(Nt):
        u[:,n+1]=spsolve(I-k*A,u[:,n]+k*f[:,n+1])
    UstFOM=u[:,1:].flatten(order='F')

    AvgRelErr_g2[test_count]=np.linalg.norm(UstROM-UstFOM)/np.linalg.norm(UstFOM) # Avg. rel. error

    rstROM=np.zeros(Ns*Nt) # ST residual
    for n in range(1,Nt+1):
        if n==1:
            rstROM[(n-1)*Ns:n*Ns]=k*f[:,n]-(I-k*A).dot(UstROM[(n-1)*Ns:n*Ns])
        else:
            rstROM[(n-1)*Ns:n*Ns]=UstROM[(n-2)*Ns:(n-1)*Ns]+k*f[:,n]-(I-k*A).dot(UstROM[(n-1)*Ns:n*Ns])
    STResROM_g2[test_count]=np.linalg.norm(rstROM)

    # store solutions for each param
    urom=np.zeros((Nt+1,N+1,N+1)); ufom=np.zeros((Nt+1,N+1,N+1))
    for n in range(1,Nt+1):
        urom[n,1:-1,1:-1]=UstROM[(n-1)*Ns:n*Ns].reshape((N-1,N-1))
        ufom[n,1:-1,1:-1]=UstFOM[(n-1)*Ns:n*Ns].reshape((N-1,N-1))
    UROM_g2[test_count]=urom; UFOM_2[test_count]=ufom
\end{Verbatim}

\subsection{Petrov--Galerkin Reduced Order Model for 2D Implicit Linear
Convection Diffusion Equation with Source Term}
\fvset{frame=single,numbers=left,fontsize=\relscale{0.7}}
\begin{Verbatim}
import numpy as np; import scipy as sp
from scipy.sparse.linalg import spsolve

# train parameter space
mu1,mu2=np.meshgrid(np.linspace(0.195,0.205,2),np.linspace(0.018,0.022,2),indexing='xy')
mu1=mu1.flatten(); mu2=mu2.flatten(); no_para=np.size(mu1,axis=0)

# space and time domain
N=70; Ns=(N-1)**2; Nt=50; h=1./N; Tfinal = 1.; k=Tfinal/Nt; t=np.linspace(0,Tfinal,Nt+1)
x1,x2=np.meshgrid(np.linspace(0,1,N+1)[1:-1],np.linspace(0,1,N+1)[1:-1],indexing='xy')
x1=x1.flatten(); x2=x2.flatten()

# basic operators
e=np.ones(N-1); I=sp.sparse.eye(Ns,format='csc')
A1D_diff=1/h**2*sp.sparse.spdiags(np.vstack((e,-2*e,e)),[-1,0,1],N-1,N-1,format='csc')
A2D_diff=sp.sparse.kron(A1D_diff,sp.sparse.eye(N-1,format='csc'),format='csc')+\
    sp.sparse.kron(sp.sparse.eye(N-1,format='csc'),A1D_diff,format='csc')
A1D_conv=1/(h)*sp.sparse.spdiags(np.vstack((-1*e,1*e,0*e)),[-1,0,1],N-1,N-1,format='csc')
A2D_conv=sp.sparse.kron(A1D_conv,sp.sparse.eye(N-1,format='csc'),format='csc') + \
         0.1*sp.sparse.kron(sp.sparse.eye(N-1,format='csc'),A1D_conv,format='csc')

# snapshot, free DoFs, No IC included
f=np.zeros((Ns,Nt+1))
for i in range(Nt+1):
    f[:,i] = 1e5*np.exp(-(((x1-0.5+(0.2*np.sin(2*np.pi*t[i])))/0.1)**2 + ((x2-0)/0.05)**2))
U=np.zeros((Ns,no_para*Nt))
for p in range(no_para):
    A=-mu1[p]*A2D_conv + mu2[p]*A2D_diff
    u=np.zeros((Ns,Nt+1))
    for n in range(Nt):
        u[:,n+1]=spsolve(I-k*A,u[:,n]+k*f[:,n+1])
    U[:,np.arange(p*Nt,(p+1)*Nt)]=u[:,1:]

# POD of solution snapshot
W,S,VT=np.linalg.svd(U)

# test parameter space
nparam1=12; nparam2=12; ParamD1,ParamD2=np.meshgrid(np.linspace(0.160,0.240,nparam1),np.linspace(0.016,0.024,nparam2))
ParamD1=ParamD1.flatten(); ParamD2=ParamD2.flatten(); num_test=np.prod(ParamD1.shape)

# variables to store (also store ParamD1 and ParamD2)
UROM_pg2=np.zeros((num_test,Nt+1,N+1,N+1)); UFOM_pg2=np.zeros((num_test,Nt+1,N+1,N+1))
AvgRelErr_pg2=np.zeros((num_test)); STResROM_pg2=np.zeros((num_test))

# Construct Phi_s
ns=19; PHIs=W[:,:ns]; PHIsT=PHIs.T

# construct Djk
nt=3; D=np.zeros((ns,nt*Nt))
for i in range(ns):
    Ri=VT[i,:]; Ri=Ri.reshape(-1,no_para,order='F')
    Wi,Si,ViT=np.linalg.svd(Ri); PHIti=Wi[:,:nt]
    for j in range(nt):
        D[i,Nt*j:Nt*(j+1)]=PHIti[:,j]

# construct PHIst for post processing, i.e., reconstruction
PHIst=np.zeros((Ns*Nt,ns*nt))
for i in range(Nt):
    for j in range(nt):
        PHIstij=np.zeros((Ns,ns))
        Dij=sp.sparse.diags(D[:,j*Nt+i],format='csc')
        PHIst[Ns*i:Ns*(i+1),ns*j:ns*(j+1)]=PHIs@Dij
PHIstT=PHIst.T

# predictive cases
for test_count in range(num_test):
    param1,param2=ParamD1[test_count],ParamD2[test_count] # set parameter mu

    A= -param1*A2D_conv + param2*A2D_diff # set model

    # construct Ast, ust (fst is zero)
    As=PHIsT@A@PHIs; fs=PHIsT@f
    ATAs=PHIsT@(I-k*A.T)@(I-k*A)@PHIs; AIs=PHIsT@(I-k*A)@PHIs; ATIs=PHIsT@(I-k*A.T)@PHIs; ATIfs = PHIsT@(I-k*A.T)@f
    Ast=np.zeros((ns*nt,ns*nt))
    for i in range(nt):
        for j in range(nt):
            Astij=np.zeros((ns,ns))
            for kk in range(Nt):
                Dik=sp.sparse.diags(D[:,i*Nt+kk],format='csc')
                Djk=sp.sparse.diags(D[:,j*Nt+kk],format='csc')
                Astij += Dik@ATAs@Djk
                if kk != Nt-1:
                    Dik_next=sp.sparse.diags(D[:,i*Nt+kk+1],format='csc')
                    Djk_next=sp.sparse.diags(D[:,j*Nt+kk+1],format='csc')
                    Astij += Dik@Djk - Dik@AIs@Djk_next - Dik_next@ATIs@Djk
            Ast[ns*i:ns*(i+1),ns*j:ns*(j+1)]=Astij
    fst=np.zeros(ns*nt)
    for j in range(nt):
        for kk in range(Nt):
            Djk=sp.sparse.diags(D[:,j*Nt+kk],format='csc')
            fst[ns*j:ns*(j+1)]+= k*Djk@ATIfs[:,kk+1]
            if kk != Nt-1:
                fs_next = fs[:,kk+1+1]
                fst[ns*j:ns*(j+1)]-=k*Djk@fs_next
    # Space-Time ROM (online phase)
    UstROM=PHIst@np.linalg.solve(Ast,fst)

    # FOM
    u=np.zeros((Ns,Nt+1))
    for n in range(Nt):
        u[:,n+1]=spsolve(I-k*A,u[:,n]+k*f[:,n+1])
    UstFOM=u[:,1:].flatten(order='F')

    AvgRelErr_pg2[test_count]=np.linalg.norm(UstROM-UstFOM)/np.linalg.norm(UstFOM) # Avg. rel. error

    rstROM=np.zeros(Ns*Nt) # ST residual
    for n in range(1,Nt+1):
        if n==1:
            rstROM[(n-1)*Ns:n*Ns]=k*f[:,n]-(I-k*A).dot(UstROM[(n-1)*Ns:n*Ns])
        else:
            rstROM[(n-1)*Ns:n*Ns]=UstROM[(n-2)*Ns:(n-1)*Ns]+k*f[:,n]-(I-k*A).dot(UstROM[(n-1)*Ns:n*Ns])
    STResROM_pg2[test_count]=np.linalg.norm(rstROM)

    # store solutions for each param
    urom=np.zeros((Nt+1,N+1,N+1)); ufom=np.zeros((Nt+1,N+1,N+1))
    for n in range(1,Nt+1):
        urom[n,1:-1,1:-1]=UstROM[(n-1)*Ns:n*Ns].reshape((N-1,N-1))
        ufom[n,1:-1,1:-1]=UstFOM[(n-1)*Ns:n*Ns].reshape((N-1,N-1))
    UROM_pg2[test_count]=urom; UFOM_pg2[test_count]=ufom
\end{Verbatim}

\bibliographystyle{unsrt}
\bibliography{references}

\begin{thebibliography}{10}

\bibitem{Adam97}
L.~M. Adams.
\newblock Subcell balance methods for radiative transfer on arbitrary grids.
\newblock {\em Transport Th. Statis. Phys.}, 26(4,5):385--431, 1997.

\bibitem{ammar2006new}
Amine Ammar, B{\'e}chir Mokdad, Francisco Chinesta, and Roland Keunings.
\newblock A new family of solvers for some classes of multidimensional partial
  differential equations encountered in kinetic theory modeling of complex
  fluids.
\newblock {\em Journal of Non-Newtonian Fluid Mechanics}, 139(3):153--176,
  2006.

\bibitem{ammar2007new}
Amine Ammar, B{\'e}chir Mokdad, Francisco Chinesta, and Roland Keunings.
\newblock A new family of solvers for some classes of multidimensional partial
  differential equations encountered in kinetic theory modelling of complex
  fluids: Part ii: Transient simulation using space-time separated
  representations.
\newblock {\em Journal of Non-Newtonian Fluid Mechanics}, 144(2-3):98--121,
  2007.

\bibitem{bai2002krylov}
Zhaojun Bai.
\newblock Krylov subspace techniques for reduced-order modeling of large-scale
  dynamical systems.
\newblock {\em Applied numerical mathematics}, 43(1-2):9--44, 2002.

\bibitem{behne2019model}
Patrick Behne, Jean Ragusa, and Jim Morel.
\newblock Model-order reduction for sn radiation transport.
\newblock In {\em ANS International Conference on Mathematics and Computation
  (M\&C)}. Portland, OR, USA, 2019.

\bibitem{berkooz1993proper}
Gal Berkooz, Philip Holmes, and John~L Lumley.
\newblock The proper orthogonal decomposition in the analysis of turbulent
  flows.
\newblock {\em Annual review of fluid mechanics}, 25(1):539--575, 1993.

\bibitem{BiBr2009}
B.~L. Bihari and P.~N. Brown.
\newblock A linear algebraic analysis of diffusion synthetic acceleration for
  the boltzmann transport equation ii: The simple corner balance method.
\newblock {\em SIAM J. Numer. Anal.}, 47(3):1782--1826, 2009.

\bibitem{bihari2009linear}
BL~Bihari and Peter~N Brown.
\newblock A linear algebraic analysis of diffusion synthetic acceleration for
  the boltzmann transport equation ii: The simple corner balance method.
\newblock {\em SIAM Journal on Numerical Analysis}, 47(3):1782--1826, 2009.

\bibitem{brand2002incremental}
Matthew Brand.
\newblock Incremental singular value decomposition of uncertain data with
  missing values.
\newblock In {\em European Conference on Computer Vision}, pages 707--720.
  Springer, 2002.

\bibitem{buchan2013pod}
AG~Buchan, CC~Pain, F~Fang, and IM~Navon.
\newblock A pod reduced-order model for eigenvalue problems with application to
  reactor physics.
\newblock {\em International Journal for Numerical Methods in Engineering},
  95(12):1011--1032, 2013.

\bibitem{buchan2015pod}
Andrew~G Buchan, AA~Calloo, Mark~G Goffin, Steven Dargaville, Fangxin Fang,
  Christopher~C Pain, and Ionel~Michael Navon.
\newblock A {POD} reduced order model for resolving angular direction in
  neutron/photon transport problems.
\newblock {\em Journal of Computational Physics}, 296:138--157, 2015.

\bibitem{CaLa68}
B.~G. Carlson and K.~D. Lathrop.
\newblock Transport theory: The method of discrete ordinates.
\newblock In H.~Greenspan et~al., editors, {\em Computing Methods in Reactor
  Physics}, pages 166--266. Gordon and Breach, New York, 1968.

\bibitem{osti_1505575}
Youngsoo Choi, William~J. Arrighi, Dylan~M. Copeland, Robert~W. Anderson,
  Geoffrey~M. Oxberry, and USDOE National Nuclear~Security Administration.
\newblock librom, 10 2019.

\bibitem{choi2019space}
Youngsoo Choi and Kevin Carlberg.
\newblock Space--time least-squares petrov--galerkin projection for nonlinear
  model reduction.
\newblock {\em SIAM Journal on Scientific Computing}, 41(1):A26--A58, 2019.

\bibitem{coale2019areduced}
Joseph Coale and Dmitriy~Y. Anistratov.
\newblock A reduced-order model for thermal radiative transfer problems based
  on multilevel quasidiffusion method.
\newblock In {\em ANS International Conference on Mathematics and Computation
  (M\&C)}. Portland, OR, USA, 2019.

\bibitem{dominesey2019reduced}
Kurt Dominesey and Wei Ji.
\newblock Reduced-order modeling of neutron transport separated in space and
  angle via proper generalized decomposition.
\newblock In {\em ANS International Conference on Mathematics and Computation
  (M\&C)}. Portland, OR, USA, 2019.

\bibitem{dominesey2019areduced}
Kurt Dominesey and Wei Ji.
\newblock A reduced-order neutron diffusion model separated in space and energy
  via proper generalized decomposition.
\newblock In {\em Transactions of the American Nuclear Society}, volume 120,
  pages 457--460. Minneapolis, Minnesota, USA, 2019.

\bibitem{fareed2018error}
Hiba Fareed and John~R Singler.
\newblock Error analysis of an incremental pod algorithm for pde simulation
  data.
\newblock {\em arXiv preprint arXiv:1803.06313}, 2018.

\bibitem{gugercin2008h_2}
Serkan Gugercin, Athanasios~C Antoulas, and Christopher Beattie.
\newblock H\_2 model reduction for large-scale linear dynamical systems.
\newblock {\em SIAM journal on matrix analysis and applications},
  30(2):609--638, 2008.

\bibitem{hardy2019dynamic}
Zachary~K Hardy, Jim~E Morel, and Cory Ahrens.
\newblock Dynamic mode decomposition for subcritical metal systems.
\newblock {\em Nuclear Science and Engineering}, pages 1--13, 2019.

\bibitem{hinze2005proper}
Michael Hinze and Stefan Volkwein.
\newblock Proper orthogonal decomposition surrogate models for nonlinear
  dynamical systems: Error estimates and suboptimal control.
\newblock In {\em Dimension reduction of large-scale systems}, pages 261--306.
  Springer, 2005.

\bibitem{hotelling1933analysis}
Harold Hotelling.
\newblock Analysis of a complex of statistical variables into principal
  components.
\newblock {\em Journal of educational psychology}, 24(6):417, 1933.

\bibitem{kunisch2002galerkin}
Karl Kunisch and Stefan Volkwein.
\newblock Galerkin proper orthogonal decomposition methods for a general
  equation in fluid dynamics.
\newblock {\em SIAM Journal on Numerical analysis}, 40(2):492--515, 2002.

\bibitem{LM93}
E.\~E.\ Lewis and W.\~F.\ Miller.
\newblock {\em Computational Methods of Neutron Transport}.
\newblock American Nuclear Society, La Grange Park, IL, 1993.

\bibitem{Liboff80}
R.~L. Liboff.
\newblock {\em Introductory Quantum Mechanics}.
\newblock Holden-Day, Inc., San Francisco, 1980.

\bibitem{loeve1955}
Michel Loeve.
\newblock {\em Probability Theory}.
\newblock D. Van Nostrand, New York, 1955.

\bibitem{mcclarren2018acceleration}
Ryan~G McClarren and Terry~S Haut.
\newblock Acceleration of source iteration using the dynamic mode
  decomposition.
\newblock {\em arXiv preprint arXiv:1812.05241}, 2018.

\bibitem{moore1981principal}
Bruce Moore.
\newblock Principal component analysis in linear systems: Controllability,
  observability, and model reduction.
\newblock {\em IEEE transactions on automatic control}, 26(1):17--32, 1981.

\bibitem{mullis1976synthesis}
C~Mullis and RA~Roberts.
\newblock Synthesis of minimum roundoff noise fixed point digital filters.
\newblock {\em IEEE Transactions on Circuits and Systems}, 23(9):551--562,
  1976.

\bibitem{prince2019separated}
Zachary Prince and Jean Ragusa.
\newblock Separated representation of spatial dimensions in $s_n$ neutron
  transport using the proper generalized decomposition.
\newblock In {\em ANS International Conference on Mathematics and Computation
  (M\&C)}. Portland, OR, USA, 2019.

\bibitem{prince2019parametric}
Zachary~M Prince and Jean~C Ragusa.
\newblock Parametric uncertainty quantification using proper generalized
  decomposition applied to neutron diffusion.
\newblock {\em International Journal for Numerical Methods in Engineering},
  pages 1--23, 1993.

\bibitem{reed2015energy}
Richard~L Reed and Jeremy~A Roberts.
\newblock An energy basis for response matrix methods based on the
  karhunen--lo{\'e}ve transform.
\newblock {\em Annals of Nuclear Energy}, 78:70--80, 2015.

\bibitem{sartori2014comparison}
A.~Sartori, D.~Baroli, A.~Cammi, D.~Chiesa, L.~Luzzi, R.~Ponciroli,
  E.~Previtali, M.E. Ricotti, G.~Rozza, and M.~Sisti.
\newblock Comparison of a modal method and a proper orthogonal decomposition
  approach for multi-group time-dependent reactor spatial kinetics.
\newblock {\em Annals of Nuclear Energy}, 71:217--229, 2014.

\bibitem{schmid2010dynamic}
Peter~J Schmid.
\newblock Dynamic mode decomposition of numerical and experimental data.
\newblock {\em Journal of fluid mechanics}, 656:5--28, 2010.

\bibitem{sirovich1987turbulence}
Lawrence Sirovich.
\newblock Turbulence and the dynamics of coherent structures. i. coherent
  structures.
\newblock {\em Quarterly of applied mathematics}, 45(3):561--571, 1987.

\bibitem{star2019pod}
S~Kelbij Star, Francesco Belloni, Gert Van~den Eynde, and Joris Degroote.
\newblock Pod-identification reduced order model of linear transport equations
  for control purposes.
\newblock {\em International Journal for Numerical Methods in Fluids},
  90(8):375--388, 2019.

\bibitem{tu2013dynamic}
Jonathan~H Tu, Clarence~W Rowley, Dirk~M Luchtenburg, Steven~L Brunton, and
  J~Nathan Kutz.
\newblock On dynamic mode decomposition: theory and applications.
\newblock {\em arXiv preprint arXiv:1312.0041}, 2013.

\bibitem{urban2014improved}
Karsten Urban and Anthony Patera.
\newblock An improved error bound for reduced basis approximation of linear
  parabolic problems.
\newblock {\em Mathematics of Computation}, 83(288):1599--1615, 2014.

\bibitem{wols2010transient}
Frank Wols.
\newblock {\em Transient analyses of accelerator driven systems using modal
  expansion techniques}.
\newblock PhD thesis, Delft University of Technology, 2010.

\bibitem{yano2014space2}
Masayuki Yano.
\newblock A space-time petrov--galerkin certified reduced basis method:
  Application to the boussinesq equations.
\newblock {\em SIAM Journal on Scientific Computing}, 36(1):A232--A266, 2014.

\bibitem{yano2014space1}
Masayuki Yano, Anthony~T Patera, and Karsten Urban.
\newblock A space-time hp-interpolation-based certified reduced basis method
  for burgers' equation.
\newblock {\em Mathematical Models and Methods in Applied Sciences},
  24(09):1903--1935, 2014.

\end{thebibliography}


\begin{thebibliography}{10}

\bibitem{choi2020space}
Youngsoo Choi, Peter Brown, William Arrighi, Robert Anderson, and Kevin Huynh.
\newblock Space--time reduced order model for large-scale linear dynamical
  systems with application to boltzmann transport problems.
\newblock {\em Journal of Computational Physics}, page 109845, 2020.

\bibitem{mullis1976synthesis}
C~Mullis and RA~Roberts.
\newblock Synthesis of minimum roundoff noise fixed point digital filters.
\newblock {\em IEEE Transactions on Circuits and Systems}, 23(9):551--562,
  1976.

\bibitem{moore1981principal}
Bruce Moore.
\newblock Principal component analysis in linear systems: Controllability,
  observability, and model reduction.
\newblock {\em IEEE transactions on automatic control}, 26(1):17--32, 1981.

\bibitem{willcox2002balanced}
Karen Willcox and Jaime Peraire.
\newblock Balanced model reduction via the proper orthogonal decomposition.
\newblock {\em AIAA journal}, 40(11):2323--2330, 2002.

\bibitem{willcox2005fourier}
Karen Willcox and Alexandre Megretski.
\newblock Fourier series for accurate, stable, reduced-order models in
  large-scale linear applications.
\newblock {\em SIAM Journal on Scientific Computing}, 26(3):944--962, 2005.

\bibitem{heinkenschloss2008balanced}
Matthias Heinkenschloss, Danny~C Sorensen, and Kai Sun.
\newblock Balanced truncation model reduction for a class of descriptor systems
  with application to the oseen equations.
\newblock {\em SIAM Journal on Scientific Computing}, 30(2):1038--1063, 2008.

\bibitem{sandberg2004balanced}
Henrik Sandberg and Anders Rantzer.
\newblock Balanced truncation of linear time-varying systems.
\newblock {\em IEEE Transactions on automatic control}, 49(2):217--229, 2004.

\bibitem{hartmann2010balanced}
Carsten Hartmann, Valentina-Mira Vulcanov, and Christof Sch{\"u}tte.
\newblock Balanced truncation of linear second-order systems: a hamiltonian
  approach.
\newblock {\em Multiscale Modeling \& Simulation}, 8(4):1348--1367, 2010.

\bibitem{petreczky2013balanced}
Mih{\'a}ly Petreczky, Rafael Wisniewski, and John Leth.
\newblock Balanced truncation for linear switched systems.
\newblock {\em Nonlinear Analysis: Hybrid Systems}, 10:4--20, 2013.

\bibitem{ma2010snapshot}
Zhanhua Ma, Clarence~W Rowley, and Gilead Tadmor.
\newblock Snapshot-based balanced truncation for linear time-periodic systems.
\newblock {\em IEEE Transactions on Automatic Control}, 55(2):469--473, 2010.

\bibitem{bai2002krylov}
Zhaojun Bai.
\newblock Krylov subspace techniques for reduced-order modeling of large-scale
  dynamical systems.
\newblock {\em Applied numerical mathematics}, 43(1-2):9--44, 2002.

\bibitem{gugercin2008h_2}
Serkan Gugercin, Athanasios~C Antoulas, and Christopher Beattie.
\newblock H\_2 model reduction for large-scale linear dynamical systems.
\newblock {\em SIAM journal on matrix analysis and applications},
  30(2):609--638, 2008.

\bibitem{astolfi2010model}
Alessandro Astolfi.
\newblock Model reduction by moment matching for linear and nonlinear systems.
\newblock {\em IEEE Transactions on Automatic Control}, 55(10):2321--2336,
  2010.

\bibitem{chiprout1992generalized}
Eli Chiprout and Michael Nakhla.
\newblock Generalized moment-matching methods for transient analysis of
  interconnect networks.
\newblock In {\em [1992] Proceedings 29th ACM/IEEE Design Automation
  Conference}, pages 201--206. IEEE, 1992.

\bibitem{pratesi2006generalized}
Marco Pratesi, Fortunato Santucci, and Fabio Graziosi.
\newblock Generalized moment matching for the linear combination of lognormal
  rvs: application to outage analysis in wireless systems.
\newblock {\em IEEE Transactions on Wireless Communications}, 5(5):1122--1132,
  2006.

\bibitem{ammar2006new}
Amine Ammar, B{\'e}chir Mokdad, Francisco Chinesta, and Roland Keunings.
\newblock A new family of solvers for some classes of multidimensional partial
  differential equations encountered in kinetic theory modeling of complex
  fluids.
\newblock {\em Journal of Non-Newtonian Fluid Mechanics}, 139(3):153--176,
  2006.

\bibitem{ammar2007new}
Amine Ammar, B{\'e}chir Mokdad, Francisco Chinesta, and Roland Keunings.
\newblock A new family of solvers for some classes of multidimensional partial
  differential equations encountered in kinetic theory modelling of complex
  fluids: Part ii: Transient simulation using space-time separated
  representations.
\newblock {\em Journal of Non-Newtonian Fluid Mechanics}, 144(2-3):98--121,
  2007.

\bibitem{chinesta2010proper}
Francisco Chinesta, Amine Ammar, and El{\'\i}as Cueto.
\newblock Proper generalized decomposition of multiscale models.
\newblock {\em International Journal for Numerical Methods in Engineering},
  83(8-9):1114--1132, 2010.

\bibitem{pruliere2010deterministic}
Etienne Pruliere, Francisco Chinesta, and Amine Ammar.
\newblock On the deterministic solution of multidimensional parametric models
  using the proper generalized decomposition.
\newblock {\em Mathematics and Computers in Simulation}, 81(4):791--810, 2010.

\bibitem{chinesta2011overview}
Francisco Chinesta, Amine Ammar, Adrien Leygue, and Roland Keunings.
\newblock An overview of the proper generalized decomposition with applications
  in computational rheology.
\newblock {\em Journal of Non-Newtonian Fluid Mechanics}, 166(11):578--592,
  2011.

\bibitem{giner2013proper}
Eugenio Giner, Brice Bognet, Juan~J R{\'o}denas, Adrien Leygue, F~Javier
  Fuenmayor, and Francisco Chinesta.
\newblock The proper generalized decomposition (pgd) as a numerical procedure
  to solve 3d cracked plates in linear elastic fracture mechanics.
\newblock {\em International Journal of Solids and Structures},
  50(10):1710--1720, 2013.

\bibitem{barbarulo2014proper}
Andrea Barbarulo, Pierre Ladev{\`e}ze, Herv{\'e} Riou, and Louis Kovalevsky.
\newblock Proper generalized decomposition applied to linear acoustic: a new
  tool for broad band calculation.
\newblock {\em Journal of Sound and Vibration}, 333(11):2422--2431, 2014.

\bibitem{amsallem2012stabilization}
David Amsallem and Charbel Farhat.
\newblock Stabilization of projection-based reduced-order models.
\newblock {\em International Journal for Numerical Methods in Engineering},
  91(4):358--377, 2012.

\bibitem{amsallem2008interpolation}
David Amsallem and Charbel Farhat.
\newblock Interpolation method for adapting reduced-order models and
  application to aeroelasticity.
\newblock {\em AIAA journal}, 46(7):1803--1813, 2008.

\bibitem{thomas2003three}
Jeffrey~P Thomas, Earl~H Dowell, and Kenneth~C Hall.
\newblock Three-dimensional transonic aeroelasticity using proper orthogonal
  decomposition-based reduced-order models.
\newblock {\em Journal of Aircraft}, 40(3):544--551, 2003.

\bibitem{hall2000proper}
Kenneth~C Hall, Jeffrey~P Thomas, and Earl~H Dowell.
\newblock Proper orthogonal decomposition technique for transonic unsteady
  aerodynamic flows.
\newblock {\em AIAA journal}, 38(10):1853--1862, 2000.

\bibitem{chinesta2011short}
Francisco Chinesta, Pierre Ladeveze, and Elias Cueto.
\newblock A short review on model order reduction based on proper generalized
  decomposition.
\newblock {\em Archives of Computational Methods in Engineering}, 18(4):395,
  2011.

\bibitem{mayo2007framework}
AJ~Mayo and AC~Antoulas.
\newblock A framework for the solution of the generalized realization problem.
\newblock {\em Linear algebra and its applications}, 425(2-3):634--662, 2007.

\bibitem{scarciotti2017data}
Giordano Scarciotti and Alessandro Astolfi.
\newblock Data-driven model reduction by moment matching for linear and
  nonlinear systems.
\newblock {\em Automatica}, 79:340--351, 2017.

\bibitem{schmid2010dynamic}
Peter~J Schmid.
\newblock Dynamic mode decomposition of numerical and experimental data.
\newblock {\em Journal of fluid mechanics}, 656:5--28, 2010.

\bibitem{chen2012variants}
Kevin~K Chen, Jonathan~H Tu, and Clarence~W Rowley.
\newblock Variants of dynamic mode decomposition: boundary condition, koopman,
  and fourier analyses.
\newblock {\em Journal of nonlinear science}, 22(6):887--915, 2012.

\bibitem{williams2015data}
Matthew~O Williams, Ioannis~G Kevrekidis, and Clarence~W Rowley.
\newblock A data--driven approximation of the koopman operator: Extending
  dynamic mode decomposition.
\newblock {\em Journal of Nonlinear Science}, 25(6):1307--1346, 2015.

\bibitem{takeishi2017learning}
Naoya Takeishi, Yoshinobu Kawahara, and Takehisa Yairi.
\newblock Learning koopman invariant subspaces for dynamic mode decomposition.
\newblock In {\em Advances in Neural Information Processing Systems}, pages
  1130--1140, 2017.

\bibitem{askham2018variable}
Travis Askham and J~Nathan Kutz.
\newblock Variable projection methods for an optimized dynamic mode
  decomposition.
\newblock {\em SIAM Journal on Applied Dynamical Systems}, 17(1):380--416,
  2018.

\bibitem{schmid2011applications}
Peter~J Schmid, Larry Li, Matthew~P Juniper, and O~Pust.
\newblock Applications of the dynamic mode decomposition.
\newblock {\em Theoretical and Computational Fluid Dynamics}, 25(1-4):249--259,
  2011.

\bibitem{kutz2016multiresolution}
J~Nathan Kutz, Xing Fu, and Steven~L Brunton.
\newblock Multiresolution dynamic mode decomposition.
\newblock {\em SIAM Journal on Applied Dynamical Systems}, 15(2):713--735,
  2016.

\bibitem{li2017extended}
Qianxiao Li, Felix Dietrich, Erik~M Bollt, and Ioannis~G Kevrekidis.
\newblock Extended dynamic mode decomposition with dictionary learning: A
  data-driven adaptive spectral decomposition of the koopman operator.
\newblock {\em Chaos: An Interdisciplinary Journal of Nonlinear Science},
  27(10):103111, 2017.

\bibitem{proctor2016dynamic}
Joshua~L Proctor, Steven~L Brunton, and J~Nathan Kutz.
\newblock Dynamic mode decomposition with control.
\newblock {\em SIAM Journal on Applied Dynamical Systems}, 15(1):142--161,
  2016.

\bibitem{tu2013dynamic}
Jonathan~H Tu, Clarence~W Rowley, Dirk~M Luchtenburg, Steven~L Brunton, and
  J~Nathan Kutz.
\newblock On dynamic mode decomposition: theory and applications.
\newblock {\em arXiv preprint arXiv:1312.0041}, 2013.

\bibitem{kutz2016dynamic}
J~Nathan Kutz, Steven~L Brunton, Bingni~W Brunton, and Joshua~L Proctor.
\newblock {\em Dynamic mode decomposition: data-driven modeling of complex
  systems}.
\newblock SIAM, 2016.

\bibitem{berkooz1993proper}
Gal Berkooz, Philip Holmes, and John~L Lumley.
\newblock The proper orthogonal decomposition in the analysis of turbulent
  flows.
\newblock {\em Annual review of fluid mechanics}, 25(1):539--575, 1993.

\bibitem{gubisch2017proper}
Martin Gubisch and Stefan Volkwein.
\newblock Proper orthogonal decomposition for linear-quadratic optimal control.
\newblock {\em Model reduction and approximation: theory and algorithms}, 5:66,
  2017.

\bibitem{kunisch2001galerkin}
Karl Kunisch and Stefan Volkwein.
\newblock Galerkin proper orthogonal decomposition methods for parabolic
  problems.
\newblock {\em Numerische mathematik}, 90(1):117--148, 2001.

\bibitem{hinze2008error}
Michael Hinze and Stefan Volkwein.
\newblock Error estimates for abstract linear--quadratic optimal control
  problems using proper orthogonal decomposition.
\newblock {\em Computational Optimization and Applications}, 39(3):319--345,
  2008.

\bibitem{kerschen2005method}
Gaetan Kerschen, Jean-claude Golinval, Alexander~F Vakakis, and Lawrence~A
  Bergman.
\newblock The method of proper orthogonal decomposition for dynamical
  characterization and order reduction of mechanical systems: an overview.
\newblock {\em Nonlinear dynamics}, 41(1-3):147--169, 2005.

\bibitem{bamer2012application}
Franz Bamer and Christian Bucher.
\newblock Application of the proper orthogonal decomposition for linear and
  nonlinear structures under transient excitations.
\newblock {\em Acta Mechanica}, 223(12):2549--2563, 2012.

\bibitem{atwell2001proper}
Jeanne~A Atwell and Belinda~B King.
\newblock Proper orthogonal decomposition for reduced basis feedback
  controllers for parabolic equations.
\newblock {\em Mathematical and computer modelling}, 33(1-3):1--19, 2001.

\bibitem{rathinam2003new}
Muruhan Rathinam and Linda~R Petzold.
\newblock A new look at proper orthogonal decomposition.
\newblock {\em SIAM Journal on Numerical Analysis}, 41(5):1893--1925, 2003.

\bibitem{kahlbacher2007galerkin}
Martin Kahlbacher and Stefan Volkwein.
\newblock Galerkin proper orthogonal decomposition methods for parameter
  dependent elliptic systems.
\newblock {\em Discussiones Mathematicae, Differential Inclusions, Control and
  Optimization}, 27(1):95--117, 2007.

\bibitem{bonnet1994stochastic}
Jean~P Bonnet, Daniel~R Cole, Jo{\"e}l Delville, Mark~N Glauser, and Lawrence~S
  Ukeiley.
\newblock Stochastic estimation and proper orthogonal decomposition:
  complementary techniques for identifying structure.
\newblock {\em Experiments in fluids}, 17(5):307--314, 1994.

\bibitem{placzek2008hybrid}
Antoine Placzek, D-M Tran, and Roger Ohayon.
\newblock Hybrid proper orthogonal decomposition formulation for linear
  structural dynamics.
\newblock {\em Journal of Sound and Vibration}, 318(4-5):943--964, 2008.

\bibitem{legresley2000airfoil}
Patrick LeGresley and Juan Alonso.
\newblock Airfoil design optimization using reduced order models based on
  proper orthogonal decomposition.
\newblock In {\em Fluids 2000 conference and exhibit}, page 2545, 2000.

\bibitem{efe2003proper}
Mehmet~Onder Efe and Hitay Ozbay.
\newblock Proper orthogonal decomposition for reduced order modeling: 2d heat
  flow.
\newblock In {\em Proceedings of 2003 IEEE Conference on Control Applications,
  2003. CCA 2003.}, volume~2, pages 1273--1277. IEEE, 2003.

\bibitem{choi2019space}
Youngsoo Choi and Kevin Carlberg.
\newblock Space--time least-squares petrov--galerkin projection for nonlinear
  model reduction.
\newblock {\em SIAM Journal on Scientific Computing}, 41(1):A26--A58, 2019.

\bibitem{urban2014improved}
Karsten Urban and Anthony Patera.
\newblock An improved error bound for reduced basis approximation of linear
  parabolic problems.
\newblock {\em Mathematics of Computation}, 83(288):1599--1615, 2014.

\bibitem{yano2014space1}
Masayuki Yano, Anthony~T Patera, and Karsten Urban.
\newblock A space-time hp-interpolation-based certified reduced basis method
  for burgers' equation.
\newblock {\em Mathematical Models and Methods in Applied Sciences},
  24(09):1903--1935, 2014.

\bibitem{yano2014space2}
Masayuki Yano.
\newblock A space-time petrov--galerkin certified reduced basis method:
  Application to the boussinesq equations.
\newblock {\em SIAM Journal on Scientific Computing}, 36(1):A232--A266, 2014.

\bibitem{yoon2010structural}
Gil~Ho Yoon.
\newblock Structural topology optimization for frequency response problem using
  model reduction schemes.
\newblock {\em Computer Methods in Applied Mechanics and Engineering},
  199(25-28):1744--1763, 2010.

\bibitem{amir2010efficient}
Oded Amir, Mathias Stolpe, and Ole Sigmund.
\newblock Efficient use of iterative solvers in nested topology optimization.
\newblock {\em Structural and Multidisciplinary Optimization}, 42(1):55--72,
  2010.

\bibitem{amsallem2015design}
David Amsallem, Matthew Zahr, Youngsoo Choi, and Charbel Farhat.
\newblock Design optimization using hyper-reduced-order modelsvd.
\newblock {\em Structural and Multidisciplinary Optimization}, 51(4):919--940,
  2015.

\bibitem{gogu2015improving}
Christian Gogu.
\newblock Improving the efficiency of large scale topology optimization through
  on-the-fly reduced order model construction.
\newblock {\em International Journal for Numerical Methods in Engineering},
  101(4):281--304, 2015.

\bibitem{choi2020gradient}
Youngsoo Choi, Gabriele Boncoraglio, Spenser Anderson, David Amsallem, and
  Charbel Farhat.
\newblock Gradient-based constrained optimization using a database of linear
  reduced-order models.
\newblock {\em Journal of Computational Physics}, page 109787, 2020.

\bibitem{choi2019accelerating}
Youngsoo Choi, Geoffrey Oxberry, Daniel White, and Trenton Kirchdoerfer.
\newblock Accelerating design optimization using reduced order models.
\newblock {\em arXiv preprint arXiv:1909.11320}, 2019.

\bibitem{white2020dual}
Daniel~A White, Youngsoo Choi, and Jun Kudo.
\newblock A dual mesh method with adaptivity for stress-constrained topology
  optimization.
\newblock {\em Structural and Multidisciplinary Optimization}, 61(2):749--762,
  2020.

\bibitem{najm2009uncertainty}
Habib~N Najm.
\newblock Uncertainty quantification and polynomial chaos techniques in
  computational fluid dynamics.
\newblock {\em Annual review of fluid mechanics}, 41:35--52, 2009.

\bibitem{walters2002uncertainty}
Robert~W Walters and Luc Huyse.
\newblock Uncertainty analysis for fluid mechanics with applications.
\newblock Technical report, NATIONAL AERONAUTICS AND SPACE ADMINISTRATION
  HAMPTON VA LANGLEY RESEARCH CENTER, 2002.

\bibitem{zang2002needs}
Thomas~A Zang.
\newblock {\em Needs and opportunities for uncertainty-based multidisciplinary
  design methods for aerospace vehicles}.
\newblock National Aeronautics and Space Administration, Langley Research
  Center, 2002.

\bibitem{petersson2020discrete}
N~Anders Petersson, Fortino~M Garcia, Austin~E Copeland, Ylva~L Rydin, and
  Jonathan~L DuBois.
\newblock Discrete adjoints for accurate numerical optimization with
  application to quantum control.
\newblock {\em arXiv preprint arXiv:2001.01013}, 2020.

\bibitem{choi2015practical}
Youngsoo Choi, Charbel Farhat, Walter Murray, and Michael Saunders.
\newblock A practical factorization of a schur complement for pde-constrained
  distributed optimal control.
\newblock {\em Journal of Scientific Computing}, 65(2):576--597, 2015.

\bibitem{choi2012simultaneous}
Youngsoo Choi.
\newblock {\em Simultaneous analysis and design in PDE-constrained
  optimization}.
\newblock PhD thesis, Stanford University, 2012.

\bibitem{sirovich1987turbulence}
Lawrence Sirovich.
\newblock Turbulence and the dynamics of coherent structures. i. coherent
  structures.
\newblock {\em Quarterly of applied mathematics}, 45(3):561--571, 1987.

\bibitem{lee2020model}
Kookjin Lee and Kevin~T Carlberg.
\newblock Model reduction of dynamical systems on nonlinear manifolds using
  deep convolutional autoencoders.
\newblock {\em Journal of Computational Physics}, 404:108973, 2020.

\bibitem{kim2020fast}
Youngkyu Kim, Youngsoo Choi, David Widemann, and Tarek Zohdi.
\newblock A fast and accurate physics-informed neural network reduced order
  model with shallow masked autoencoder, 2020.

\end{thebibliography}

\end{document}